\documentclass[preprint]{elsarticle}
\biboptions{longnamesfirst,semicolon}

\usepackage[utf8]{inputenc}
\usepackage[top=2cm, bottom=4.5cm,left=2.5cm, right=2.5cm]{geometry}
\usepackage{amsmath,amsthm,amsfonts,amssymb,amscd}
\usepackage{enumerate}
\usepackage{fancyhdr}
\usepackage{xcolor}
\usepackage{graphicx}
\usepackage{listings}
\usepackage{hyperref}
\usepackage{pgfplots}
\usepackage{tikz}
\usepackage{amsmath}
\usepackage{pgfplotstable}
\usepackage{booktabs}
\usepackage{array}
\usepackage{colortbl}
\usepackage{caption}
\usepackage{subcaption}
\usepackage{verbatim}
\usepgflibrary{shadings}
\graphicspath{{./images/}}

 \pgfplotsset{compat=1.18}
\begin{document}
	\begin{frontmatter} 
        \title{Structure-preserving schemes conserving entropy and kinetic energy} 
	\author[1]{Kunal Bahuguna\corref{cor1}\fnref{fn1}}
	\ead{kunalb@iisc.ac.in}
	\author[1]{Ramesh Kolluru\fnref{fn2}}
	\ead{kollurur@alum.iisc.ac.in}
	\author[1]{S.V. Raghurama Rao} 
	\ead{raghu@iisc.ac.in}
	\cortext[cor1]{Corresponding author}
	\fntext[fn1]{Research Scholar}
	\fntext[fn2]{Present affilication: Senior CFD Consultant, BosonQ PSI (BQP), India Office, 2743, 15th cross, 27th main road, HSR 1st sector, Bangalore -560102} 
	\address[1]{Department of Aerospace Engineering, Indian Institute of Science, Bangalore-560012, India}
		\begin{abstract}
			This paper presents a novel structure-preserving scheme for Euler equations, focusing on the numerical conservation of entropy and kinetic energy.   Explicit flux functions engineered to conserve entropy are introduced within the finite-volume framework. Further, discrete kinetic energy conservation too is introduced.   A systematic inquiry is presented, commencing with an overview of numerical entropy conservation and formulation of entropy-conserving and kinetic energy-preserving fluxes, followed by the study of their properties and efficacy.   A novelty introduced is to associate numerical entropy conservation to the discretization of the energy conservation equation.  Furthermore, an entropy-stable shock-capturing diffusion method and a hybrid approach utilizing the entropy distance to manage smooth regions effectively are also introduced.  The addition of artificial viscosity in apprproiate regions ensures entropy generation sufficient to prevent numerical instabilities. Various test cases, showcasing the efficacy and stability of the proposed methodology, are presented.  
		\end{abstract}
	\end{frontmatter}
	\newpageafter{abstract}
        \section{Introduction}
            The current thrust in algorithm development is directed towards devising suitable numerical methods that adhere to the additional constraints inherent in Euler equations, such as preserving entropy and kinetic energy. Owing to their non-linear hyperbolic nature, Euler equations permit discontinuous solutions like shock waves and contact discontinuities. These solutions lack uniqueness, prompting the pursuit of solutions that adhere to the Euler equations in a weak sense. This is in contrast to that of modelling with the Navier-Stokes equations, where physical viscosity guarantees adherence to the second law of thermodynamics, yielding physically meaningful solutions.  Numerical methods devised for Euler equations must also ensure the conservation or production of numerical entropy, thereby converging towards the correct solution. Furthermore, the stability of numerical methods for non-linear partial differential equations hinges on entropy stability, rendering it a highly desirable property.  
            
            Tadmor \cite{tadmor_1987,tadmor_2003} introduced the mathematical conditions required for numerical entropy conservation and stability for systems of hyperbolic equations, within the framework of the finite-volume method. The resulting conservative fluxes of entropy did not have an explicit form and required an expensive numerical quadrature for calculations \cite{barth_1999}.  Roe \cite{roe_2006} introduced affordable and explicit entropy conservative flux functions for Euler equations.  Ranocha \cite{ranocha_2018} introduced various parameterization of entropy variables in different ways.   All of these fluxes are found to be second-order accurate spatially. LeFloch and Rohde \cite{lefloch_2000} introduced arbitrary higher-order entropy conservative fluxes.  
            
            Ensuring accurate preservation of discrete kinetic energy is another crucial objective for numerical flux functions, particularly in flows dominated by kinetic energy \cite{quirk_1994, pirozzoli_2011}. Jameson \cite{jameson_2008} established conditions on the numerical momentum flux to maintain kinetic energy preservation, and a slightly modified version was presented by Ranocha \cite{ranocha_2020}.  Combining this with entropy conservation, Chandrashekar introduced an entropy-conserving and kinetic energy-preserving flux \cite{chandrashekar_2013}, further addressing kinetic energy stability to dissipate discrete kinetic energy akin to the effect of physical viscosity in Navier-Stokes equations. Entropy conservative schemes can also be obtained through an optimization problem \cite{abgrall_2022}.
            
            In the numerical simulation of viscous flows too, additional diffusion becomes necessary due to the inadequacy of physical viscosity in yielding oscillation-free solutions on coarse grids. However, any artificial diffusion introduced must invariably support entropy generation to prevent nonphysical solutions from entropy violations. Notably, improper entropy generation, as observed in Roe's approximate Riemann solver \cite{roe_1981}, leads to issues such as expansion shocks and shock instabilities such as carbuncle phenomenon \cite{ismail_2009a}. To circumvent such anomalies, various schemes, especially the family of Riemann solvers, need to explicitly ensure the satisfaction of entropy inequalities, thereby averting numerical aberrations. Ismail and Roe \cite{ismail_2009} proposed a stabilizing matrix diffusion method to achieve a `consistent' entropy across shocks.  Nevertheless, while numerical diffusion tailored for shock waves proves effective, it might not be suitable for smooth regions. Consequently, many schemes adopt a hybrid approach in which numerical diffusion is adapted based on flow gradients \cite{harten_1978, jameson_1981}. Furthermore, achieving higher-order entropy stability is feasible through scaled entropy variable reconstruction or sign-preserving ENO reconstruction \cite{fjordholm_2012}.  
            
            In this work, a novel methodology is presented, combining entropy conservation and kinetic energy preservation strategies.  a special feature of Euler equations, that the energy equation actually contains an entropy conservation equation and a kinetic energy conservation equation, is exploited by associating numerical entropy conservative fluxes with the energy equation. Further, entropy stability is ensured by appropriately adding numerical diffusion as required.  The structure of this paper is as follows: In Section 2, a concise introduction is given to numerical entropy conservation.   In Section 3, a diffusion-based strategy for conservation, detailing the formulation of entropy-conserving and kinetic-energy-preserving numerical fluxes is outlined. The numerical investigation of these fluxes, focusing on aspects like numerical entropy conservation, computational efficiency and accurate contact discontinuity capturing, is presented in Section 4, including the demonstration of the experimental order of convergence.  Then, an entropy-stable shock-capturing numerical diffusion method based on Rankine-Hugoniot conditions, highlighting its ability to generate appropriate entropy, is introduced in Section 5. Addressing the insufficiency of this numerical diffusion in smooth regions, a hybrid scheme based on entropy distance is proposed in Section 6..  In Section 7, the numerical results for one-dimensional and two-dimensional benchmark test cases for Euler equations are presented. Finally, the conclusions drawn from the work are presented in Section 8. 
    \section{Structures within Euler Equations: Entropy and Kinetic Energy}   
    Governing equations for compressible inviscid flow are the unsteady Euler equations, with the closure based on the assumption of a perfect gas.  
    Focusing on the one-dimensional Euler equations with conserved variable vector $\mathbf{U}=[\rho,\rho u, \rho E]^T$ and the flux vector $\mathbf{F}(\mathbf{U})=[\rho u,\rho u^2+p,\rho u E + pu]^T$, they can be expressed as conserving mass, momentum and energy, as in \eqref{euler_eq}. 
    \begin{subequations}\label{euler_eq}
        \begin{align} 
            \frac{\partial \mathbf{U}}{\partial t} + \frac{\partial \mathbf{F(\mathbf{U})}}{\partial x} = 0, \ \textrm{or,} \ \textrm{individually} \\  
            \frac{\partial \rho}{\partial t} +\frac{\partial}{\partial x}  \left(\rho u\right)=&0 \label{euler_eq1} \\
            \frac{\partial }{\partial t} (\rho u)+\frac{\partial }{\partial x} (\rho u^2+p)=&0 \label{euler_eq2} \\
            \frac{\partial }{\partial t} (\rho E)+\frac{\partial}{\partial x}(\rho u E+pu)=&0 \label{euler_eq3}
        \end{align}
    \end{subequations}
    Here $\rho$ is the density, $u$ is the velocity, $p$ is the thermodynamic pressure and $E=e+ \frac{1}{2} u^{2}$ is the total specific (internal + kinetic) energy. The system is closed using the ideal gas law $p=(\gamma-1)\rho e$ where $e$ is the specific internal energy, and $\gamma$ is the ratio of specific heats of the gas. This system of equations is hyperbolic in nature, with the flux Jacobian matrix $\tilde{\mathbf{A}}=\partial \bf{F}(\bf{U})/\partial \bf{U}$ having real and distinct eigenvalues. Since the system is non-linear, Euler equations may generate discontinuous solutions even when the initial conditions are smooth. These solutions are also not unique, and it is important for the numerical schemes to converge to the right physically admissible solutions ({\em i.e.}, the vanishing viscosity solutions).  
    \subsection{Entropy inequality, entropy variables, symmetric form and entropy potential functions}  
    Since the system of Euler equations is hyperbolic, it admits an additional conservation law \cite{friedrichs_1971}, representing the conservation of entropy in smooth regions.   Consider a convex entropy function $\eta(\mathbf{U})$.  Multiplying the Euler equations by its derivative, we obtain 
    \begin{equation}
        \frac{\partial \eta(\mathbf{U})}{\partial \mathbf{U}} \left[ \frac{\partial \mathbf{U}}{\partial t} + \frac{\partial \mathbf{F(\mathbf{U})}}{\partial x} = 0 \right] 
    \end{equation} 
    or 
    \begin{equation*}
        \frac{\partial \eta(\mathbf{U})}{\partial t} + \mathbf{\eta^{'}} (\mathbf{U}) \mathbf{F^{'}} (\mathbf{U}) \frac{\partial \mathbf{U}}{\partial x} = 0  
    \end{equation*} 
    If we now introduce a consistency condition as 
    \begin{equation}\label{ent_comp_eq}
		\mathbf{\zeta^{'}} (\mathbf{U})=\mathbf{\eta^{'}} (\mathbf{U})^T \mathbf{F^{'}} (\mathbf{U})
	\end{equation} 
    we obtain 
    \begin{equation*}
        \frac{\partial \eta(\mathbf{U})}{\partial t} + \mathbf{\zeta^{'}} (\mathbf{U}) \frac{\partial \mathbf{U}}{\partial x} = 0  
    \end{equation*} 
    or 
    \begin{equation}
        \frac{\partial \eta(\mathbf{U})}{\partial t} + \frac{\partial \zeta(\mathbf{U}) } {\partial x} = 0  
    \end{equation} 
    which is the entropy conservation equation, with $\zeta(\mathbf{U})$ being the entropy flux function.   
    The above entropy conservation equation is valid only in smooth regions of the flow.  If we include the regions of discontinuities (shock waves), then we have the inequality 
        \begin{equation} \label{entropy_inequality} 
		\frac{\partial \eta (\mathbf{U})}{\partial t}+\frac{\partial \zeta (\mathbf{U})}{\partial x} \leq 0 
	\end{equation}
    For 1-D Euler equations, the entropy pair fulfilling \eqref{ent_comp_eq} is given by 
	\begin{equation*}\label{ent_pair}
		\eta (\mathbf{U})=-\frac{\rho s}{\gamma-1}, \qquad \zeta (\mathbf{U})=-\frac{\rho u s}{\gamma-1},
	\end{equation*}
	with $s=ln(p/\rho^\gamma)$, representing the specific entropy. Interestingly, this entropy pair is valid not only for Euler equations but also for Navier-Stokes equations.  
    Mock \cite{mock_1980} showed that the Euler equations with an entropy equation implies the existence of entropy variable vector $\bf{V}$ with a one-to-one mapping to $\bf{U}$, and is given by 
    \begin{equation*}\label{ent_var}
        \mathbf{V}=\frac{d \eta (\mathbf{U})}{d \mathbf{U}}=\left[\frac{\gamma-s}{\gamma-1}-\frac{\rho u^2}{2p},\frac{\rho u}{p},-\frac{\rho}{p} \right]^T
    \end{equation*} 
    Entropy variable symmetrizes the hyperbolic system\eqref{euler_eq}, {\em i.e.}, changing the variables from $\mathbf{U}$ to $\mathbf{V}$ in the Euler equations gives the following symmetric system
    \begin{equation}
        \tilde{\mathbf{A}} \frac{\partial \mathbf{V}}{\partial t} +\tilde{\mathbf{B}} \frac{\partial \mathbf{V}}{\partial x} =0
    \end{equation}
    where $\tilde{\mathbf{A}}$ is a symmetric positive definite matrix and $\tilde{\mathbf{B}}$ is a symmetric matrix. Further, entropy variable facilitates the introduction of {\em entropy potential function} $\phi$ and {\em entropy flux potential function} $\psi$ (together called as entropy pair) as 
    \begin{equation}\label{pot_fun}
	\phi=\mathbf{V}^T \cdot \mathbf{U} - \eta(\mathbf{U}) = \rho \qquad \text{and} \qquad \psi=\mathbf{V}^T \cdot \mathbf{F}(\mathbf{U}) - \zeta(\mathbf{U})=\rho u
    \end{equation}
    This formulation of entropy inequality together with the hyperbolic system, and thus the concept of entropic solutions, extends to arbitrary systems of conservation laws in $n$ dimensions, as elucidated by Harten \cite{harten_1983}. 
    \subsection{Structures in the total energy equation: entropy conservation and kinetic energy preservation}  
    The two specific structures within the Euler equations worth preserving in a numerical scheme are: {\em entropy conservation} and {\em kinetic energy preservation}.  Both of them are parts of the total energy conservation equation and this fact is exploited in our work.  

     Let us first obtain the kinetic energy equation separately; we multiply \eqref{euler_eq2} by $u$, \eqref{euler_eq1} by $-u^2/2$, and then sum them to obtain:
    \begin{equation}\label{ke_eq}
        \frac{\partial}{\partial t} \left(\frac{\rho u^2}{2} \right) + \frac{\partial}{\partial x} \left(\frac{\rho u^3}{2} \right) +u\frac{\partial p}{\partial x}=0
    \end{equation}
    The kinetic energy equation is thus implicit in the mass and momentum conservation equations, and the numerical discretizations of mass and momentum may not necessarily enforce kinetic energy preserving.  Thus, developing {\em kinetic energy preserving schemes} is a separate endeavor. 
  
    Now consider the energy conservation equation  \eqref{euler_eq3}. 
    $$ \frac{\partial }{\partial t} (\rho E)+\frac{\partial}{\partial x}(\rho u E+pu)=0 $$  
    Using the relation $E=e+\frac{u^2}{2}$ in the energy equation, the internal energy and kinetic energy parts can be separated. 
	\begin{equation*}
		\frac{\partial}{\partial t} \left(\rho e \right) + \frac{\partial}{\partial x} \left(\rho u e \right) +p\frac{\partial u}{\partial x} + \frac{\partial}{\partial t} \left(\frac{\rho u^2}{2} \right) + \frac{\partial}{\partial x} \left(\frac{\rho u^3}{2} \right) +u\frac{\partial p}{\partial x}=0
	\end{equation*} 
    or 
        \begin{equation}
	\underbrace{\frac{\partial}{\partial t} \left(\rho e \right) + \frac{\partial}{\partial x} \left(\rho u e \right) + \rho e (\gamma -1) \frac{\partial u}{\partial x}}_{\rm internal \ energy \ consevation \ + \ source \ term}  +  
        \underbrace{\frac{\partial}{\partial t} \left(\frac{\rho u^2}{2} \right) + \frac{\partial}{\partial x} \left(\frac{\rho u^3}{2} \right) +u\frac{\partial p}{\partial x}=0}_{\rm kinetic \ energy \ conservation \ + \ source \ term}  
	\end{equation} 
	Now, subtracting \eqref{ke_eq} from the above equation, one can focus only on the internal energy equation.   
	\begin{equation}
		\frac{\partial}{\partial t} \left(\rho e \right) + \frac{\partial}{\partial x} \left(\rho u e \right) 
        + p\frac{\partial u}{\partial x}=0
	\end{equation}
	Utilizing $e=\frac{p}{\rho(\gamma-1)}$, we obtain 
        \begin{equation*} 
        \frac{\partial p}{\partial t} + \frac{\partial (up)}{\partial x} + (\gamma -1) p \frac{\partial u}{\partial x} = 0 
        \end{equation*} 
        or 
        \begin{equation*}
		\frac{\partial p}{\partial t}  + u \frac{\partial p}{\partial x}  +\gamma p\frac{\partial u}{\partial x}=0
	\end{equation*}
	Using $ds=\frac{dp}{p}- \frac{\gamma}{\rho} d \rho$, we obtain from the above equation 
        \begin{equation*} 
            p \frac{\partial s}{\partial t} + \frac{\gamma p}{\rho} \frac{\partial \rho}{\partial t} + u p \frac{\partial s}{\partial x} 
            +  \frac{u \gamma p}{\rho} \frac{\partial \rho}{\partial x} + \gamma p \frac{\partial u}{\partial x} = 0 
        \end{equation*}  
        or 
        \begin{equation} 
            \frac{\gamma p}{\rho}  \left[ \frac{\partial \rho}{\partial t} + u \frac{\partial \rho}{\partial x} 
            + \rho \frac{\partial u}{\partial x} \right]  
            + p \left[ \frac{\partial s}{\partial t} + u \frac{\partial s}{\partial x} \right] = 0 
        \end{equation} 
        or (after utilizing the continuity equation) 
        \begin{equation} 
            \frac{\partial s} {\partial t} + u \frac{\partial s}{\partial x} = 0 
        \end{equation}   
        Thus, the entropy convection equation is a part of the internal energy conservation equation.  Together with the mass conservation equation, we can obtain the entropy conservation equation as 
	\begin{equation}\label{ent_eq2}
		\frac{\partial}{\partial t} \left(\rho s\right)+\frac{\partial}{\partial x} \left(\rho u s\right)=0.
	\end{equation}
	Thus, the entropy conservation equation is implicit in the internal energy conservation equation. Again, like the kinetic energy preservation, numerical discretizations of Euler equations do not necessarily satisfy the discrete entropy conservation.  Numerical schemes satisfying these additional conservation equations are often termed as \it{structure preserving schemes} \normalfont and are discussed in the next subsections. The fact that the total energy equation contains both the entropy conservation and the kinetic energy conservation is used in section 2 to construct entropy conservative flux by only modifying the energy flux. Note that the equations \eqref{ke_eq} and \eqref{ent_eq2} are only valid for smooth regions and not valid across discontinuities like shocks.  
  \subsection{Numerical Entropy Conservation}


In semi-discrete form, the update formula for a cell  (with cell-centre at $x_{j}$ and cell-interfaces $x_{j\pm \frac{1}{2}}$) can be expressed as 
\begin{equation}
    \frac{d\mathbf{U}_j}{dt}=-\frac{1}{\Delta x_j}\left( \mathbf{F}_{j+\frac{1}{2}}-\mathbf{F}_{j-\frac{1}{2}} \right)
\end{equation}
When multiplied by the entropy variables at cell $j$, denoted as $\mathbf{V}_j$, it yields the semi-discrete entropy conservation equation:
\begin{equation}
    \frac{d\eta(\mathbf{U})_j}{dt}=-\frac{1}{\Delta x_j} \mathbf{V}_j \cdot\left( \mathbf{F}_{j+\frac{1}{2}}-\mathbf{F}_{j-\frac{1}{2}} \right)
\end{equation}
The right-hand side of the above equation is not in conservation form.  Therefore, we modify it as  
\begin{equation}
    \Delta x_j \frac{d\eta(\mathbf{U})_j}{dt}=-\left(\overline{\mathbf{V}}_{j+\frac{1}{2}}-\frac{1}{2}\Delta \mathbf{V}_{j+\frac{1}{2}}\right) \cdot \mathbf{F}_{j+\frac{1}{2}}+\left(\overline{\mathbf{V}}_{j-\frac{1}{2}}+\frac{1}{2}\Delta \mathbf{V}_{j-\frac{1}{2}}\right) \cdot \mathbf{F}_{j-\frac{1}{2}}
\end{equation}
by introducing $\overline{\mathbf{V}}_{j+\frac{1}{2}}=(\mathbf{V}_{j+1}+\mathbf{V}_j)/2$ and $\Delta{\mathbf{V}}_{j+\frac{1}{2}}=\mathbf{V}_{j+1}-\mathbf{V}_j$. The numerical entropy flux consistent with a given potential function $\psi$ at an interface is defined as $\zeta_{j+\frac{1}{2}}=\overline{\mathbf{V}}_{j+\frac{1}{2}} \cdot \mathbf{F}_{j+\frac{1}{2}} - \overline{\psi}_{j+\frac{1}{2}}$, where $\overline{\psi}_{j+\frac{1}{2}}=(\psi_{j+1}+\psi_{j})/2$. Thus, the above equation simplifies to 
\begin{equation}
    \Delta x_j\frac{d\eta(\mathbf{U})_j}{dt}+\zeta_{j+\frac{1}{2}}-\zeta_{j-\frac{1}{2}}=\frac{1}{2} \Delta \mathbf{V}_{j+\frac{1}{2}} \cdot \mathbf{F}_{j+\frac{1}{2}} - \frac{1}{2}\Delta \psi_{j+\frac{1}{2}}+ \frac{1}{2} \Delta \mathbf{V}_{j-\frac{1}{2}} \cdot \mathbf{F}_{j-\frac{1}{2}} -\frac{1}{2}\Delta\psi_{j-\frac{1}{2}}
\end{equation}
To ensure the numerical conservation of entropy, the right hand side must vanish.  Therefore, the entropy conservative numerical fluxes $\mathbf{F}^c$ should satisfy the condition \cite{tadmor_1987,tadmor_2003}: 
\begin{equation}\label{ec_cond}
    (\mathbf{V}_{j+1}-\mathbf{V}_{j}) \cdot \mathbf{F}^c_{j+\frac{1}{2}}=\psi_{j+1}-\psi_{j}
\end{equation}
 
 \subsection{Numerical Kinetic Energy Preservation}
	 A semi-discrete form of the kinetic energy equation can be obtained by multiplying the semi-discrete mass and momentum equations by $-u_j^2/2$ and $u_j$, respectively, and adding them.
	\begin{equation}\label{sdke_eq}
		\Delta x_j \frac{d}{dt} \left( \rho_j \frac{u_j^2}{2} \right) = \frac{u_j^2}{2}\left( (\rho u)_{j+\frac{1}{2}} - (\rho u)_{j-\frac{1}{2}} \right) -u_j \left( (\rho u^2)_{j+\frac{1}{2}} - (\rho u^2)_{j-\frac{1}{2}} \right) - u_j \left( p_{j+\frac{1}{2}}-p_{j-\frac{1}{2}} \right)
	\end{equation}
	Jameson \cite{jameson_2008} introduced the following form of convective momentum.
	\begin{equation}\label{kep_cond}
		(\rho u^2)_{j+\frac{1}{2}}=(\rho u)_{j+\frac{1}{2}}(u_{j+1}+u_j)/2
	\end{equation}
	We get the local kinetic energy preservation equation using this flux in \eqref{sdke_eq}.
	\begin{equation}\label{local_kep_eq}
		\Delta x_j \frac{d}{dt}\left(\frac{\rho_j u_j^2}{2}\right)+\left(\frac{\rho u^3}{2}\right)_{j+\frac{1}{2}}-\left(\frac{\rho u^3}{2}\right)_{j-\frac{1}{2}}=-u_j(p_{j+\frac{1}{2}}-p_{j-\frac{1}{2}})
	\end{equation}
	Here the kinetic energy flux is defined as $\left(\frac{\rho u^3}{2}\right)_{j+\frac{1}{2}}=\left(\frac{u_{j+1}u_j}{2}\right)(\rho u)_{j+\frac{1}{2}}$ which is of the form $F(U_L,U_R)=-F(U_R,U_L)$ (given mass flux satisfies the same property). Summing equation (\ref{sdke_eq}) on all $j=1$ to $n$, using the form of convective momentum flux as in the above equation, and taking the boundary mass fluxes as $(\rho u)_{1/2}=\rho_L u_L$ and $(\rho u)_{n+1/2}=\rho_R u_R$ and momentum fluxes as $(\rho u^2)_{1/2}=\rho_L u^2_L+p_L$ and $(\rho u^2)_{n+1/2}=\rho_R u^2_R+p_R$, we recover the global kinetic energy preservation equation.
	\begin{equation}
		\sum_{j=1}^{n}\Delta x_j \frac{d}{dt}\left(\frac{\rho_j u_j^2}{2}\right)+u_R\left(\frac{\rho_R u_R^2}{2}+p_R\right)-u_L\left(\frac{\rho_L u_L^2}{2}+p_L\right)=\sum_{j=1}^n p_{j+\frac{1}{2}}(u_{j+1}-u_j)
	\end{equation}
	Note that we are still to choose the forms of mass flux $(\rho u)_{j+\frac{1}{2}}$, pressure flux $p_{j+\frac{1}{2}}$ and the energy flux $(\rho u E+pu)_{j+\frac{1}{2}}$.  Gassner et al. \cite{gassner_2016} observed that the arithmetic pressure average preserves kinetic energy better than the schemes using different pressure averages. Ranocha \cite{ranocha_2020} showed that an arithmetic pressure average is well suited as pressure flux for Euler equations at the incompressible limit. This flux formulation still leaves mass flux $(\rho u)_{j+\frac{1}{2}}$, and energy flux $(\rho u E + pu)_{j+\frac{1}{2}}$ undefined and allows us to fix these using other criteria.
	\section{Entropy Conservative and Kinetic Energy Preserving Flux}
	\subsection{Entropy Conservative Numerical Flux}
	Entropy conservative fluxes can be constructed by assuming $\mathbf{F}^c_{j+\frac{1}{2}}$ to be of the following form (average flux + diffusive flux) as 
	\begin{equation}\label{ec_flux_gen}
		\mathbf{F}^c_{j+\frac{1}{2}}=\frac{1}{2}(\mathbf{F}_{j+1}+\mathbf{F}_{j})-\frac{1}{2}\mathbf{\alpha} (\mathbf{V}_{j+1}-\mathbf{V}_j)=\overline{\mathbf{F}}_{j+1/2}-\frac{1}{2}\mathbf{\alpha} (\mathbf{V}_{j+1}-\mathbf{V}_j) = \overline{\mathbf{F}}_{j+1/2}-\frac{1}{2}\mathbf{\alpha} \Delta \mathbf{V}
	\end{equation}   
	Here, $\alpha$ is a {\em scalar diffusion coefficient}, adding diffusion to satisfy semi-discrete entropy conservation. It can be obtained by substituting the flux \eqref{ec_flux_gen} into \eqref{ec_cond}. Thus, 
        $\Delta \mathbf{V} \cdot \mathbf{F}^c_{j+\frac{1}{2}} = \Delta \psi $ gives  
    \begin{equation}
        \Delta \mathbf{V} \cdot \left( \overline{\mathbf{F}}_{j+1/2}-\frac{1}{2}\mathbf{\alpha} \Delta \mathbf{V} \right) = \Delta \psi 
    \end{equation}
    This gives one form of diffusion coefficient, termed as $\alpha_{1}$ here, as
	\begin{equation}\label{alpha_1}
		\frac{\alpha_1}{2}=\frac{\overline{\mathbf{F}} \cdot \Delta \mathbf{V}-\Delta \psi}{\Delta \mathbf{V} \cdot \Delta \mathbf{V}}=\frac{\overline{F^1}\Delta V_1+\overline{F^2}\Delta V_2+\overline{F^3}\Delta V_3-\Delta \psi}{\Delta V_1 \Delta V_1+\Delta V_2 \Delta V_2+\Delta V_3 \Delta V_3}
	\end{equation}
	where $\overline{F^1}$, $\overline{F^2}$ and $\overline{F^3}$ are the average mass, momentum and energy fluxes of $j^{th}$ and $(j+1)^{th}$ cells. $\Delta(X)$ of any quantity $X$ is defined as $\Delta(X) = X_{j+1}-X_{j}$. This flux, termed EC1, can be defined by
	\begin{equation}\label{EC1_flux}
		\begin{aligned}
		F^1_{EC1}=&\overline{F^1}-\frac{1}{2}\alpha_{1} (\Delta V_1)\\
		F^2_{EC1}=&\overline{F^2}-\frac{1}{2}\alpha_{1} (\Delta V_2)\\
		F^3_{EC1}=&\overline{F^3}-\frac{1}{2}\alpha_{1} (\Delta V_3)
		\end{aligned}
	\end{equation}
	It can be noted that $\alpha_{1}$ is a scalar diffusion that acts on mass, momentum, and energy equations, resulting in an entropy conservative flux. Note that this flux can also be obtained through an optimization procedure (see appendix A) as given by Abgrall \cite{abgrall_2018}. This flux, however, is not kinetic energy preserving. It was shown in section 2.2 that the entropy equation is a part of the energy equation, and just adding a corresponding diffusion to the energy flux can ensure entropy conservation. This can be achieved by taking $\alpha_{2}=[0,0,\alpha_{2}]^T$ in equation \eqref{ec_flux_gen} and we get the following form of $\alpha_2$
	\begin{equation} \label{alpha_2} 
		\frac{\alpha_2}{2}=\frac{\overline{\mathbf{F}} \cdot \Delta \mathbf{V}-\Delta \psi}{\Delta V_3 \Delta V_3}=\frac{\overline{F^1}\Delta V_1+\overline{F^2}\Delta V_2+\overline{F^3}\Delta V_3-\Delta \psi}{\Delta V_3 \Delta V_3}
	\end{equation}
	This flux, termed EC2, can be defined by
	\begin{equation}\label{EC2_flux}
		\begin{aligned}
			F^1_{EC2}=&\overline{F^1}\\
			F^2_{EC2}=&\overline{F^2}\\
			F^3_{EC2}=&\overline{F^3}-\frac{1}{2}\alpha_{2} (\Delta V_3)
		\end{aligned}
	\end{equation}\\
	The entropy conservation of a numerical flux thus can theoretically be achieved by adding diffusion only to the energy equation. This makes mass and momentum flux available for other modifications, as shown in the following subsection.
	\subsection{Entropy Conservative and Kinetic Energy Preserving Numerical Flux}
	Jameson \cite{jameson_2008} gives the following numerical condition on mass and momentum fluxes for the discrete kinetic energy preservation condition to be satisfied.
	\begin{equation}\label{ke_cond}
		F^2_{j+\frac{1}{2}}=u_{j+\frac{1}{2}} F^1_{j+\frac{1}{2}} + \overline{p}
	\end{equation}
	Here $u_{j+\frac{1}{2}}=\frac{1}{2}(u_{j+1}+u_j)$ and $\overline{p}$ is any consistent pressure flux. An entropy conservative and kinetic energy preserving flux can be obtained by taking the correct form of momentum flux given by \eqref{ke_cond} in equation \eqref{ec_flux_gen}. Like EC2 flux, the diffusion term $\alpha_3$ is only applied to the energy equation. Thus, the entropy conservative and kinetic energy preserving numerical flux (termed ECKEP here) is given as
	\begin{equation}\label{ECKEP_flux}
		\begin{aligned}
			F^1_{ECKEP}=&\overline{F^1}\\
			F^2_{ECKEP}=&\overline{F^1} \overline{u}+\overline{p}\\
			F^3_{ECKEP}=&\overline{F^3}-\frac{1}{2}\alpha_{3} (\Delta V_3)
		\end{aligned}
	\end{equation}
	where $\alpha_3$ is 
	\begin{equation} \label{alpha_3}
		\frac{\alpha_3}{2}=\frac{\overline{F^1}\Delta V_1+(\overline{F^1}\overline{u}+\overline{p})\Delta V_2+\overline{F^3}\Delta V_3-\Delta \psi}{\Delta V_3 \Delta V_3}
	\end{equation}
	The convective momentum flux in \eqref{ECKEP_flux} is taken such that it results in the preservation of kinetic energy, and the pressure flux is taken as an arithmetic average.
    
	Note that to prevent denominators of diffusion terms from becoming zero in fluxes (\ref{alpha_1}), (\ref{alpha_2}) and (\ref{alpha_3}), a small parameter $\delta$ is added to the denominators ($\delta = 10^{-16}$ in our computations).
	
        \section{Numerical investigation of the entropy conservative and kinetic energy preserving fluxes}
        
	\subsection{Exact Capturing of Steady Contact Discontinuity} 
    
\begin{figure}[!h]
	\begin{minipage}{.95\textwidth}
		\centering
		\includegraphics[scale=0.4]{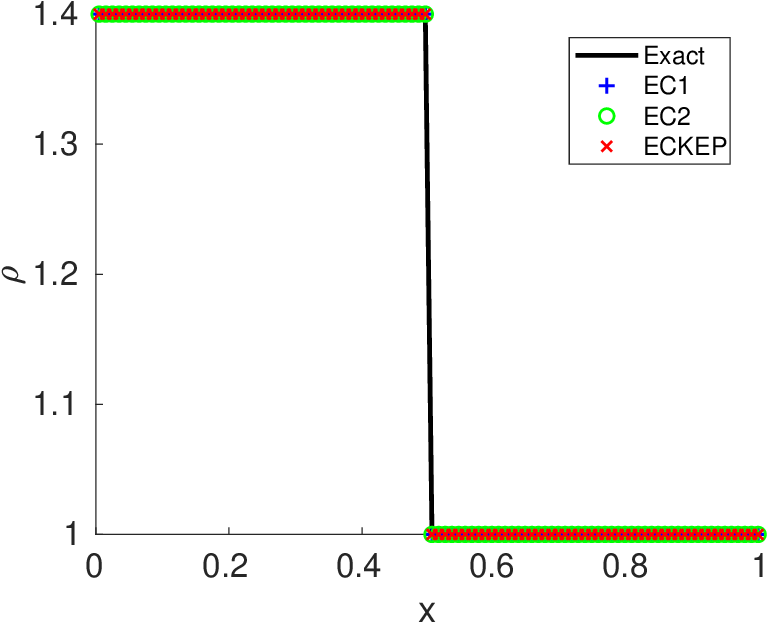}\quad
		\includegraphics[scale=0.4]{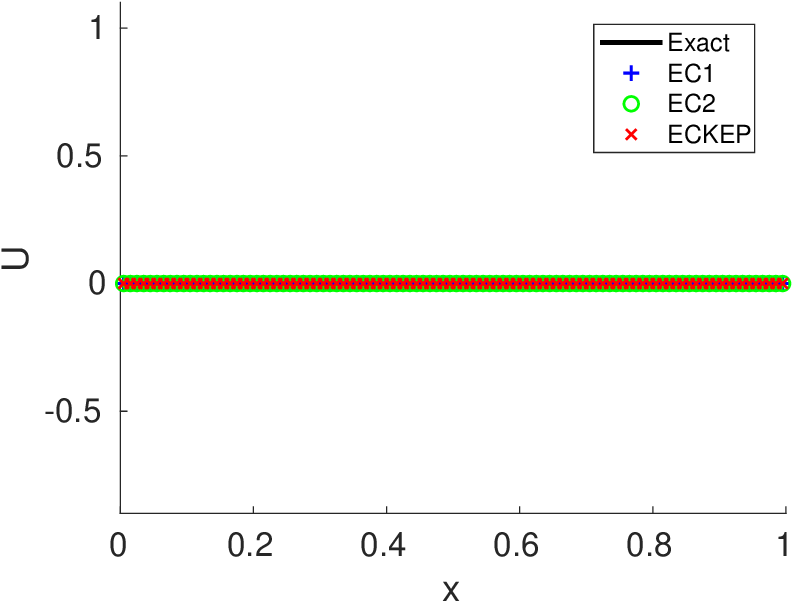}\quad
	\end{minipage}\\[1em]
	\caption{Density and velocity plots for a stationary contact discontinuity at T=2 s by the three entropy conservative fluxes}  
	\label{contact_fig}
	\end{figure}

 	Across a stationary  contact discontinuity in one dimension, we have primitive variable vectors $\mathbf{W}$ on the left and right sides given by
	\begin{equation*}
		\mathbf{W}_L=\left [\begin{matrix} \rho_L \\ 0 \\ p \end{matrix} \right] \qquad and \qquad  \mathbf{W}_R=\left [\begin{matrix} \rho_R \\ 0 \\ p \end{matrix} \right] 
	\end{equation*}
	Resolving contact disctontinuities accurately is important for good resolution of shear and boundary layers in multidimensional flow \cite{leer_1990}. We can easily check whether a numerical scheme will capture a steady contact discontinuity exactly. Flux at the interface across a steady contact disctontinuity is given by
	\begin{equation*}
		\mathbf{F}_{j+1/2}=\left [\begin{matrix} \rho u \\ \rho u^2 +p \\ \rho u E +pu \end{matrix} \right]_{j+\frac{1}{2}}=\left [\begin{matrix} 0 \\ p \\ 0 \end{matrix} \right]
	\end{equation*}
	The above condition is satisfied by EC1, EC2, and ECKEP fluxes.  A Reimann problem with $\mathbf{W}_L=[1.4,0,1]^T$ and $\mathbf{W}_R=[1,0,1]^T$ is considered and the solution is shown in figure \ref{contact_fig}.   The steady contact discontinuities are preserved exactly for all the three numerical fluxes presented.  
	\subsection{Experimental order of convergence (EOC)}
	Entropy conservative fluxes \eqref{EC1_flux} and \eqref{EC2_flux} and the entropy conservative and kinetic energy preserving flux \eqref{ECKEP_flux} are all second-order accurate in space. This can be demonstrated by taking a sinusoidal test case with the following initial conditions.  
	\begin{subequations}\label{den_wave}
		\begin{align}
			\rho (x,0)=&1+0.2\sin({2\pi x}) \\
			u(x,0)=&0.1 \\
			p(x,0)=&1
		\end{align}
	\end{subequations}
	EOC was computed on a domain of $[0,1]$. A periodic boundary was applied at both ends. The exact solution to the above problem is known. The velocity and pressure remain constant and the density wave advances with time given by the following equation.
	
	\begin{equation}
		\rho(x,t)=1+0.2 \sin\left( 2\pi(x-0.1t) \right)
	\end{equation}
	The computational domain was divided into N=40,80,160,320,640 and 1280 points, and $L_1$ and $L_2$ errors in density were computed at $t=10$. The SSPRK-3 \cite{shu_1988}  method was used for time integration. For a $k^{th}$ order accurate scheme, the error on N points $(dx=1/N)$ and 2N points $(dx=1/2N)$ are given by
	\begin{subequations}
		\begin{align}
		\|\epsilon_{N}\|=&C(\Delta x)^k+O (\Delta x)^{k+1} \\
		\|\epsilon_{2N}\|=&C\left(\frac{\Delta x}{2}\right)^k+O \left(\frac{\Delta x}{2}\right)^{k+1}
		\end{align}
	\end{subequations}
	Thus, the order of convergence $k$ can be calculated using
	\begin{equation}
		k=\log_2\left( \frac{ \| \epsilon_{N} \|}{\| \epsilon_{2N} \|} \right)
	\end{equation}
	$L_1$ and $L_2$ errors with computed order of accuracy for schemes EC1, EC2 and ECKEP are given in tables \ref{table:tab_ec1}, \ref{tab_ec2} and \ref{tab_eckep} respectively. As shown, all fluxes were found to be of the order $O(\Delta x)^2$.  \\
	\pgfplotstableset{
		columns/N/.style={
			fixed,precision=1,
			column name={N},
		},
		columns/l1err/.style={
			fixed,precision=8,
			column name={$L_1 \text{ Error}$},
			dec sep align,
		},
		columns/EOC1/.style={
			fixed, fixed zerofill,precision=8,
			string replace={0}{}, 
			column name={EOC},
			dec sep align,
		},
		columns/l2err/.style={
			fixed, fixed zerofill,precision=8,
			column name=$L_2 \text{ Error}$,
			dec sep align,
		},
		columns/EOC2/.style={
			fixed, fixed zerofill,precision=8,
			string replace={0}{}, 
			column name={EOC},
			dec sep align={c|},{column type/.add={}{|}}
		},
		empty cells with={--} 
	}
	\begin{table}[!h]
		\centering
		\begin{filecontents*}{EOC_results_EC1_m.txt}
N,l1err,EOC1,l2err,EOC2
40,0.00328453059071451,0,0.00364680752844891,0
80,0.000821765000727281,1.99888941024011,0.000912587278308923,1.99859961147585
160,0.00020545818582347,1.99988107568862,0.000228201215655253,1.99965606534347
320,5.13652071665826e-05,1.99998144249278,5.70536825417229e-05,1.99991456337534
640,1.28413366959398e-05,1.99999607857642,1.42636314876086e-05,1.99997867319833
1280,3.21033590379724e-06,1.99999922263826,3.56592103356624e-06,1.99999467506364
\end{filecontents*}
\pgfplotstabletypeset[
		col sep=comma,
		string type,
		every head row/.style={before row=\hline,after row=\hline},
		every last row/.style={after row=\hline},
		every first column/.style={column type/.add={|}{}},
		every last column/.style={column type/.add={}{|}}
		]{EOC_results_EC1_m.txt}
		\caption{EOC using $L_1$ and $L_2$ errors for EC1 flux}\label{table:tab_ec1}
	\end{table}
	
	\begin{table}[!h]
		\centering
		\begin{filecontents*}{EOC_results_EC2_m.txt}
N,l1err,EOC1,l2err,EOC2
40,0.00325341579408761,0,0.00361227494539429,0
80,0.000813920332634532,1.99899571990661,0.000904047729465203,1.99843685929566
160,0.000203508590120433,1.9997978968203,0.000226070607762138,1.99962550775331
320,5.08787825246704e-05,1.99995363811139,5.65212164507168e-05,1.99990901357669
640,1.27197922014311e-05,1.99998904683607,1.41305244482758e-05,1.99997750405145
1280,3.17995415987041e-06,1.99999722820783,3.53264484651407e-06,1.99999439098601
\end{filecontents*}
\pgfplotstabletypeset[
		col sep=comma,
		string type,
		every head row/.style={before row=\hline,after row=\hline},
		every last row/.style={after row=\hline},
		every first column/.style={column type/.add={|}{}},
		every last column/.style={column type/.add={}{|}}
		]{EOC_results_EC2_m.txt}%
		\caption{EOC using $L_1$ and $L_2$ errors for EC2 flux}\label{tab_ec2}
	\end{table}
	\begin{table}[!h]
		\centering
		\begin{filecontents*}{EOC_results_ECKEP_m.txt}
N,l1err,EOC1,l2err,EOC2
40,0.00325341575819553,0,0.00361227493977574,0
80,0.000813920329474811,1.99899570959131,0.000904047729090038,1.99843685765038
160,0.00020350859014594,1.99979789103879,0.000226070607691269,1.99962550760687
320,5.08787830914081e-05,1.99995362222207,5.65212164325322e-05,1.99990901358859
640,1.27197928894372e-05,1.99998898487169,1.41305244697578e-05,1.99997750139403
1280,3.17995334022537e-06,1.99999767810243,3.53264480010014e-06,1.99999441213423
\end{filecontents*}
\pgfplotstabletypeset[
		col sep=comma,
		string type,
		every head row/.style={before row=\hline,after row=\hline},
		every last row/.style={after row=\hline},
		every first column/.style={column type/.add={|}{}},
		every last column/.style={column type/.add={}{|}}
		]{EOC_results_ECKEP_m.txt}%
		\caption{EOC using $L_1$ and $L_2$ errors for ECKEP flux}\label{tab_eckep}
	\end{table}
	\begin{figure}[!h] \label{EOC}
		\begin{minipage}{.95\textwidth}
			\centering
			\includegraphics[scale=0.4]{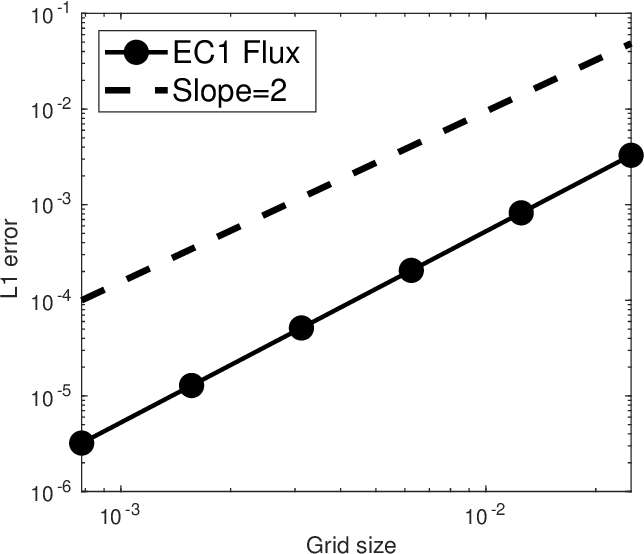}\quad
			\includegraphics[scale=0.4]{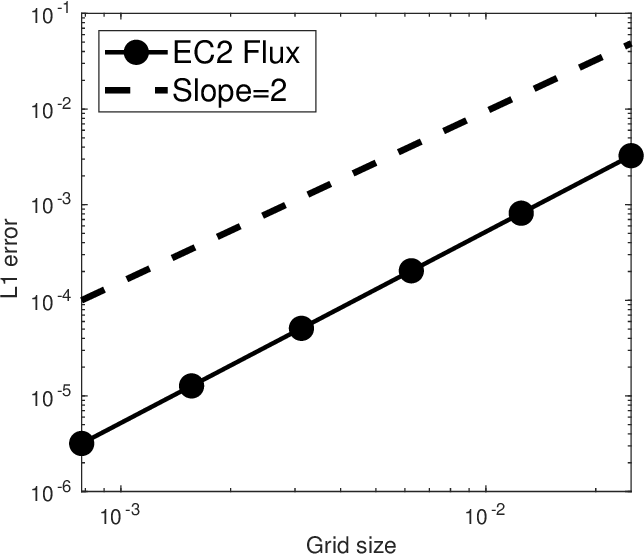}\quad
			\includegraphics[scale=0.4]{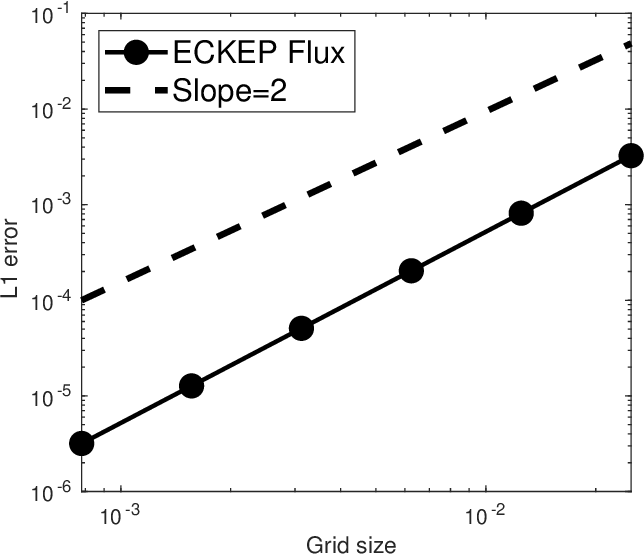}\\
			\includegraphics[scale=0.4]{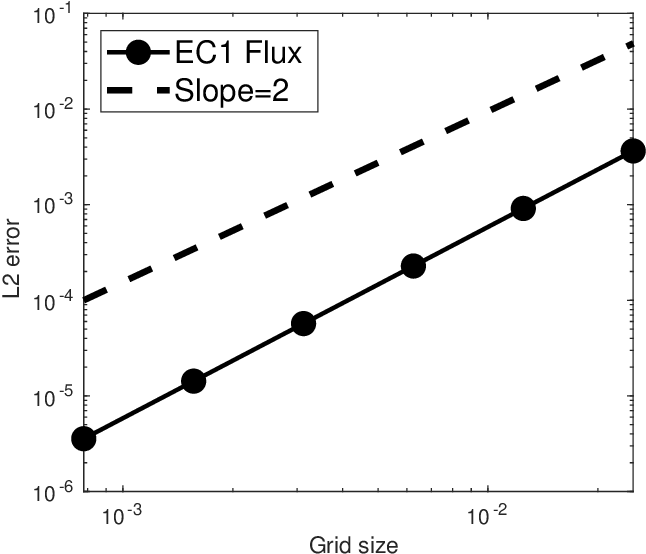}\quad
			\includegraphics[scale=0.4]{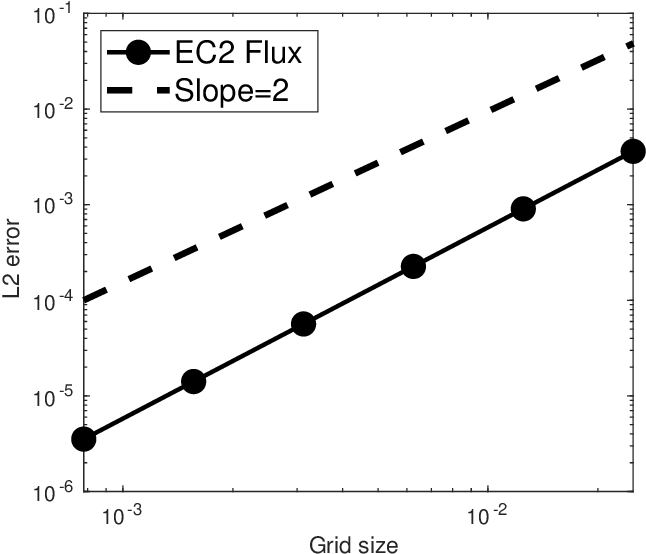}\quad
			\includegraphics[scale=0.4]{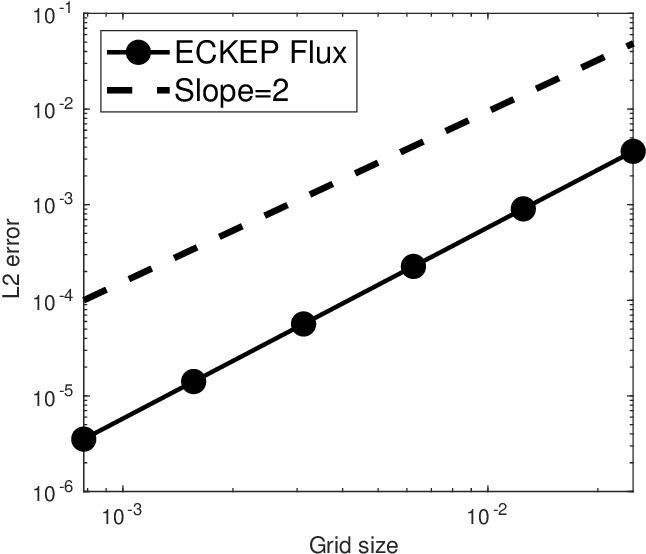}\\
		\end{minipage}\\[1em]
		\caption{Log-log plot of errors with grid size for EC1, EC2 and ECKEP fluxes respectively} 
	\end{figure}
	\newpage
	
        \subsection{Conservation properties of the new flux functions}
	\subsubsection{Taylor-Green Vortex}
\begin{figure}[!h]
		\centering
		\begin{minipage}{.45\textwidth}
			\centering
			\includegraphics[scale=0.15]{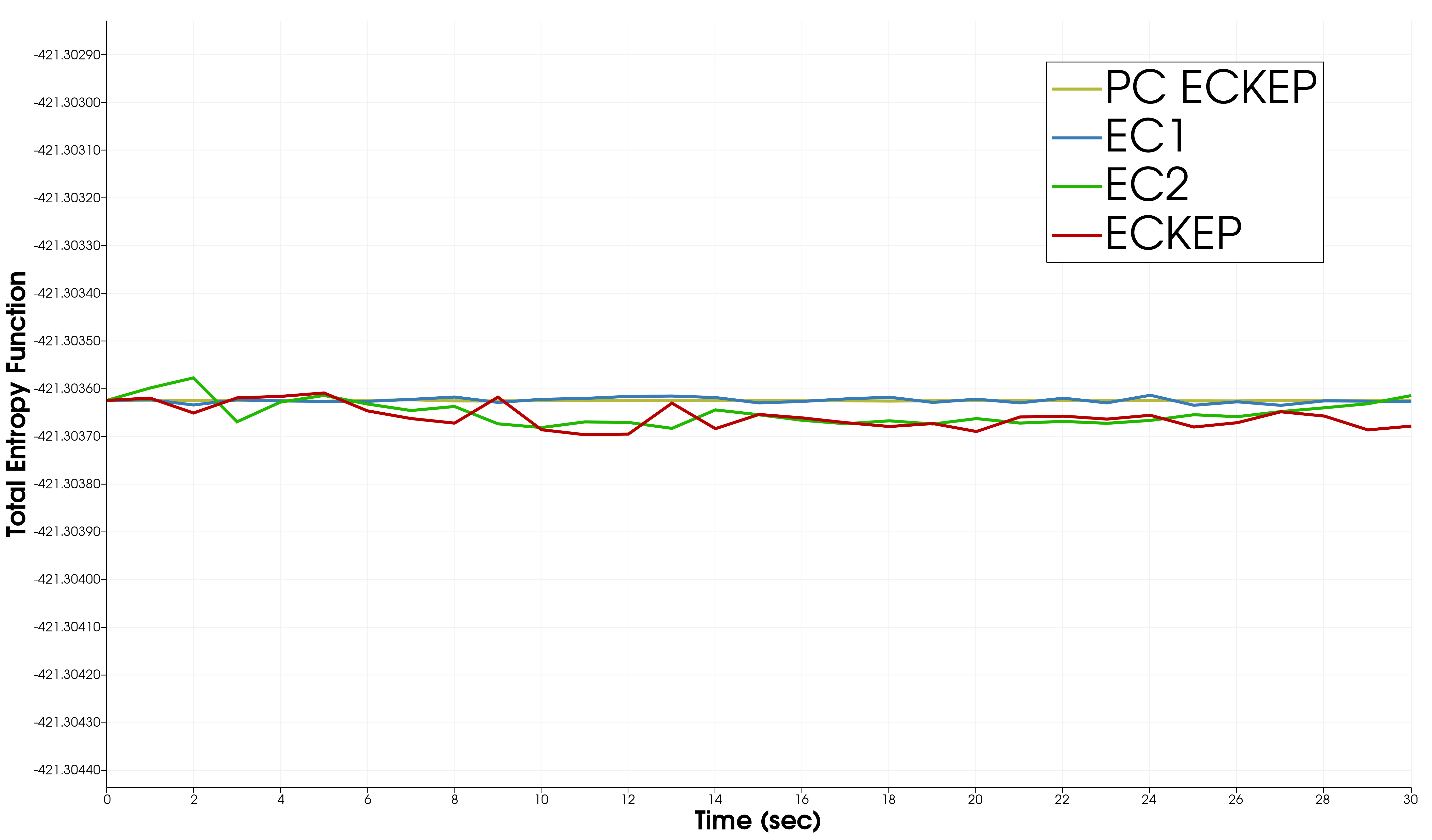} 
			\caption{Total Entropy for Taylor-Greens vortex} \label{TG_figure1}
		\end{minipage}
		\begin{minipage}{.45\textwidth}
			\centering
			\includegraphics[scale=0.15]{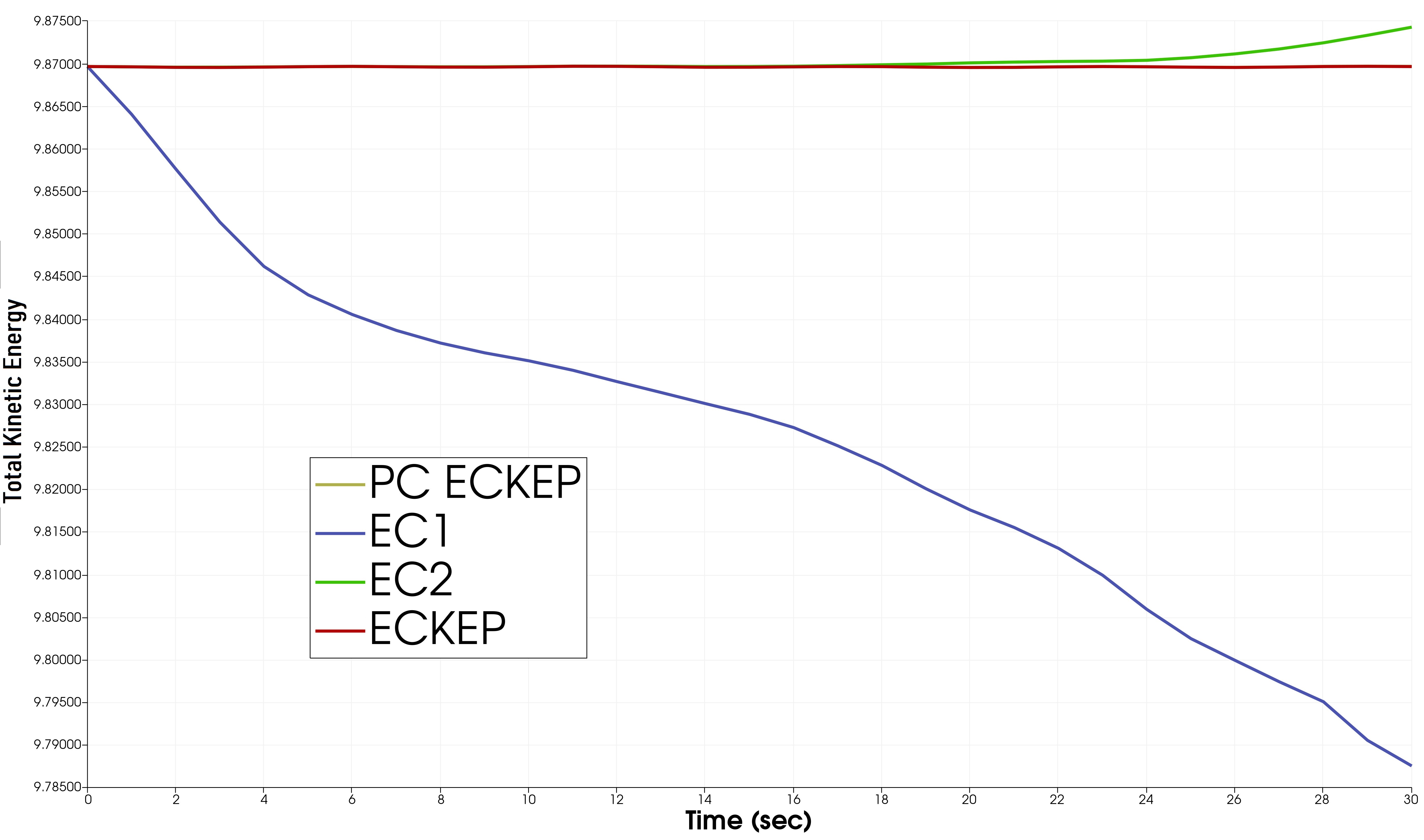}
			\caption{Total Kinetic Energy for Taylor-Greens vortex} \label{TG_figure2}
		\end{minipage}\\[1em]
	\end{figure}

	Taylor-Green vortex problem is a shock-free test case originally introduced for low-speed compressible turbulent flows. It has also been used to test the inviscid portion of Navier-Stokes equations and shows conservation properties of numerical schemes \cite{yee_2020}. A $[2\pi \times 2\pi]$ computational domain is taken with $100\times 100$ grid size with periodic boundary conditions. The problem was run for an extended time of 20 seconds. Initial conditions are given by
	\begin{equation}
		\begin{aligned}
			\rho&=1, \\
			u&=\sin(x)\times cos(y),\\
			v&=-\cos(x)\times sin(y), \\
			p&=\frac{100}{\gamma}+\frac{\cos(2x)+\cos(2y)}{4}
		\end{aligned}
	\end{equation}
	Conservation quantities are only spatially conserved, so it is desirable to have a high order time-stepping with a low CFL number. Thus, third-order SSPRK \cite{shu_1988} was implemented for temporal integration with a CFL of 0.1. The objective was to compare numerically different conservative fluxes. Comparison is made of EC1, EC2, ECKEP, and the entropy conservative and kinetic energy preserving flux of Chandrasekar \cite{chandrashekar_2013} (referred to as \textit{PC ECKEP}). It can be seen that all fluxes preserve the total entropy in domain, $\sum\limits_{\forall j} \eta_j$, with sufficient accuracy for the given time, as shown in figure \ref{TG_figure1}. Total kinetic energy, $\sum\limits_{\forall j} \rho_j u_j^2/2$, keeps reducing for EC1 from the very start of the solution and diverges for EC2 at a later time (approximately 25 sec) which is expected as these schemes are not designed to be necessarily kinetic energy preserving. Schemes ECKEP and PC ECKEP are kinetic energy preserving and no change is observed as evident from figure \ref{TG_figure2}.     Only the semi-discrete cases are considered for preserving the structures in this study.  Extending these ideas to the fully discrete case is more involved, is some times known to generate oscillations (for example, see \cite{Gousami_Murman_Duraisamy}) and is not attempted in this work.  
	\subsubsection{Isentropic vortex convection}
	Isentropic vortex is another shock-free test case with the following initial data, $\mathbf{W}_o$ \cite{spiegel_2015}.
	\begin{equation}\label{isen_vortex}
		\begin{aligned}
			\rho(x,y)&=(1-\frac{\gamma-1}{2}\omega^2)^{\frac{1}{\gamma-1}}, \\
			u(x,y)&=M \cos(\theta)-\frac{(y-y_c)\omega}{R},\\
			v(x,y)&=M \sin(\theta)+\frac{(x-x_c)\omega}{R}, \\
			p(x,y)&=\frac{1}{\gamma}\rho^{\gamma}
		\end{aligned}
	\end{equation}
	Here $\omega$ is a Gaussian function of the form $\omega=\beta \exp(f) $. The strength of the vortex is determined by $\beta$, and the perturbation function $f$ is given by
	\begin{equation*}
		f(x,y)=-\frac{1}{2\sigma^2}\left( \left(\frac{x}{R}\right)^2+ \left(\frac{y}{R}\right)^2 \right).
	\end{equation*}
	$(x_c,y_c)$ is initially the centre of the vortex.
	For this study, a square domain of 200x200 cells was taken with periodic boundaries on all sides. $\beta=1$, $R=1$ and $\sigma=1$ were taken in \eqref{isen_vortex}. Mach number is $2/\gamma$ and the flow angle $(\theta)$ is 45 degrees. The exact solution of the above problem is the above initial data translated with $M \cos(\theta)$ and $M \sin(\theta)$, {\em i.e.}, $\mathbf{W}(x,y)=\mathbf{W}_{o}(x-M \cos(\theta) t,y-M \sin(\theta) t)$.
	\begin{figure}[!h] 
		\begin{minipage}{.95\textwidth}
			\centering
			\includegraphics[scale=0.15]{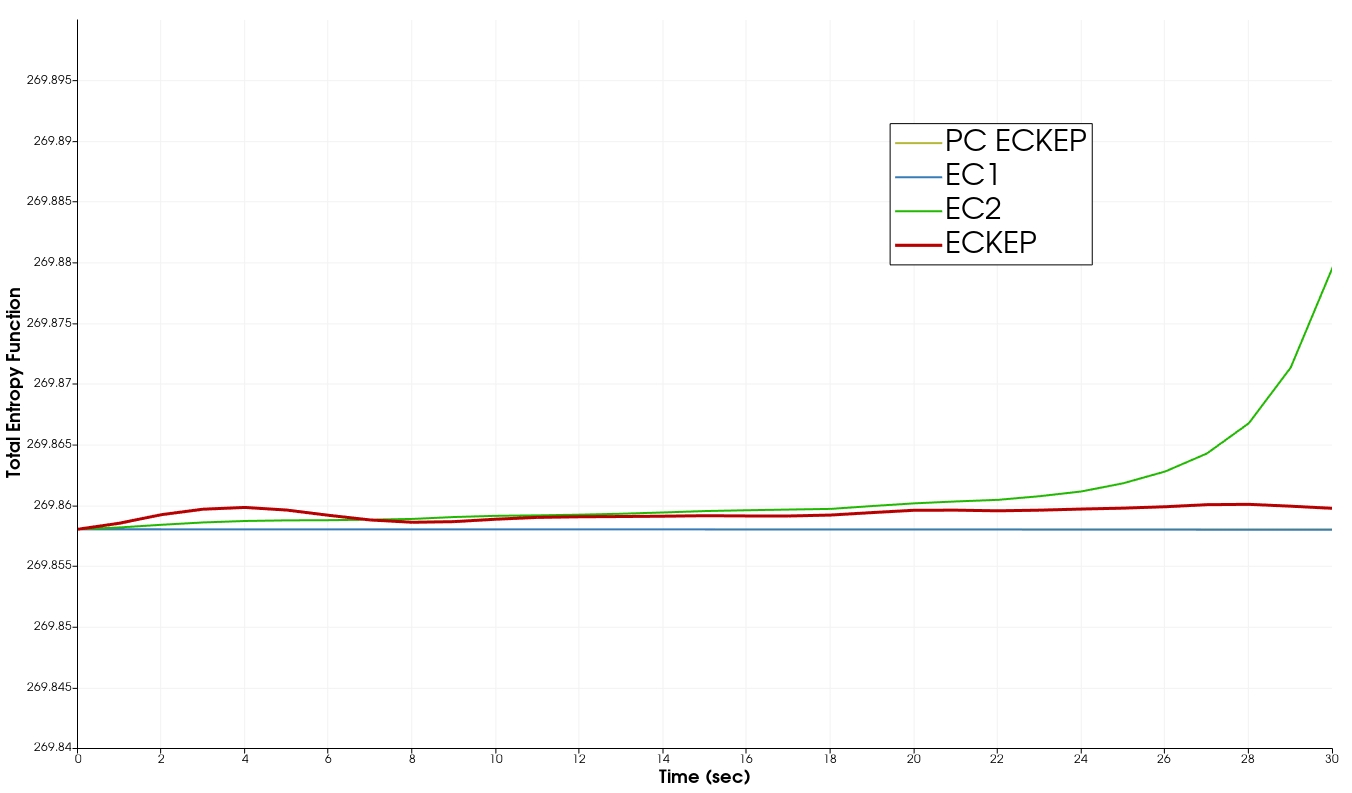} \label{IV_figure1}
			\includegraphics[scale=0.21]{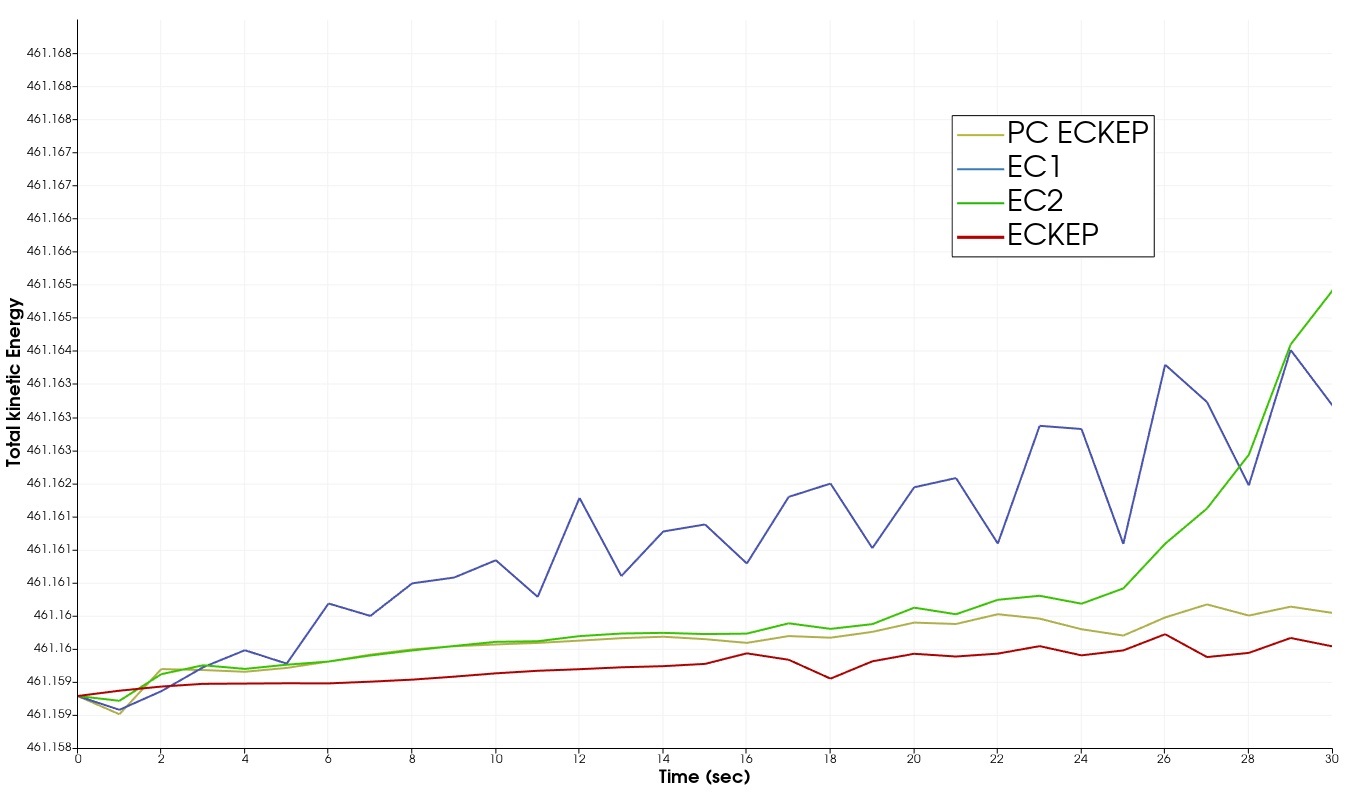}\\
		\end{minipage}\\[1em]
		\caption{Total entropy and total kinetic energy v/s Time for Isentropic vortex test case for different entropy conservative fluxes}\label{IV_figure} 
	\end{figure}  
    As shown in figure \ref{IV_figure}, entropy remains constant throughout the simulations for PC ECKEP, EC1, and ECKEP schemes and for EC2 scheme up to t=25 s. Divergence of EC2 scheme is unexpected but can be attributed to the accumulation of errors in temporal discretization (as the entropy conservation formulation is forced only for the semi-discrete case).  Kinetic energy preserving schemes, PC ECKEP and ECKEP perform as expected and kinetic energy remains constant as shown in figure \ref{IV_figure}. Kinetic energy diverges with time for EC1 and EC2 schemes since these schemes do not necessarily preserve kinetic energy numerically.   
	\subsection{Computational Cost}
	The entropy conservative fluxes presented do not need to evaluate computationally expensive averages (unlike in Riemann solvers) and are expected to be more efficient. To test this, we ran the smooth test case problem given in \eqref{den_wave} for 1000 computational cells and with SSPRK-2 discretization of temporal derivatives for a final time of 10 seconds. This resulted in 47,477,430 flux function calls for all the schemes, and the average time per call is given in table \ref{cpu_time}. Out of all entropy conservative schemes, the EC2 scheme is consecutively 29\%, 11\%, and 5\% more efficient than ROE EC (referring to Roe's entropy conservative scheme \cite{ismail_2009}), ECKEP, and EC1 schemes, respectively. Also, the new ECKEP scheme is computationally 8\% more efficient than the PC ECKEP scheme.
	\begin{table}[!h]
		\centering
		\begin{tabular}{ l c}
			\toprule
			\textbf{Flux} & \textbf{Average time per call  (microseconds) } \\
			\midrule
			ROE EC & 0.12947 \\
			PC ECKEP & 0.10401 \\
			EC1 & 0.09648  \\
			EC2 & 0.09215  \\
			ECKEP & 0.09532  \\
			\bottomrule
		\end{tabular}
		\caption{Wall clock times for various entropy conservative schemes}
		\label{cpu_time}
	\end{table}
	\section{Rankine–Hugoniot condition satisfying Entropy Stable Scheme}
	Entropy conservative flux is only valid at smooth regions, and additional diffusion must be added at discontinuities such as shocks. Additionally, the diffusion must generate entropy such that inequality in \eqref{entropy_inequality} is satisfied. Ismail and Roe \cite{ismail_2009} gave one possible entropy stable diffusion, which linearizes Euler equations about an average state.  More recently, Chandrashekar \cite{chandrashekar_2013} introduced a diffusion which satisfies both entropy stability and kinetic energy preservation.  However, the above types of diffusion depend on the underlying eigen-structure of Euler equations. With the motivation to avoid complicated Riemann solvers and their strong dependence on the eigen-structure, the simple central scheme (which still catpures steady discontinuitires exactly) introduced by Jaisankar and Raghurama Rao \cite{jaisankar_2009} (Method of optimal viscosity for enhanced resolution of shocks, MOVERS) is utilized as a basic framework for the additional diffusion proposed here. The above scheme satisfies the Rankine-Hugonoit jump conditions across an interface in determining numerical diffusion. Further advantage of this approach is that this strategy can be extended to any system of hyperbolic conservation laws with an R-H-like condition. The general form of diffusive flux for this central scheme is given by
	\begin{equation}
	\mathbf{F}^d_{j+\frac{1}{2}}=-\frac{1}{2}\mathbf{\tilde{D}}_{j+\frac{1}{2}} (\mathbf{U}_{j+1}-\mathbf{U}_j)
	\end{equation}
	For MOVERS-n the diffusion is given by   $\mathbf{\tilde{D}}_{j+\frac{1}{2}}=diag(s^1,s^2,s^3)_{j+\frac{1}{2}}$ where wave-speeds $s^k$ are computed using $\Delta \mathbf{F}=\mathbf{\tilde{D}} \cdot \Delta \mathbf{U}$. $\mathbf{\tilde{D}}_{j+\frac{1}{2}}$ is given by
	\begin{equation*}
		\mathbf{\tilde{D}}_{j+\frac{1}{2}}=
		 \left[ \begin{matrix}
		\left| \frac{F^1_{j+1}-F^1_{j}}{U^1_{j+1}-U^1_{j}} \right|  &0 &0 \\ 
		0& \left| \frac{F^2_{j+1}-F^2_{j}}{U^2_{j+1}-U^2_{j}} \right| &0 \\
		0& 0&\left| \frac{F^3_{j+1}-F^3_{j}}{U^3_{j+1}-U^3_{j}} \right| \end{matrix}  \right]
	\end{equation*}
		To prevent unphysical values of wave speeds when $\Delta U^k$ is very small, a wave-speed correction is introduced such that 
	\begin{equation*}
		s^k_{j+\frac{1}{2}}=sign(s^k_{j+\frac{1}{2}}) \lambda_{max} \quad \text{if} \quad |s^k_{j+\frac{1}{2}}| \geq \lambda_{max}
	\end{equation*}
	\begin{equation*}
		s^k_{j+\frac{1}{2}}=sign(s^k_{j+\frac{1}{2}}) \lambda_{min} \quad \text{if} \quad |s^k_{j+\frac{1}{2}}| \leq \lambda_{min}
	\end{equation*} 
        where $\lambda_{min}=\min(|u-a|,|u|,|u+a|)$ and $\lambda_{max}=\max(|u-a|,|u|,|u+a|)$.  
	When $U_{j+1}^k \rightarrow U_{j}^k$ wave speed is taken as $\lambda_{min}$.  
	We propose an accurate discontinuity capturing version of MOVERS flux (RH) by taking the following form of diffusion flux.  
	\begin{equation}\label{SM}  
	\mathbf{F}^{RH}_{j+\frac{1}{2}}=-\frac{1}{2} \min(s_1,s_2,s_3)_{j+\frac{1}{2}} \tilde{I} (\mathbf{U}_{j+1}-\mathbf{U}_j)=-\frac{1}{2} (\alpha^{S})_{j+1/2} \tilde{I} (\mathbf{U}_{j+1}-\mathbf{U}_j)
	\end{equation}
    A fix is applied to the diffusion to enable a smooth transition at the sonic point and prevent a sonic glitch. 
    \begin{equation}
        \alpha^S=\frac{(\alpha^S)^2+\Theta^2}{2 \Theta} \quad \text{where} \quad \Theta=0.1
    \end{equation}
    The above fix is only used if $\alpha^S$ is not zero.
        Across any isolated contact discontinuity moving with a speed of $v_d$ in one dimension, all the wave speeds computed using $\Delta \mathbf{F}/\Delta \mathbf{U}$ collapse to the speed of the contact discontinuity, {\em i.e.}, $s_1=s_2=s_3=v_d$. This allows us to preserve exactly steady (and grid-aligned in multi-dimensions) contact discontinuities since  d $v_d=0$ for them and thus no diffusion is added. For moving contact discontinuities, the diffuion is proportional to the speed of discontinuity. This gives a better wave speed estimate than the Local-Lax Friedrichs (Rusnaov) scheme and additionally preserves steady (and grid-aligned)  discontinuities.  
	An entropy stable flux at interface $x_{j+\frac{1}{2}}$ is now given by 
	\begin{equation}\label{ES_Scheme1}
		\mathbf{F}_{j+\frac{1}{2}}=\mathbf{F}^{EC}_{j+\frac{1}{2}}+\mathbf{F}^{RH}_{j+\frac{1}{2}}
	\end{equation}
        The above combination will preserve steady contact discontinuities exactly since both the entropy conservative fluxes (as discussed in section 4.1) and the entropy stable dissipation (based on R-H conditions) preserves steady contact discontinuities exactly.  Even though the entropy stable dissipation also preserves steady shocks, this property is not satisfied by the entropy conservative fluxes.   
	We can show that this leads to an entropy stable scheme by taking the dot product of the semi-discrete conservation law by entropy variable $\mathbf{V}_j$.
	\begin{equation*}
		\Delta x \frac{d\eta(\mathbf{U})_j}{dt}=-\mathbf{V}_j \cdot (\mathbf{F}_{j+\frac{1}{2}}-\mathbf{F}_{j-\frac{1}{2}})	
	\end{equation*}
	\begin{equation*}
		\Delta x \frac{d\eta(\mathbf{U})_j}{dt}=-\left(\frac{1}{2}(\mathbf{V}_{j+1}+\mathbf{V}_j)-\frac{1}{2}(\mathbf{V}_{j+1}-\mathbf{V}_{j})\right) \cdot \mathbf{F}_{j+\frac{1}{2}}+\left(\frac{1}{2}(\mathbf{V}_{j}+\mathbf{V}_{j-1})+\frac{1}{2}(\mathbf{V}_{j}-\mathbf{V}_{j-1})\right)\mathbf{F}_{j-\frac{1}{2}}
	\end{equation*}
	\begin{equation*}
		\begin{aligned}
		\Delta x \frac{d\eta(\mathbf{U})_j}{dt}=&-\left[
		\frac{1}{2}(\mathbf{V}_{j+1}+\mathbf{V}_j)\cdot \mathbf{F}^{EC}_{j+\frac{1}{2}}-\frac{1}{2}(\mathbf{V}_{j+1}-\mathbf{V}_j)\cdot\mathbf{F}^{EC}_{j+\frac{1}{2}}-\frac{1}{2}(\mathbf{V}_{j+1}+\mathbf{V}_j)\cdot (\alpha^{S})_{j+\frac{1}{2}} (\mathbf{U}_{j+1}-\mathbf{U}_j) \right. \\ &
		-\frac{1}{2}(\mathbf{V}_{j}+\mathbf{V}_{j-1})\cdot \mathbf{F}^{EC}_{j-\frac{1}{2}}+\frac{1}{2}(\mathbf{V}_{j}-\mathbf{V}_{j-1})\cdot\mathbf{F}^{EC}_{j-\frac{1}{2}}
		+\frac{1}{2}(\mathbf{V}_{j}+\mathbf{V}_{j-1}) (\alpha^{S})_{j-\frac{1}{2}} (\mathbf{U}_{j}-\mathbf{U}_{j-1}) \\ &
		\left. +\frac{1}{4}(\mathbf{V}_{j+1}-\mathbf{V}_j)(\alpha^{S})_{j+\frac{1}{2}}(\mathbf{U}_{j+1}-\mathbf{U}_{j})
		+\frac{1}{4}(\mathbf{V}_{j}-\mathbf{V}_{j-1})(\alpha^{S})_{j-\frac{1}{2}}(\mathbf{U}_{j}-\mathbf{U}_{j-1})
		 \right]
		\end{aligned}
	\end{equation*}
	Using the entropy conservation relation $\Delta V \cdot F^{EC}=\Delta \psi$ we get
	\begin{equation}\label{ent_eq_sd}
		\begin{aligned}
			\Delta x \frac{d\eta(\mathbf{U})_j}{dt}+\zeta_{j+\frac{1}{2}}-\zeta_{j-\frac{1}{2}}=- 
			\left[ +\frac{1}{4}(\mathbf{V}_{j+1}-\mathbf{V}_j)(\alpha^{S})_{j+\frac{1}{2}}(\mathbf{U}_{j+1}-\mathbf{U}_{j}) \right. &\\ \left.
			+\frac{1}{4}(\mathbf{V}_{j}-\mathbf{V}_{j-1})(\alpha^{S})_{j-\frac{1}{2}}(\mathbf{U}_{j}-\mathbf{U}_{j-1})
			\right]
		\end{aligned}
	\end{equation}
	where the numerical entropy flux $\zeta_{j+\frac{1}{2}}$ is given by
	\begin{equation}
		\zeta_{j+\frac{1}{2}}=\frac{1}{2}(\mathbf{V}_{j+1}+\mathbf{V}_j)\cdot \mathbf{F}^{EC}_{j+\frac{1}{2}}-\frac{1}{2}(\psi_{j+1}+\psi_j)-\frac{1}{2}(\mathbf{V}_{j+1}+\mathbf{V}_j)\cdot (\alpha^{S})_{j+\frac{1}{2}} (\mathbf{U}_{j+1}-\mathbf{U}_j)
	\end{equation}
	In equation \ref{ent_eq_sd}, the RHS is always negative given that $\alpha^{S}$ is always positive and since $\Delta \mathbf{V} \cdot \Delta \mathbf{U}>=0$ (See appendix B). Thus, the scheme always leads to reduction in mathematical entropy and is thus entropy stable. Scheme (given by \eqref{ES_Scheme1}) can also be shown to be kinetic energy stable, {\em i.e.}, it prevents kinetic energy from growing spuriously. Taking ECKEP flux as the central flux in \eqref{ES_Scheme1}, we get the discrete kinetic energy preservation equation as 
	\begin{equation}\label{semi_dis_ke_eq}
            \begin{aligned}
		\Delta x_j \frac{d}{dt}\left(\frac{\rho_j u_j^2}{2}\right)+\left(\frac{\rho u^3}{2}\right)_{j+\frac{1}{2}}-\left(\frac{\rho u^3}{2}\right)_{j-\frac{1}{2}}&+u_j(p_{j+\frac{1}{2}}-p_{j-\frac{1}{2}})=\\&-\frac{1}{4}\left[(\alpha^{S})_{j+\frac{1}{2}}\rho_{j+1}(u_{j+1}-u_j)^2+(\alpha^{S})_{j-\frac{1}{2}}\rho_{j-1}(u_j-u_{j-1})^2\right]
            \end{aligned}
	\end{equation}
	Here, the numerical kinetic energy flux can be written as follows.
        \begin{equation}
            \left(\frac{\rho u^3}{2}\right)_{j+\frac{1}{2}}=\overline{(\rho u)}_{j+\frac{1}{2}}\frac{u_{j+1} u_j}{2}-\frac{1}{4} (\alpha^{S})_{j+\frac{1}{2}}(\rho_{j+1} u_{j+1}^2-\rho_{j} u_{j}^2)
        \end{equation}
 The above flux is consistent and conservative. The terms on the RHS of \eqref{semi_dis_ke_eq} are always negative (given $\alpha^{S} \geq 0$) and thus lead to kinetic energy stability.
 In smooth regions where $\Delta \mathbf{U} \rightarrow 0$, the diffusion coefficient for MOVERS type flux tends to $\lambda_{max}$ because of wave speed correction. This results in unnecessary high diffusion in smooth regions like expansion fans. To mitigate this problem, in the next section, we use a shock sensor to apply RH condition-based diffusion only at shocks. This allows low diffusion in smooth regions while maintaining accurate shock capturing.
	\section{Second order hybrid entropy stable scheme}
\begin{figure}[!h]
		\begin{subfigure}{0.5\textwidth}
			\includegraphics[scale=0.2]{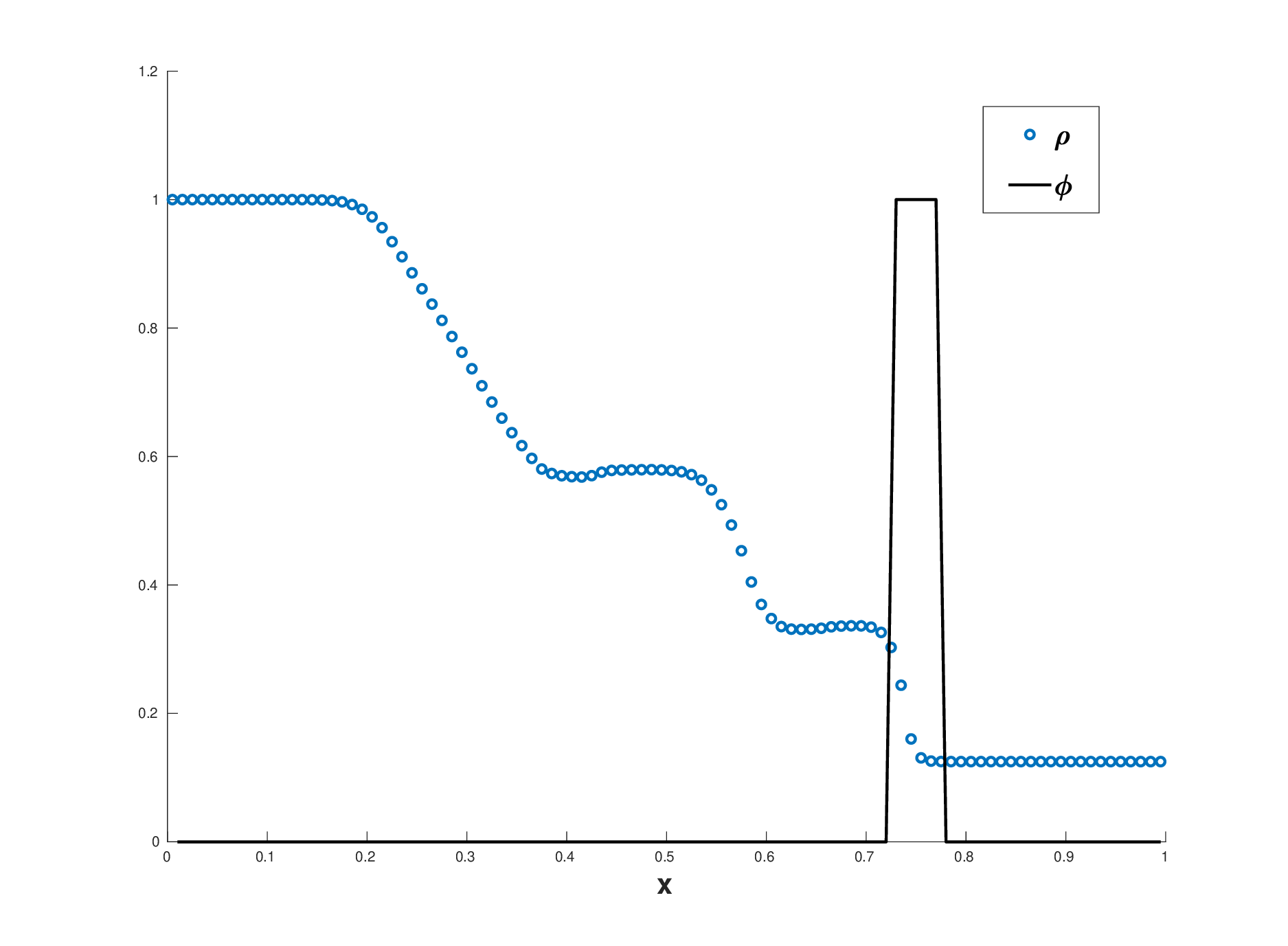}
			\caption{} \label{shock_sensora}
		\end{subfigure}%
		\hspace{-2.1cm}
		\begin{subfigure}{0.5\textwidth}
			\includegraphics[scale=0.25]{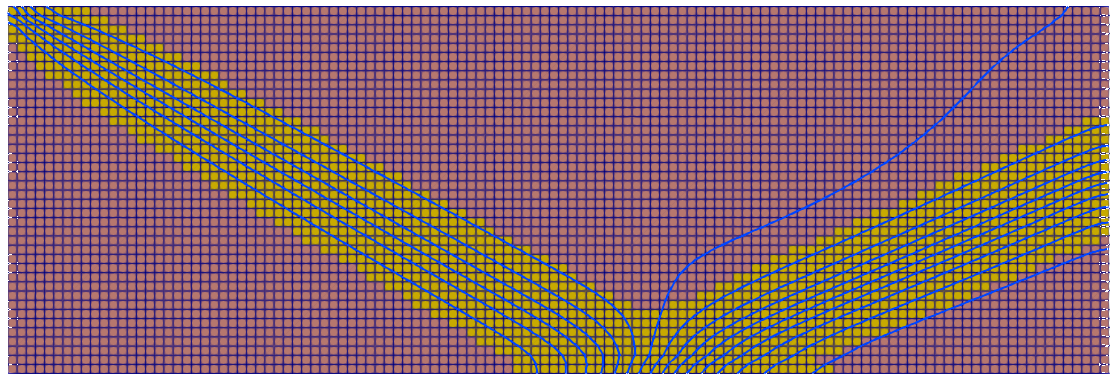}
			\caption{} \label{shock_sensorb}
		\end{subfigure}
		\caption{(a): Shock sensor value for SOD tube test case. (b): Shock sensor value and contour of density variation for oblique shock reflection test case (cells with shocks marked with yellow).} \label{shock_sensor}
	\end{figure}

		 This section proposes a hybrid scheme based on a higher-order diffusion term with entropy conservative flux, which keeps diffusion low in smooth regions. At shocks, diffusion is switched to RH condition based on entropy distance, as given by \eqref{SM}.
		 Entropy distance is defined as the dot product between the change in the entropy variable vector and the change in the conserved variable vector between two points. This senses gradients of flow variables and thus can be used to identify regions of high gradients and shocks. It has been used as a means for mesh movement and applying entropy fixes \cite{zaide_2009} \cite{shrinath_2023} \cite{roy_2023}. The entropy distance for a given interface can be given as
		 \begin{equation*}
		 	ED_{j+\frac{1}{2}}=(\mathbf{U}_{j+1}-\mathbf{U}{j}) \cdot (\mathbf{V}_{j+1}-\mathbf{V}{j})
		 \end{equation*}
		  Here, we propose an exponentially scaled entropy distance-based sensor to identify shocks, which is given as
		 \begin{equation}
		 	\phi_{j+\frac{1}{2}}=\begin{aligned}
		 		1-\left|\exp(-qSED_{j+\frac{1}{2}})\right| 
		 	\end{aligned}
		 \end{equation}
		 Here $SED$ is the scaled entropy distance given as $ED_{j+\frac{1}{2}}/\underset{\forall j}{max} (ED_{j+\frac{1}{2}})$ and q is the scaling factor. Since exponent is computationally expensive to compute, we can use quadratic approximation given as follows
		 \begin{equation}\label{ss_1}
		 	\phi_{j+\frac{1}{2}}=\begin{aligned}
		 		1-\left|1-qSED_{j+\frac{1}{2}}+(qSED_{j+\frac{1}{2}})^2\right| 
		 	\end{aligned}
		 \end{equation}
		 $\phi$ is then clipped and flattened as follows
		 \begin{equation}\label{ss_2}
		 	\phi_{j+\frac{1}{2}}=
		 	\left\{ \begin{matrix}
		 		\begin{aligned} 1 & \quad \text{if} \quad \phi_{j+\frac{1}{2}} > \epsilon \\ 0 & \quad \text{otherwise} \end{aligned}
		 	\end{matrix} \right.
		 \end{equation}
		 For most problems, $q$ is taken between 8 to 16 and $\epsilon$ is taken between $10^{-1}$ to $10^{-2}$. The above function might not detect complete shock profiles correctly for solutions on coarse grids where shock profiles get smeared over multiple cells. Thus a "flattening" of was done {\em i.e},  $\phi_{j+\frac{1}{2}}=max(\phi_{j-\frac{1}{2}},\phi_{j+\frac{1}{2}},\phi_{j+\frac{3}{2}})$. For two-dimensional problems, all the direct neighbours of a cell are considered while calculating $\phi$. Entropy distance and the shock sensor value are plotted for a sod tube test case, and as seen in figure \ref{shock_sensora}, it matches with the shock location. Figure \ref{shock_sensorb} shows the highlighted cells where the shock is present for an oblique shock reflection test case. \\
	Through numerical experiments, it was found that the entropy conservative scheme does not add sufficient stabilizing diffusion in smooth regions on coarse meshes. Additional fourth-order numerical diffusion is added based on the JST scheme \cite{jameson_1981} to prevent oscillations. The diffusion is given as
	\begin{equation}\label{4tho_diff}
		\mathbf{F}_{j+\frac{1}{2}}^R=\frac{1}{2} \alpha_R (\mathbf{U}_{j+2}-3 \mathbf{U}_{j+1}+3\mathbf{U}_{j}-\mathbf{U}_{j-1})
	\end{equation}
	The diffusion coefficient is given by
	\begin{equation}
		\alpha_R=\frac{1}{32} \tilde{\lambda}_{j+\frac{1}{2}}
	\end{equation}
	and $\tilde{\lambda}_{j+\frac{1}{2}}$ is the coefficient of numerical diffusion based on that of the Riemann invariants based exact contact discontinuity capturing scheme (RICCA) \cite{kolluru_2022}. The following equation defines it.
	\begin{equation}
		\tilde{\lambda}_{j+\frac{1}{2}}=\left\{ \begin{matrix}
			\begin{aligned} \frac{|u_{j+1}|+|u_j|}{2} & \quad \quad \quad \qquad \text{if} \quad |\mathbf{F}_{j+1}-\mathbf{F}_j| \leq \delta \quad \text{and} \quad  |\mathbf{U}_{j+1}-\mathbf{U}_j| \leq \delta \\
			\max(|u_{j+1}|,&|u_{j}|) +sign(|\Delta p|)   a_{j+\frac{1}{2}} \quad \quad \quad \text{otherwise} 
			\end{aligned}
		\end{matrix} \right.
	\end{equation}
    Here $\delta$ is small parameter taken as $10^{-16}$.
	The above diffusion and the entropy conservative scheme result in a numerical scheme that can capture steady contact discontinuities exactly.
	The final scheme can be written as
	\begin{equation}\label{es_final}
		\mathbf{F}=\mathbf{F}^{EC}+(1-\phi)\mathbf{F}^{R}+\phi \mathbf{F}^{RH}
	\end{equation}
	The above hybrid scheme is spatially second-order accurate in smooth regions as the central entropy conservative flux is second-order accurate, as shown in section 4.2. In the following sections, the schemes given by \eqref{ES_Scheme1} and \eqref{es_final} will be referred to as ES (Entropy stable) and HES (Hybrid entropy stable) schemes, respectively. The following sections test them for various one-dimensional and two-dimensional test cases.
	
	\section{Numerical Results}
	\subsection{One Dimensional Results}
	Numerous 1-D Euler test cases have been solved to determine the robustness and accuracy of the entropy stable scheme given in \eqref{es_final}. Note that EC3 flux, kinetic energy preserving, was taken as the entropy conservative flux in ES \eqref{ES_Scheme1} and  HES \eqref{es_final} schemes. For the one-dimensional problems, 100 points with a CFL number of 0.1 is taken. Neumann boundary conditions are applied applied on both ends. Time-step is computed using $\Delta t = CFL \times \frac{\Delta x}{\underset{\forall j}{\max}(|u_j|+a_j)}$. The initial conditions, final time and location of initial discontinuity for all test cases are given in table \ref{1d_cases}.
		\begin{table}[!h]
		\centering
		\begin{tabular}{c c c c c c c c c}
			\toprule
			No.&$x_o$ & $\rho_L$ & $u_L$ &$p_L$& $\rho_R$ & $u_R$ &$p_R$& $t_f$ \\
			\midrule
			1&0.3 & 1 & 0.75 & 1 & 0.125 & 0 & 0.1 &0.2 \\
			2&0.5 & 1 & 0 & 1000 & 1 & 0 & 0.01 & 0.012 \\
			3&0.4 & 5.9924 & 19.5975 & 460.894 & 5.9924 & -6.19633 & 46.0950 & 0.035 \\
			4&0.5 & 1 & 1 & $\frac{1}{\gamma M^2}$ & $\frac{\frac{\gamma+1}{\gamma-1}PR+1}{\frac{\gamma+1}{\gamma-1} +PR}$ & $\frac{PR}{\gamma M^2}$ &$\sqrt{\frac{\gamma(2+(\gamma-1)M^2)u_R}{(2\gamma M^2+(1-\gamma))\rho_R}} $ &5 \\
			5&0.5&1.4&0&1&1&0&1&2\\
			6&0.1&3.86&-0.81&10.33&1&-3.44&1&4\\
			7&.5&1.4&0.1&1.0&1.0&0.1&1.0&1\\
			\bottomrule
		\end{tabular}
		\caption{Initial conditions for 1D test cases}
		\label{1d_cases}
	\end{table} \\
	Test case 1 is a sod shock tube problem with a right-moving shock, right-moving contact discontinuity and a left-moving rarefaction with a sonic point. The ES scheme is compared with the Roe EC and PC ECKEP schemes. ECKEP works well with this test case as shock and contact discontinuity are captured with less diffusion, and no expansion shock or sonic glitch is observed in the expansion fan region. It is comparable to the existing entropy stable and kinetic preserving schemes. Test cases 2 and 3 involve shocks with high-pressure gradients and are used to test the robustness of the numerical schemes. The schemes are able to give non-oscillatory results for these test cases. Test case 4 depicts a stationary shock, captured with one interior point by the ES scheme. Test case 5 represents a stationary contact discontinuity which the new schemes can capture exactly. Test case 6 deals with a slowly moving shock wave moving to the right. Numerical schemes with low numerical diffusion suffer from post-shock oscillations for this test. However, the new entropy stable schemes capture the shock without oscillations. Test case 7 represents a slowly moving contact discontinuity moving to the right. The entropy stable schemes are seen to capture the contact discontinuity with minimal diffusion and without oscillations. 
	\begin{figure}
		\begin{minipage}{1\textwidth}
			\centering
			\includegraphics[scale=0.395]{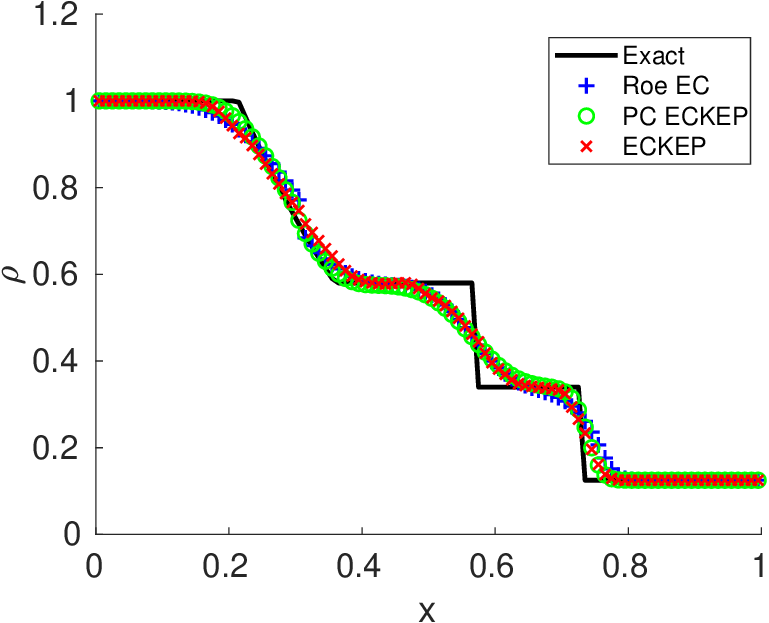}\quad
			\includegraphics[scale=0.395]{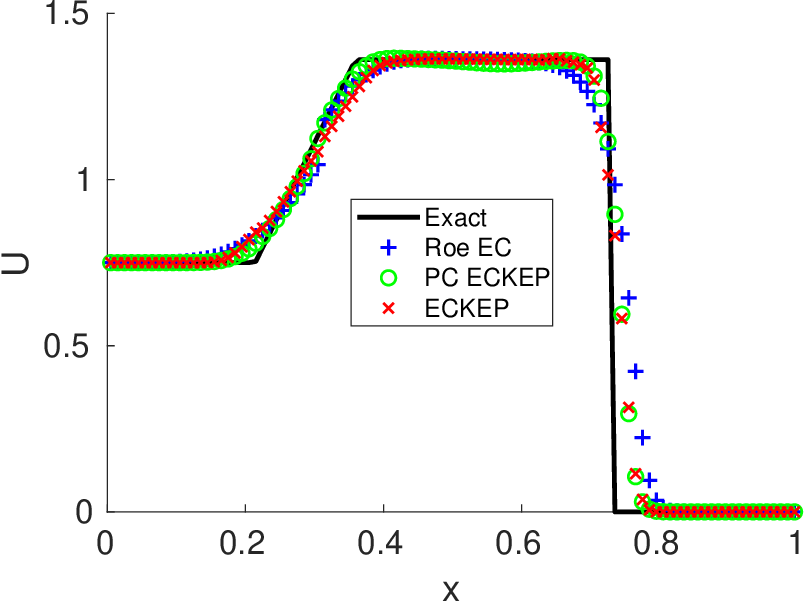}\quad
			\includegraphics[scale=0.395]{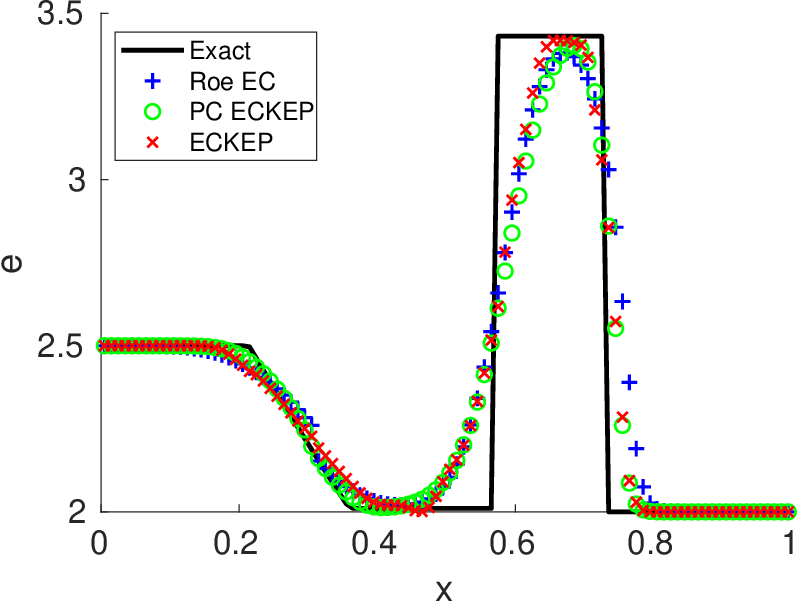}\\
		\end{minipage}
		\vspace{-2mm}
		\caption{Plots of density, velocity and internal energy for various schemes for the first test case}
		\label{One_D_1}
	\end{figure}
	\begin{figure}
		\begin{minipage}{1\textwidth}
			\includegraphics[scale=0.395]{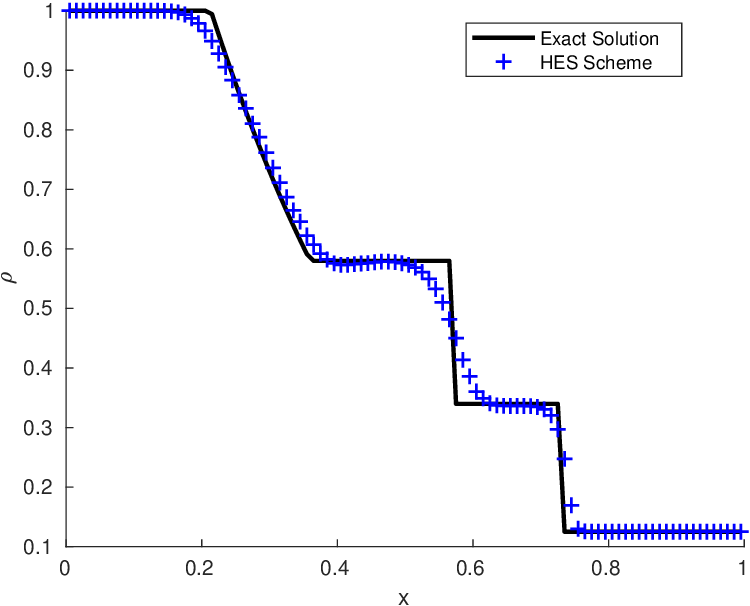}\quad
			\includegraphics[scale=0.395]{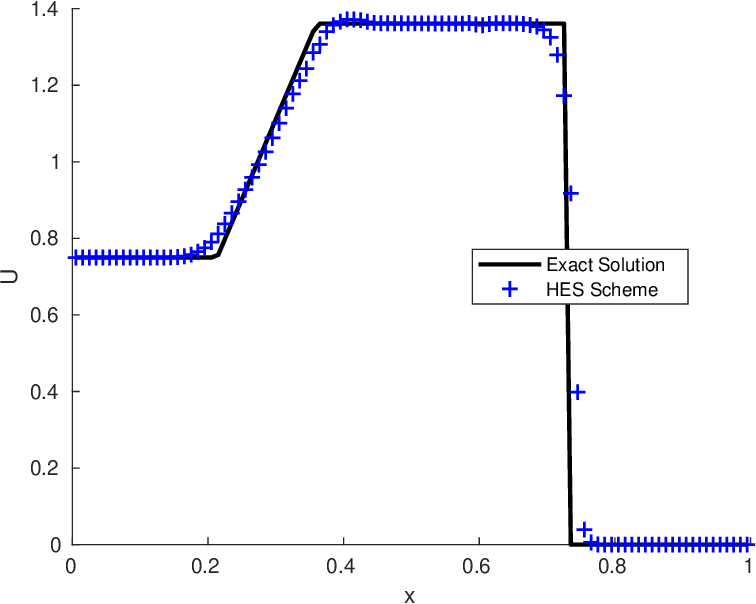}\quad
			\includegraphics[scale=0.395]{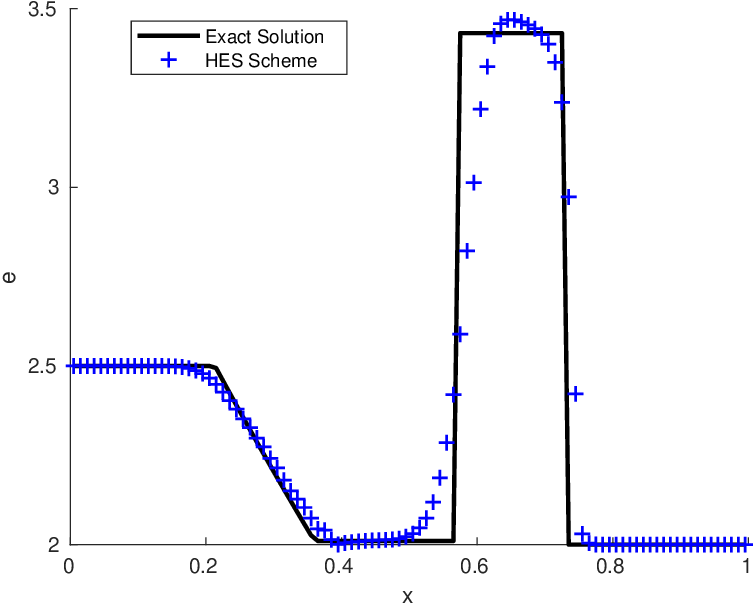}\\
		\end{minipage}
		\vspace{-6mm}
		\caption{Plots of density, velocity and internal energy using HES scheme for the first test case}
		\label{One_D_12}
	\end{figure}
	\begin{figure}
		\begin{minipage}{1\textwidth}
			\centering
			\includegraphics[scale=0.395]{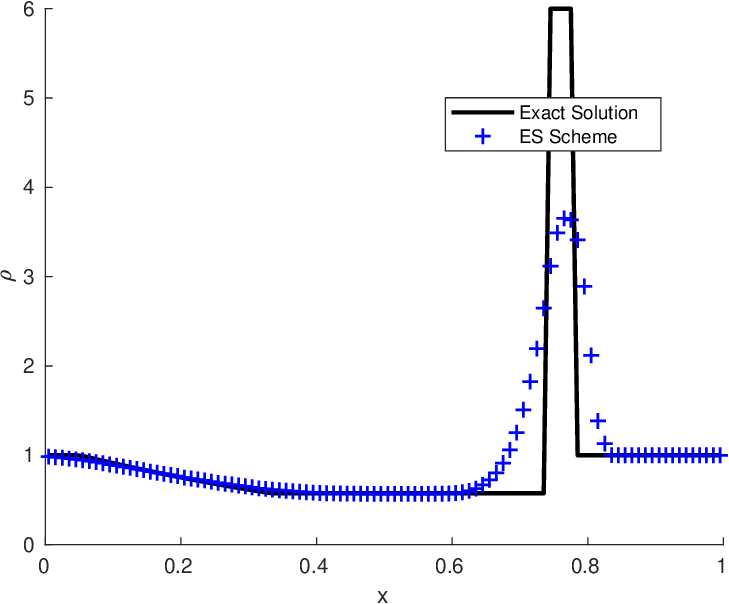}\quad
			\includegraphics[scale=0.395]{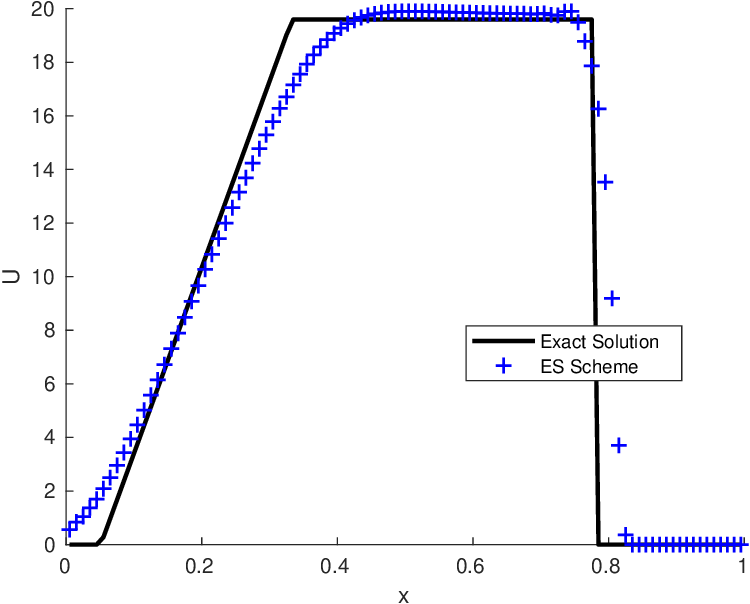}\quad
			\includegraphics[scale=0.395]{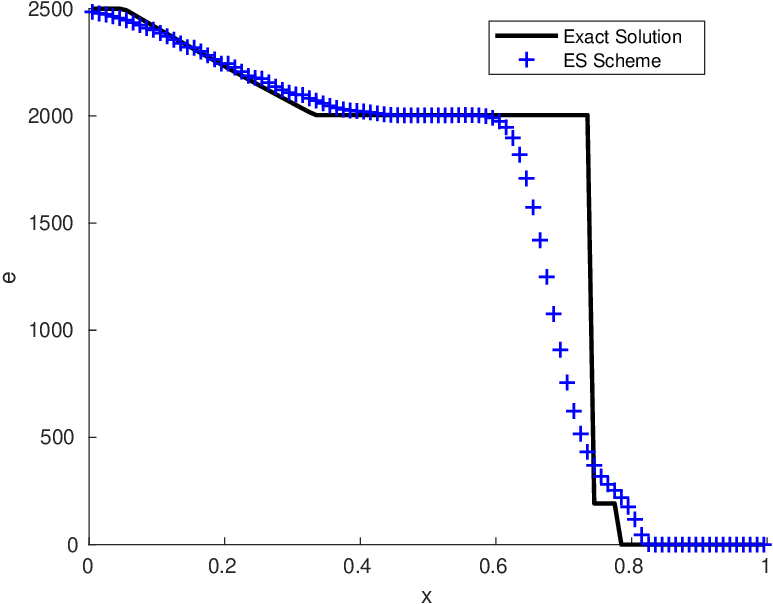}\\
		\end{minipage}
		\begin{minipage}{1\textwidth}
			\centering
			\includegraphics[scale=0.395]{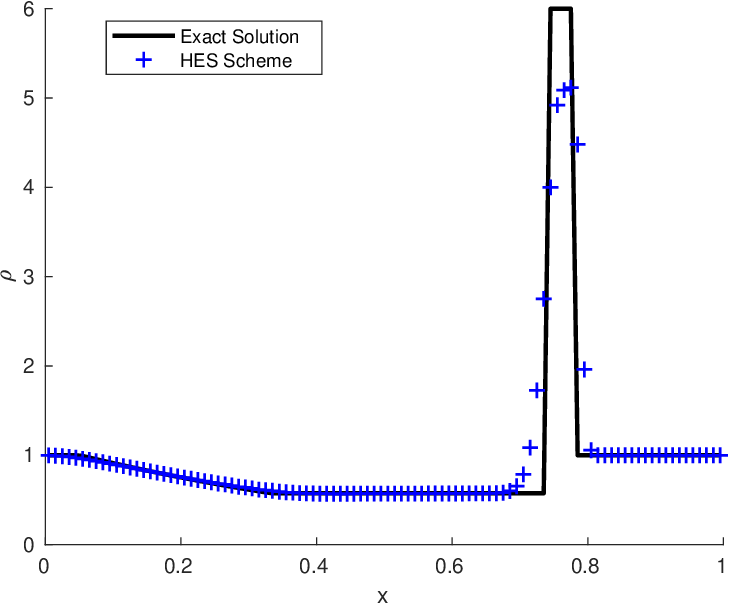}\quad
			\includegraphics[scale=0.395]{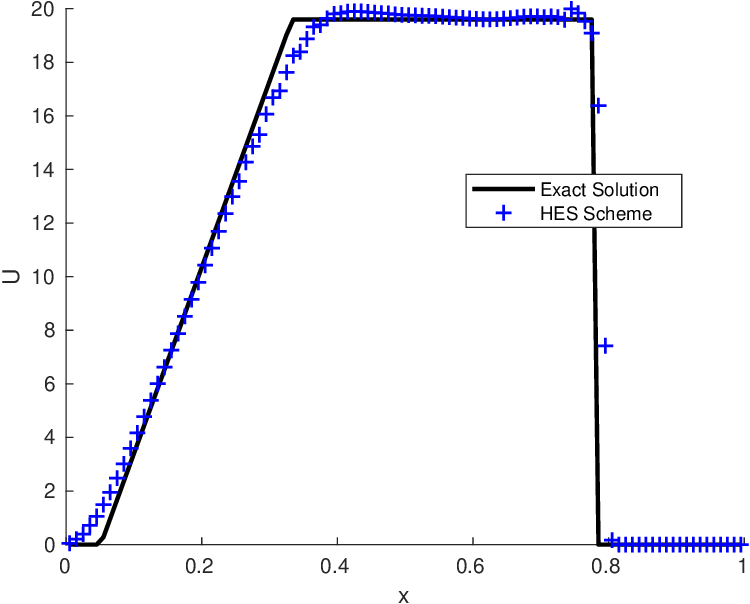}\quad
			\includegraphics[scale=0.395]{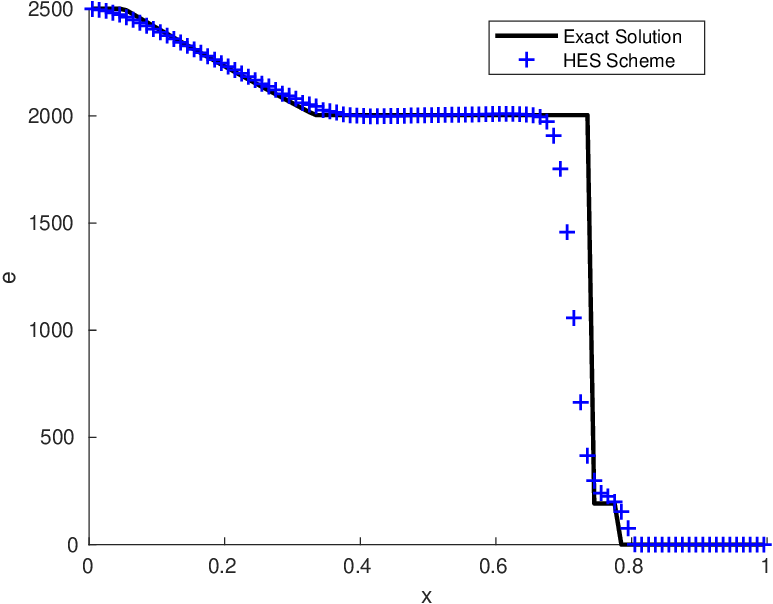}\\
		\end{minipage}
		\caption{Plots of density, velocity and internal energy using ES (top) and HES (bottom) schemes for the third test case}
		\label{One_D_3}
	\end{figure}
	\begin{figure}
		\begin{minipage}{1\textwidth}
			\centering
			\includegraphics[scale=0.395]{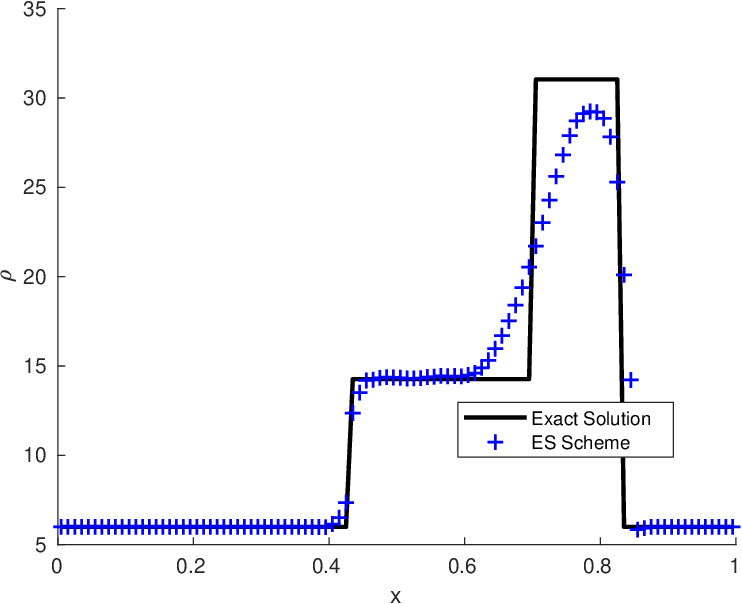}\quad
			\includegraphics[scale=0.395]{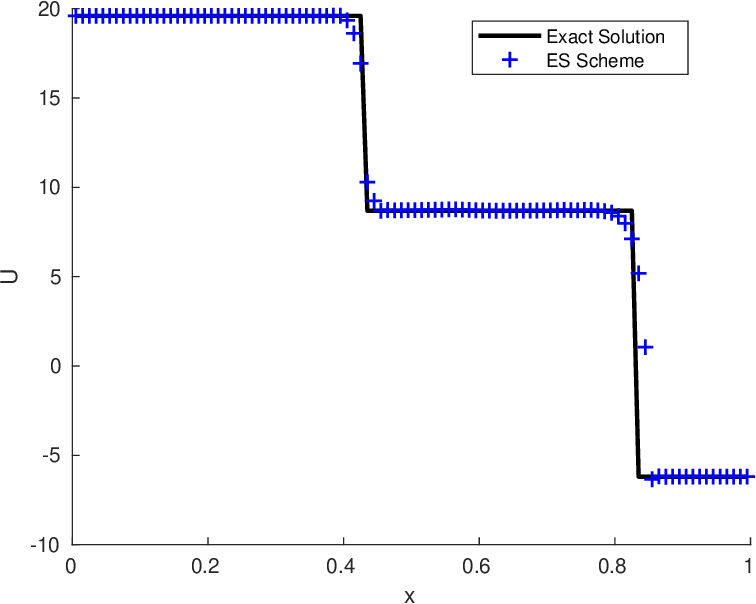}\quad
			\includegraphics[scale=0.395]{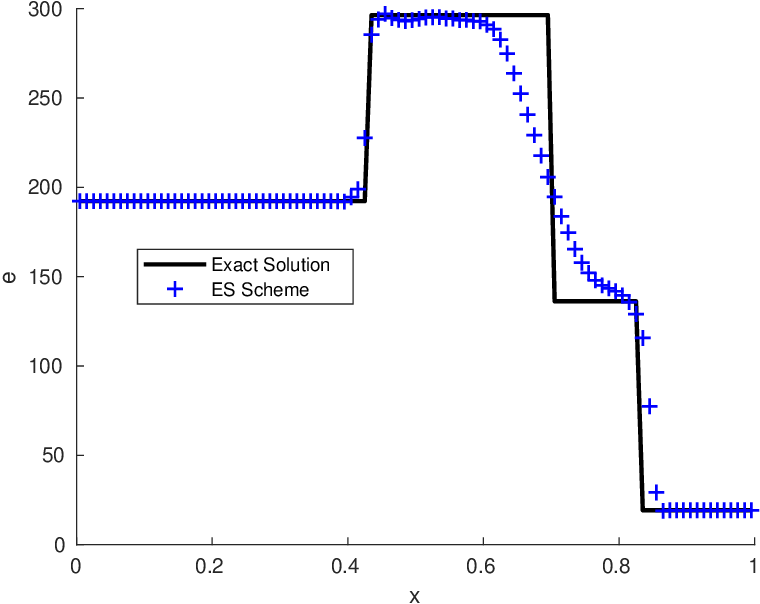}\\
		\end{minipage}
		\begin{minipage}{1\textwidth}
			\centering
			\includegraphics[scale=0.395]{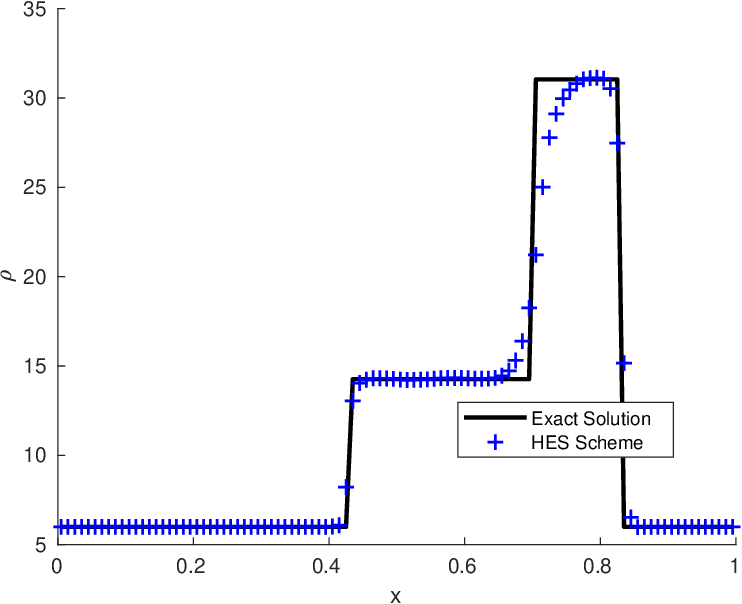}\quad
			\includegraphics[scale=0.395]{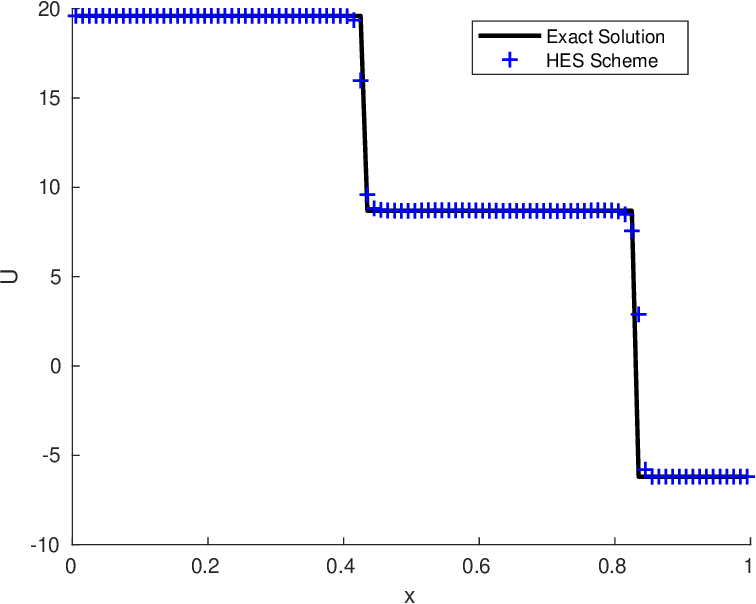}\quad
			\includegraphics[scale=0.395]{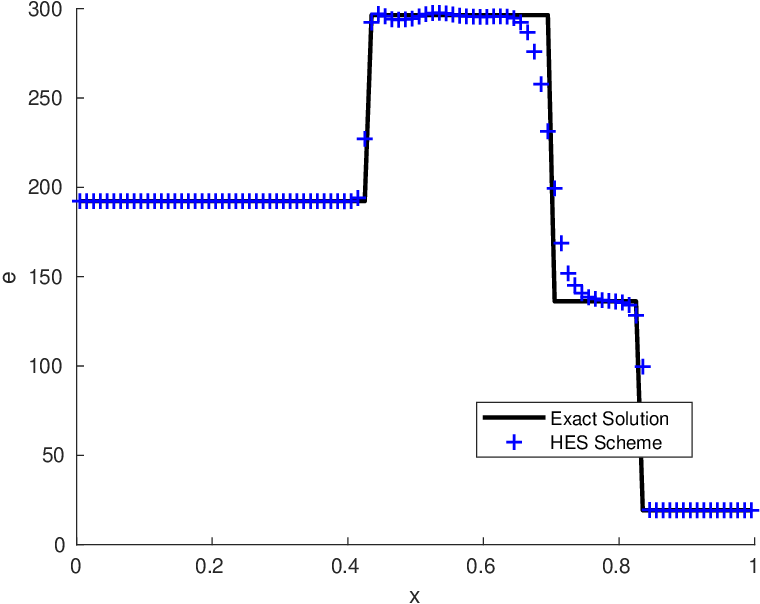}\\
		\end{minipage}
		\caption{Plots of density, velocity and internal energy using ES (top) and HES (bottom) schemes for the fourth test case}
		\label{One_D_4}
	\end{figure}
	\begin{figure}
		\begin{minipage}{1\textwidth}
			\centering
			\includegraphics[scale=0.395]{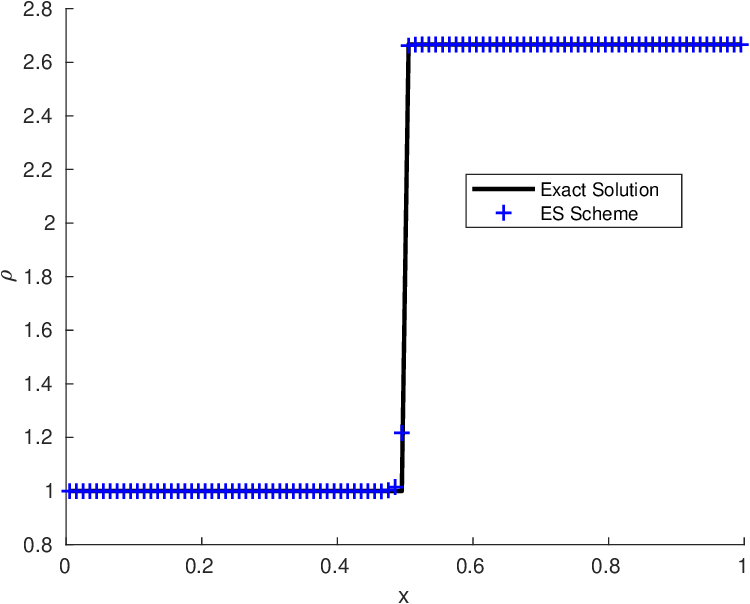}\quad
			\includegraphics[scale=0.395]{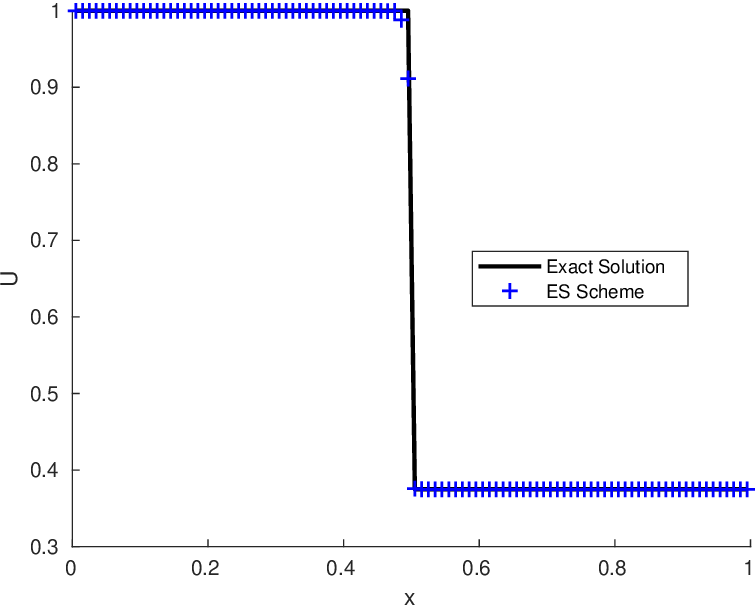}\quad
			\includegraphics[scale=0.395]{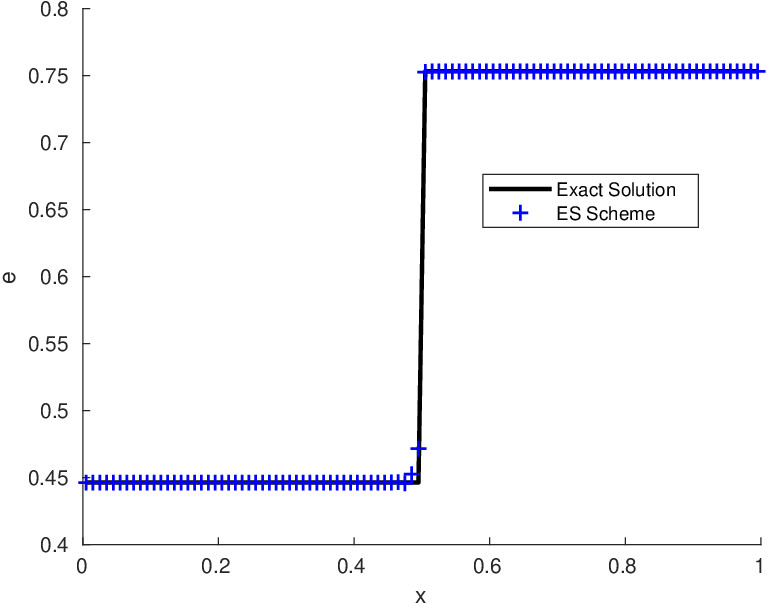}\\
		\end{minipage}
		\begin{minipage}{1\textwidth}
			\centering
			\includegraphics[scale=0.395]{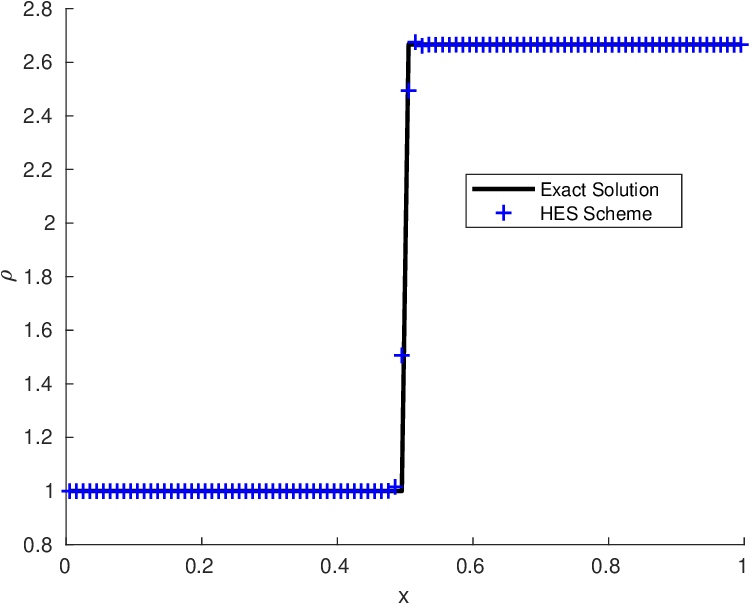}\quad
			\includegraphics[scale=0.395]{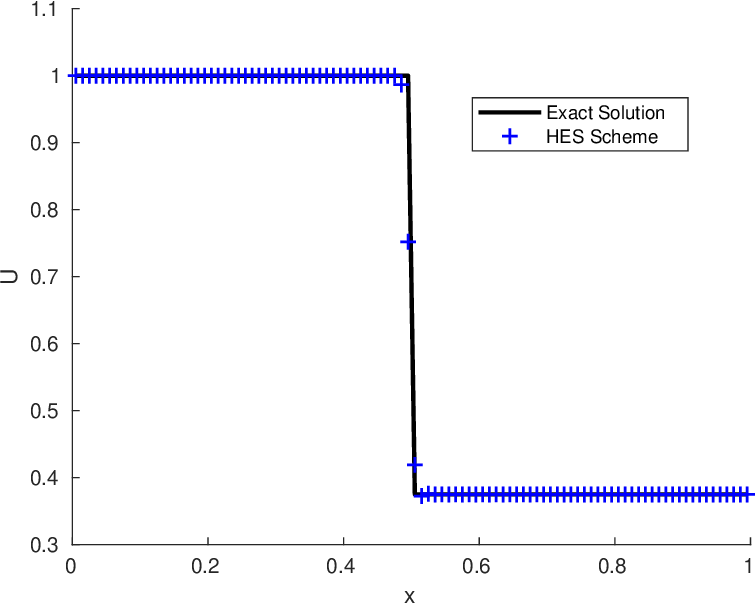}\quad
			\includegraphics[scale=0.395]{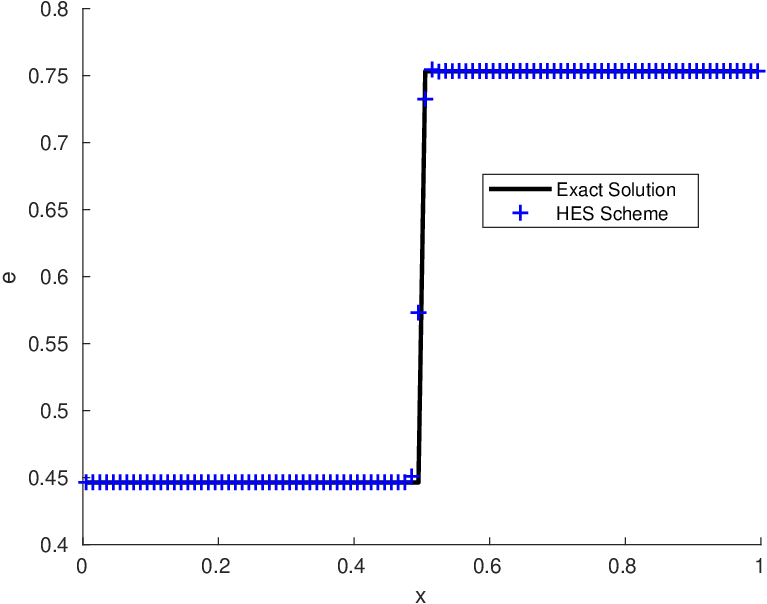}\\
		\end{minipage}
		\caption{Plots of density, velocity and internal energy using ES (top) and HES (bottom) schemes for the stationary shock wave.}
		\label{One_D_6}
	\end{figure}
		\begin{figure}
			\centering
			\includegraphics[scale=0.394]{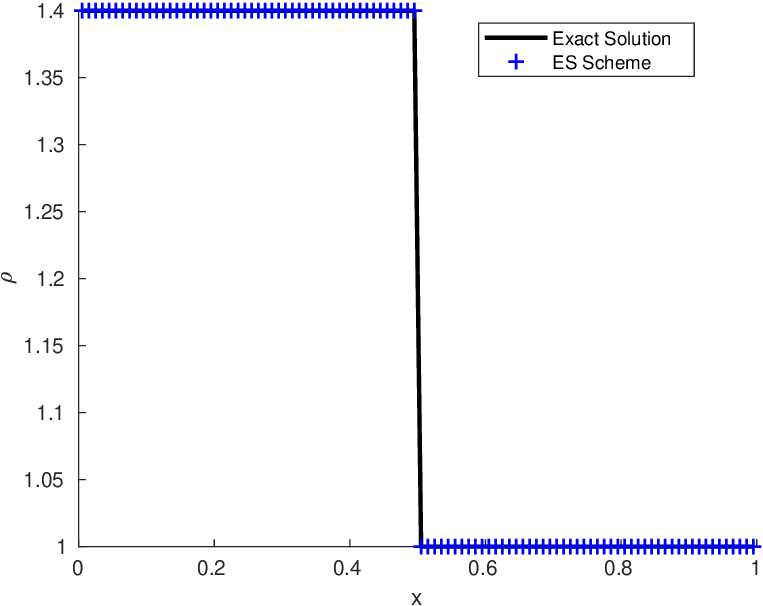}\quad
			\includegraphics[scale=0.394]{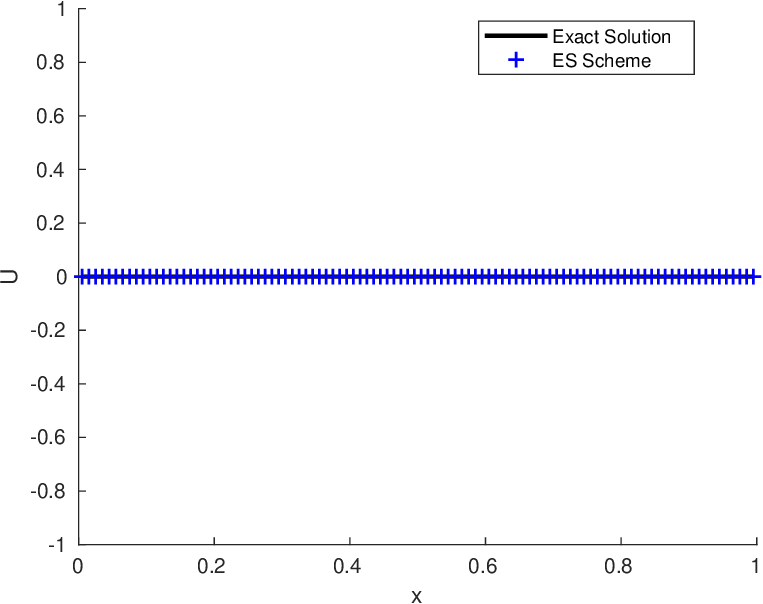}\quad
			\includegraphics[scale=0.394]{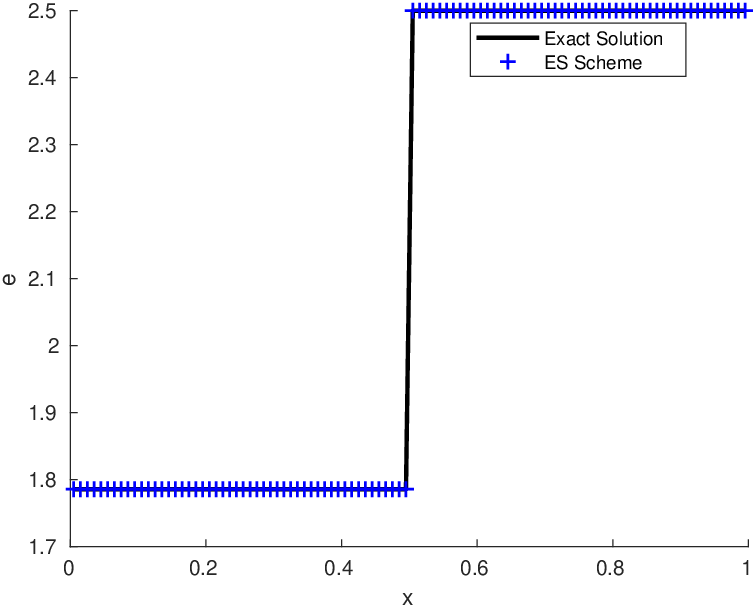}\\
		\caption{Plots of density, velocity and internal energy using ES scheme for stationary contact wave.}
		\label{One_D_7}
	\end{figure}\\
	\vspace{-20cm}
	\begin{figure}
		\begin{minipage}{1\textwidth}
			\centering
			\includegraphics[scale=0.395]{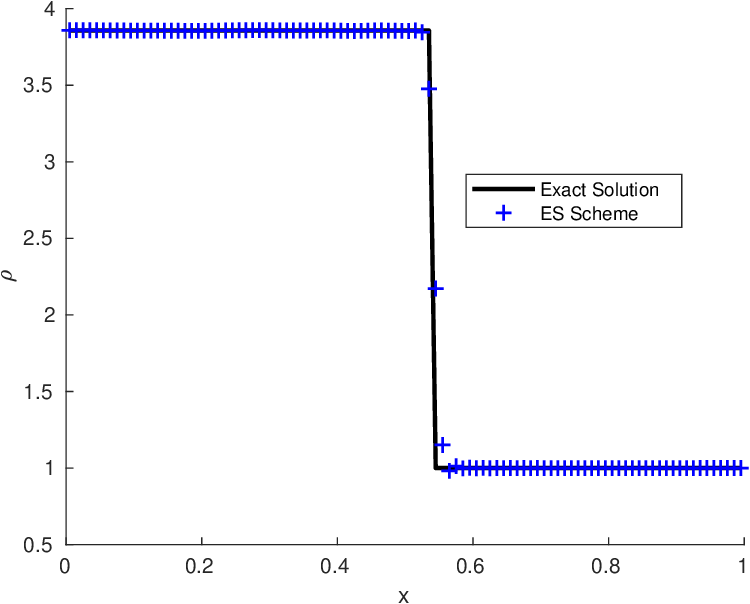}\quad
			\includegraphics[scale=0.395]{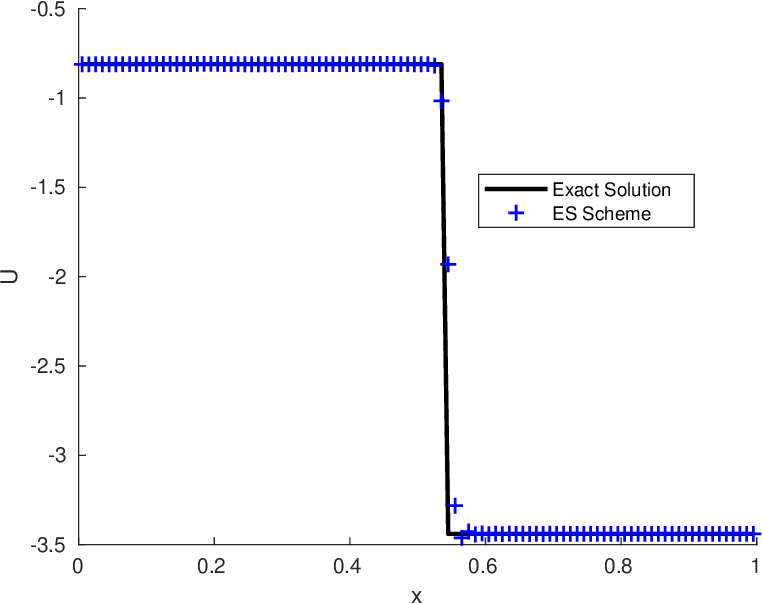}\quad
			\includegraphics[scale=0.395]{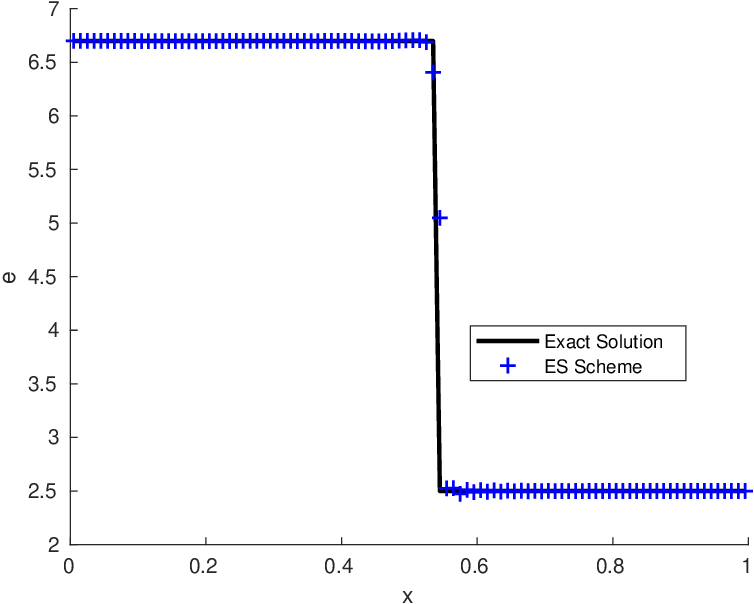}\\
		\end{minipage}
		\begin{minipage}{1\textwidth}
			\centering
			\includegraphics[scale=0.395]{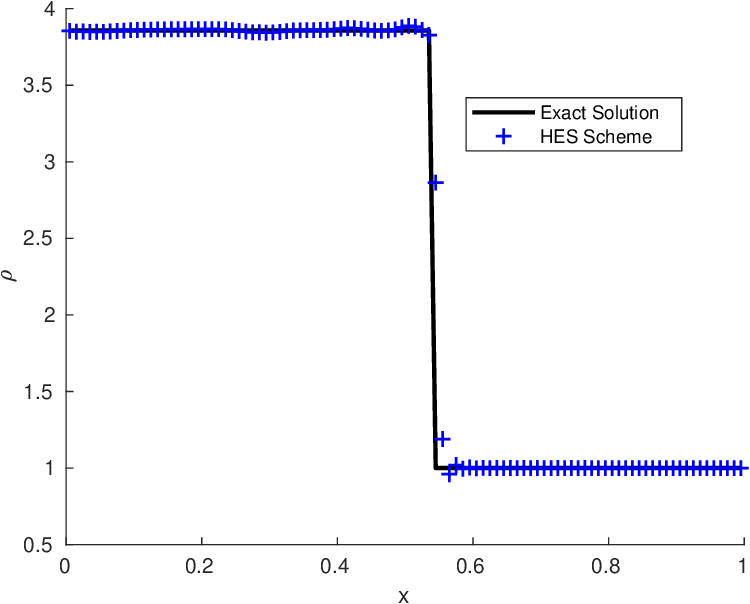}\quad
			\includegraphics[scale=0.395]{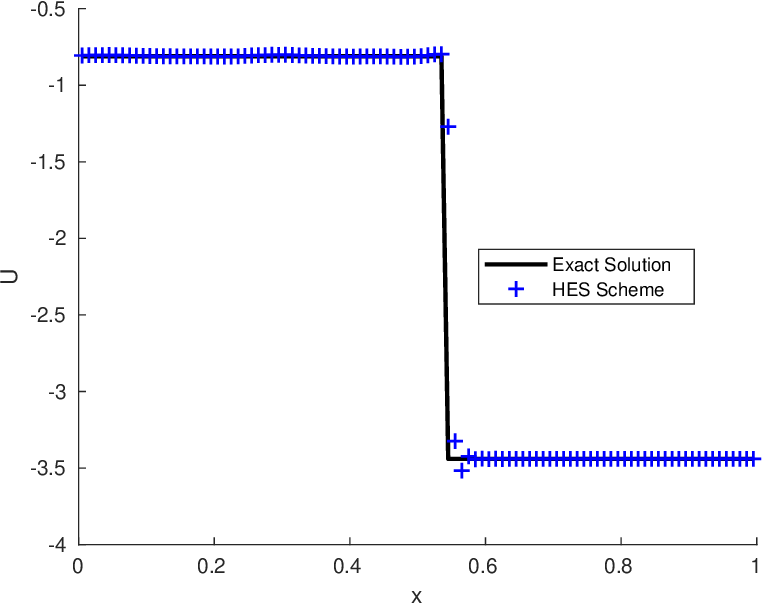}\quad
			\includegraphics[scale=0.395]{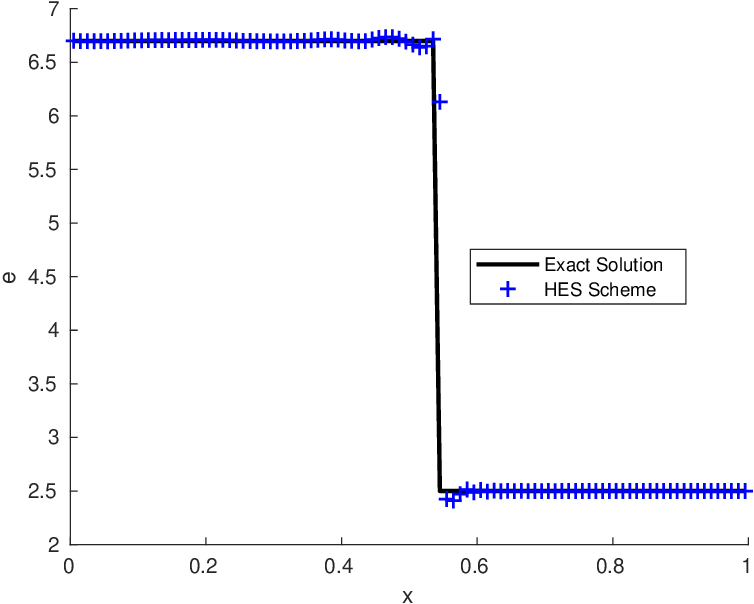}\\
		\end{minipage}
		\caption{Plots of density, velocity and internal energy using ES (top) and HES (bottom) schemes slowly moving shock wave.}
		\label{One_D_8}
	\end{figure}
	\begin{figure}
		\begin{minipage}{1\textwidth}
			\centering
			\includegraphics[scale=0.39]{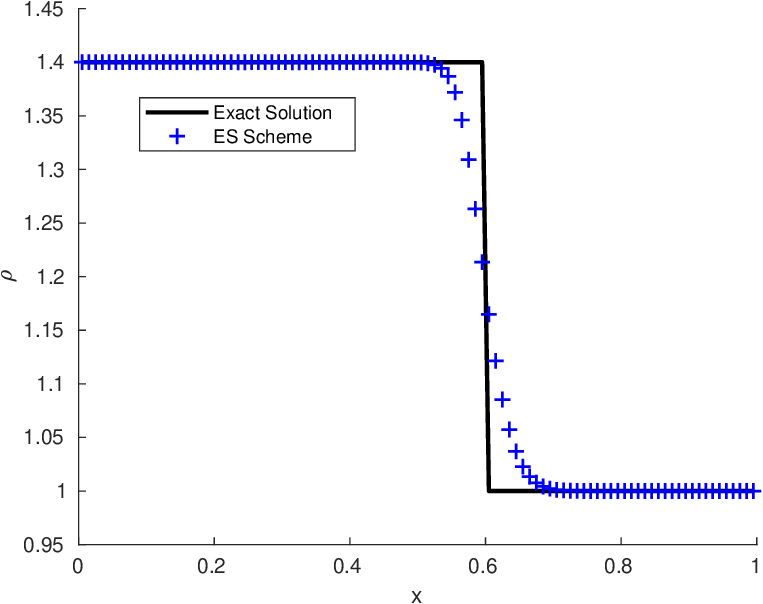}\quad
			\includegraphics[scale=0.39]{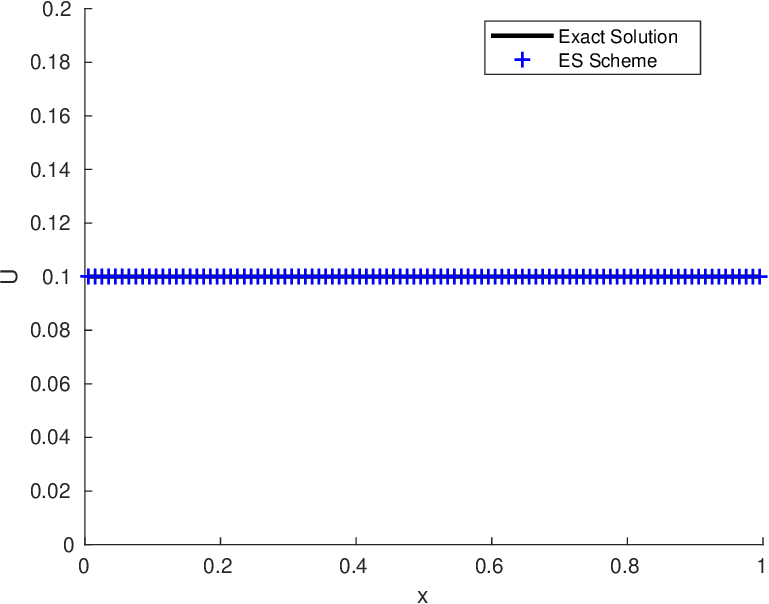}\quad
			\includegraphics[scale=0.39]{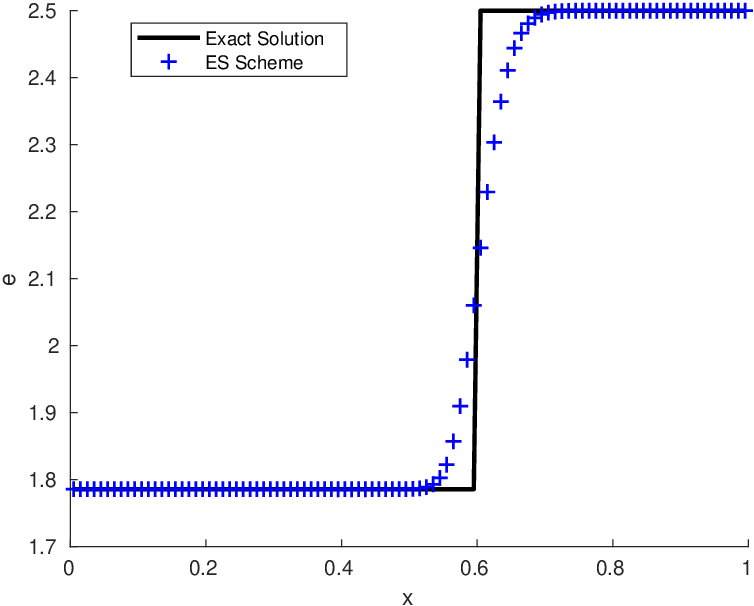}\\
		\end{minipage}
		\begin{minipage}{1\textwidth}
			\centering
			\includegraphics[scale=0.39]{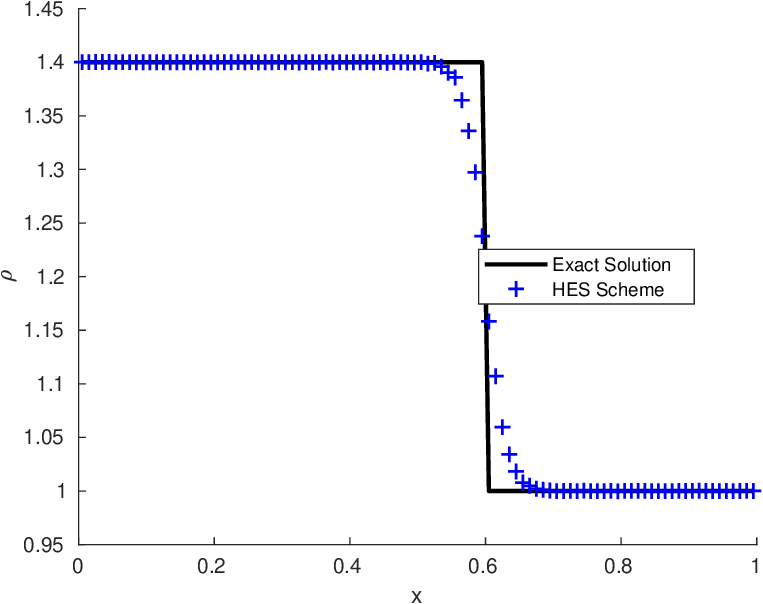}\quad
			\includegraphics[scale=0.39]{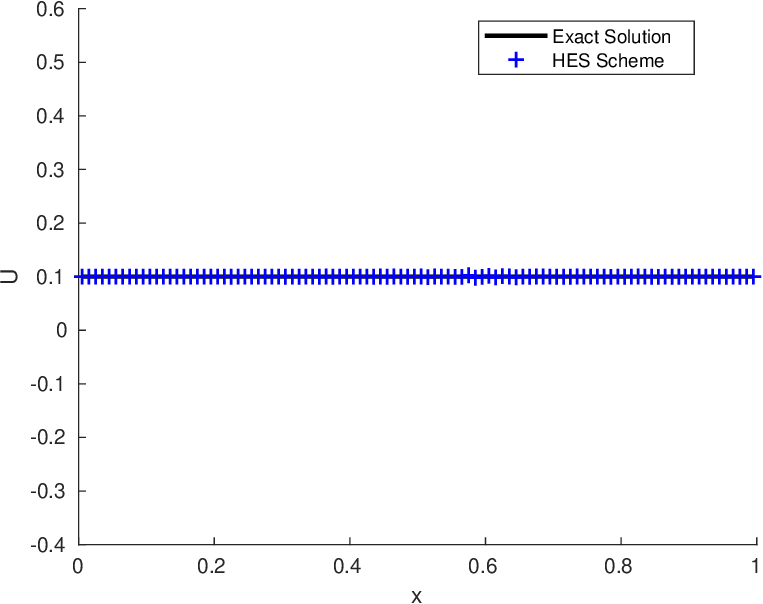}\quad
			\includegraphics[scale=0.39]{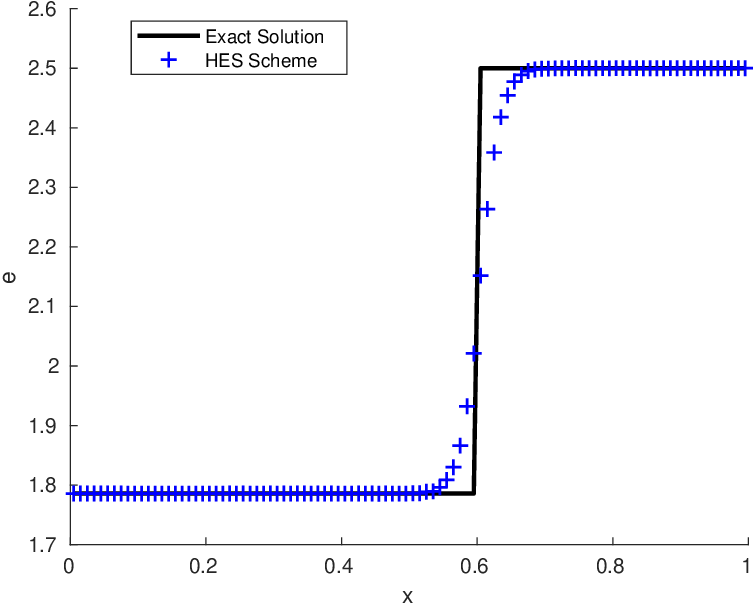}\\
		\end{minipage}
		\caption{Plots of density, velocity and internal energy using ES (top) and HES (bottom) schemes slowly moving contact wave.}
		\label{One_D_9}
	\end{figure}
	\newpage
	
 \subsection{Two Dimensional Results}
	
	This section presents solutions to several two-dimensional problems using the ES scheme \eqref{ES_Scheme1} and HES scheme \eqref{es_final}. The diffusion operator \eqref{4tho_diff} is computed utilizing the formulation for unstructured grids outlined in \cite{jameson_2015}. These test cases serve to validate the solutions for both accuracy and robustness. A CFL number of 0.1 is utilized for all instances, and the time-step is determined using the equation:  
\begin{equation}
    \Delta t=\min_{\forall Cells}{\frac{CFL\times Area}{(|\mathbf{V}\cdot \mathbf{n}_{\eta} |+c)\Delta s_{\eta}+(|\mathbf{V}\cdot \mathbf{n}_{\zeta} |+c)\Delta s_{\zeta}}}
\end{equation}
Here, $\eta$ and $\zeta$ represent the grid coordinate directions, and the corresponding normals $\mathbf{n}_{\zeta}$ and $\mathbf{n}_{\eta}$ are computed by averaging the normals from opposite sides of the control volume. For more details, see chapter 6 of \cite{blazek_2001}.  All steady test cases are run until all residuals (density, x-momentum, y-momentum, and total energy) fell below $10^{-8}$. The value of $q$ in equations \ref{ss_1} is taken as 10, and $\epsilon$ in \ref{ss_2} is taken between 0.1 and 0.5.  
 \subsubsection{Oblique Shock Reflection}
	
\begin{figure}\textbf{}
		\begin{subfigure}{1\textwidth}
			\centering
			\includegraphics[width=.48\linewidth]{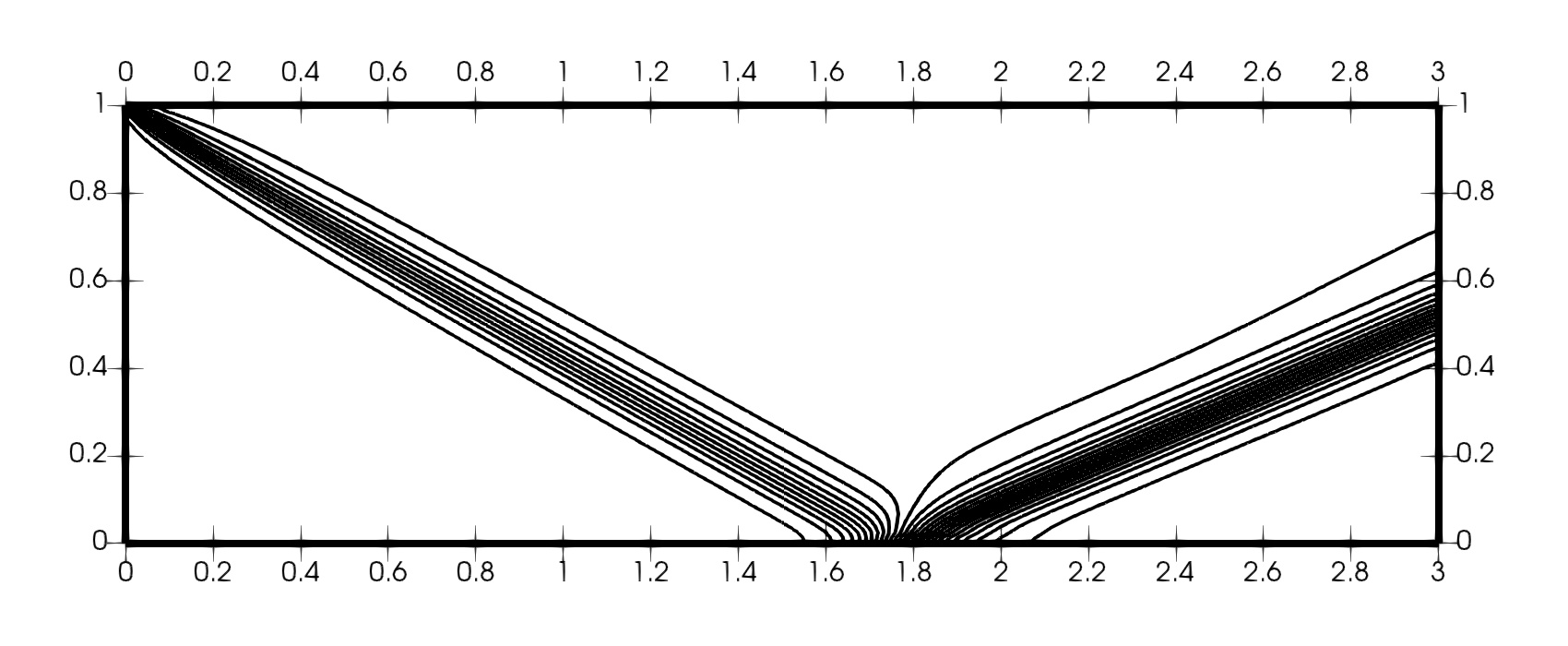}\quad
			\includegraphics[width=.48\linewidth]{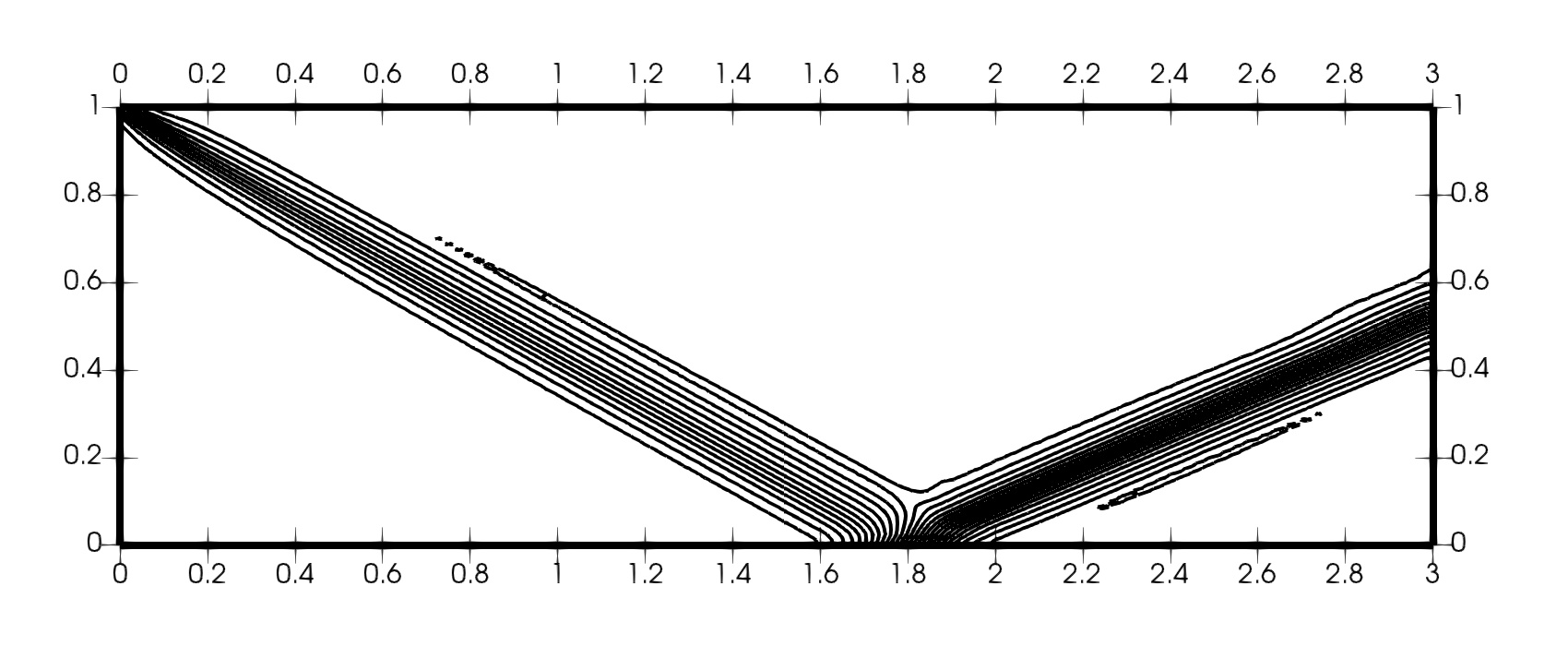}\quad
			\caption{Density contours for (30 in 0.9-2.7) for 240x80 grid with ES(left) and HES(right) Scheme}
		\end{subfigure}
		\begin{subfigure}{1\textwidth}
			\centering
			\includegraphics[width=.48\linewidth]{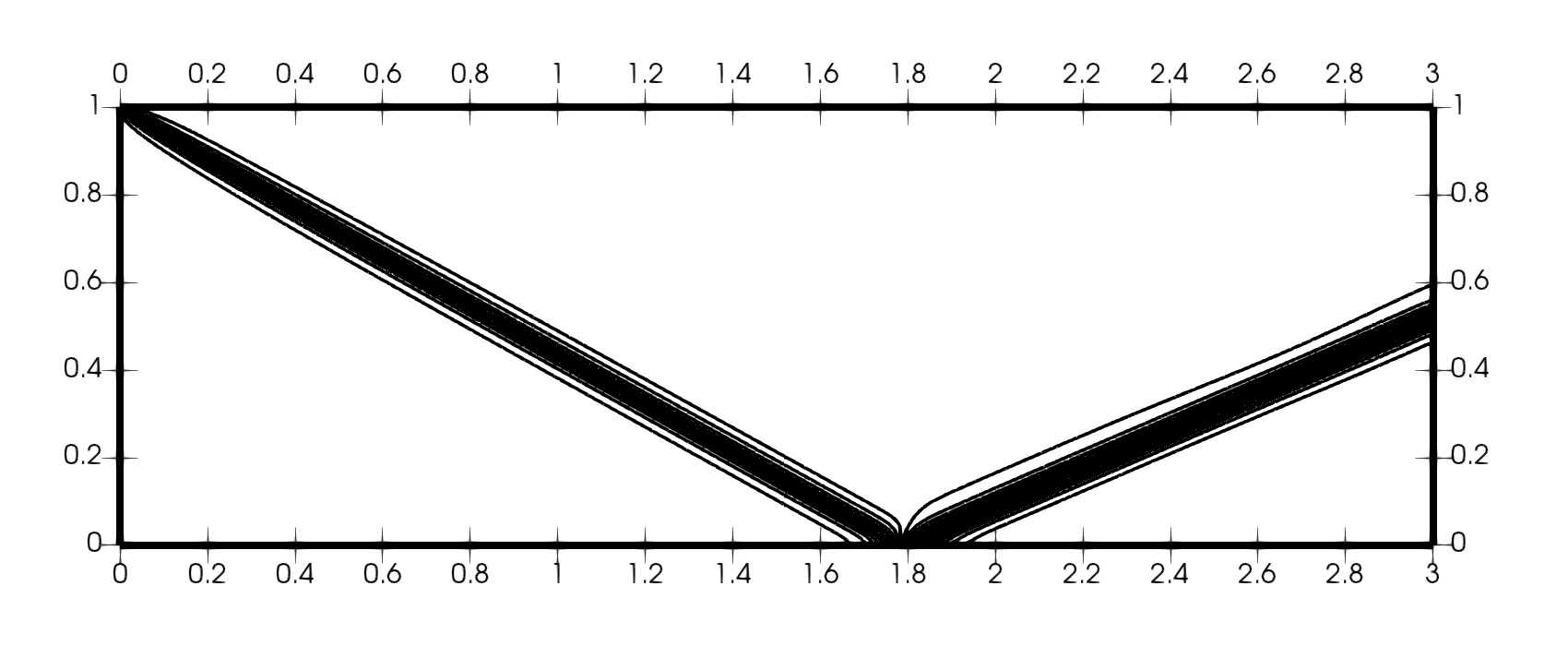}\quad
			\includegraphics[width=.48\linewidth]{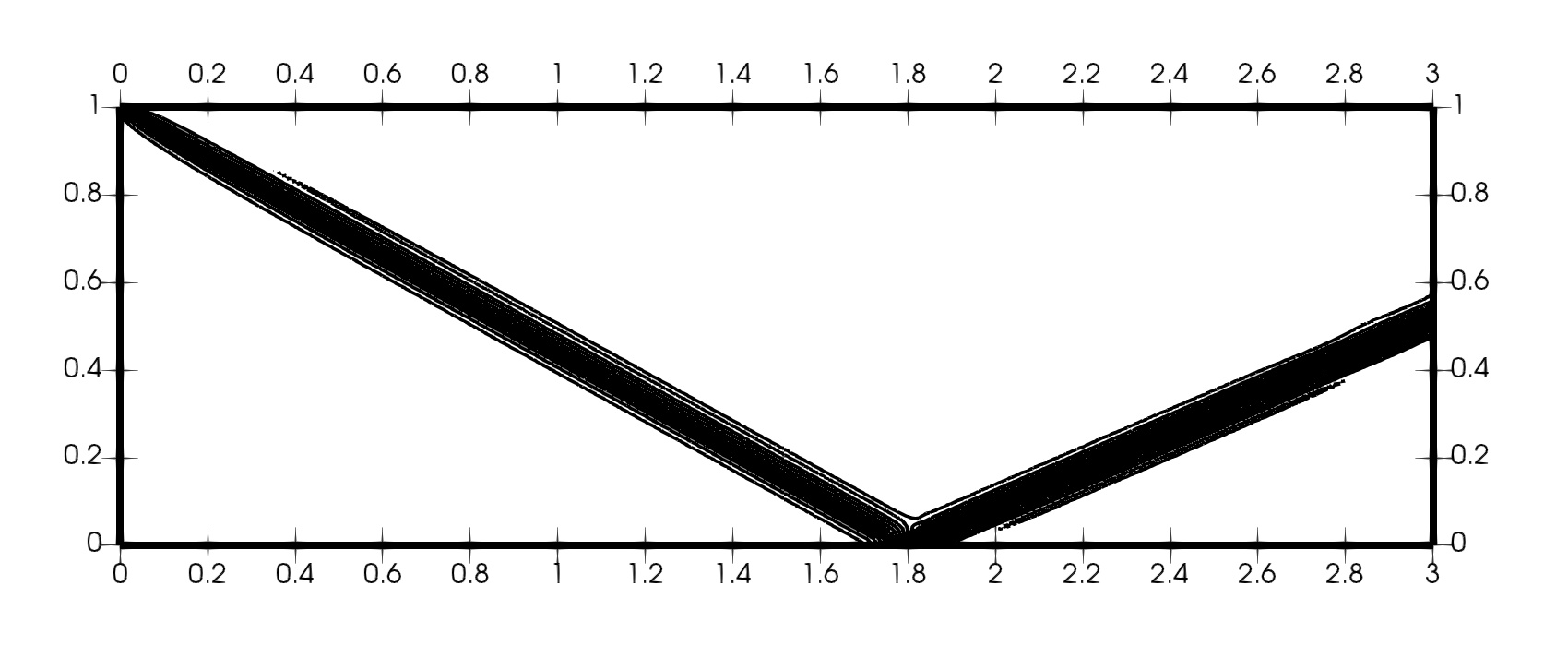}\quad
			\caption{Density contours for (30 in 0.9-2.7) for 480x160 grid with ES(left) and HES(right) Scheme}
		\end{subfigure}
		\caption{Oblique shock reflection test case}
		\label{obl_shock_ref}
	\end{figure}

\begin{figure}
		\begin{minipage}{1\textwidth}
			\centering
			\includegraphics[scale=0.4]{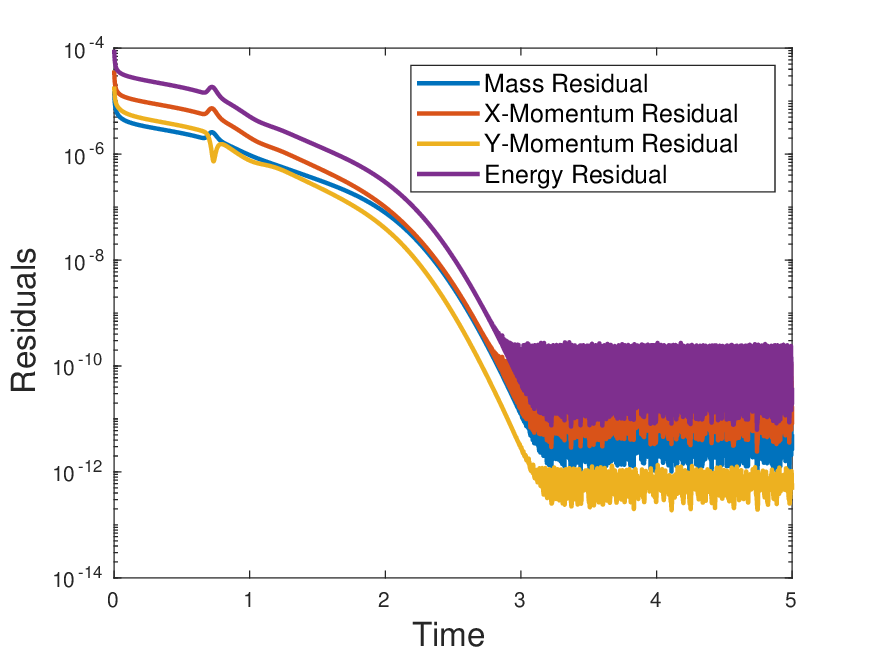}\quad
			\includegraphics[scale=0.4]{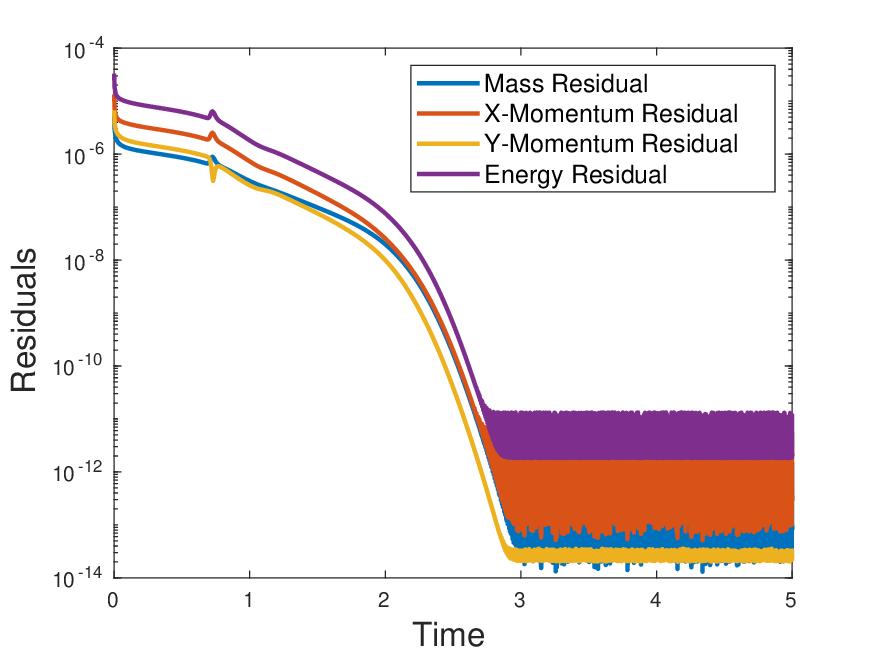}\quad
		\end{minipage}
		\caption{Semi-log plot of residuals for ES scheme for 240x80 grid (left) and 480x160 grid (right)}
		\label{obl_shock_ref_res}
	\end{figure}

\begin{figure}
            \begin{minipage}{1\textwidth}
			\centering
			\includegraphics[scale=0.4]{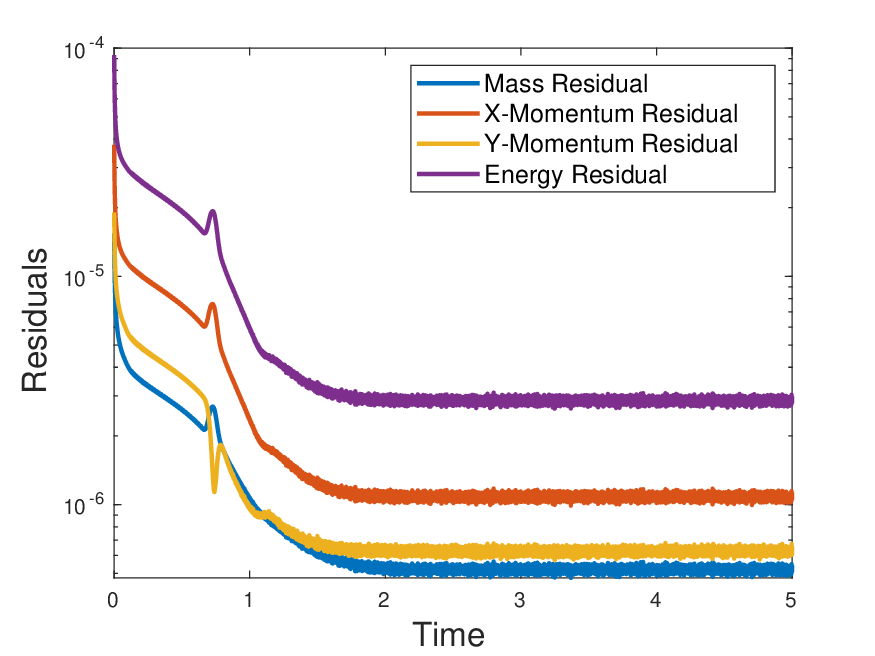}\quad
			\includegraphics[scale=0.4]{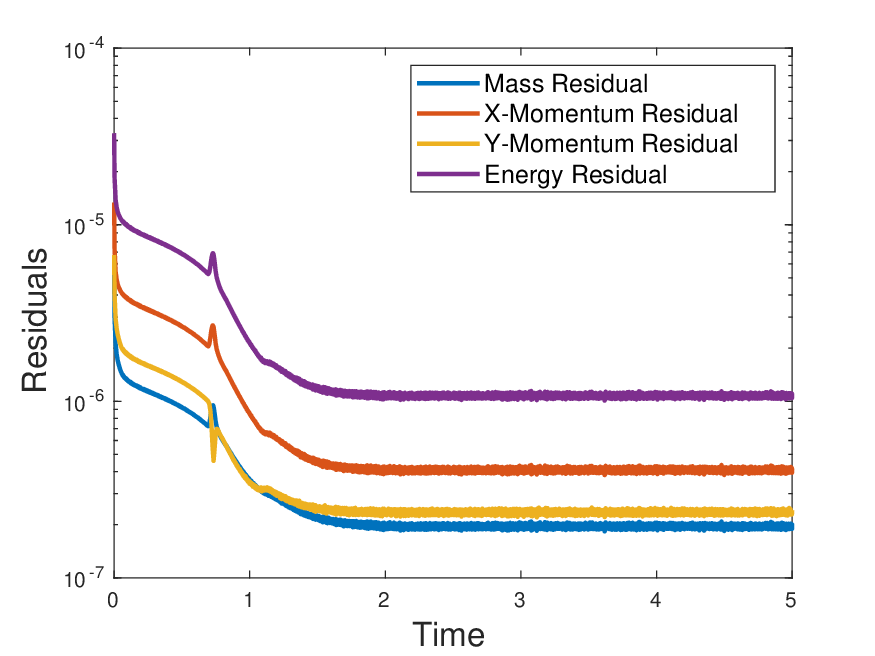}\quad
		\end{minipage}
		\caption{Semi-log plot of residuals for HES scheme for 240x160 grid (left) and 480x160 grid (right)}
		\label{obl_shock_ref_res_2}
	\end{figure}

As described by Yee \cite{yee_1982}, this test case involves an oblique shock encountering a wall, resulting in a reflected shock wave. The free stream Mach number ($M_{\infty}$) is 2.9, and the incident shock angle is $29^o$. The domain is discretized with Cartesian mesh cells in the $[0,3] \times [0,1]$ region. Supersonic inflow boundary conditions are applied on the left boundary with inflow values ($\mathbf{W}=[1.0,2.9,0,1.0/1.4]$)  and post-shock values are prescribed at the top boundary. A flow tangency boundary condition is used at the bottom boundary (solid wall), while a supersonic outflow boundary condition is applied on the right boundary. The flow field is initialized with the inflow conditions. Simulations are conducted with grid sizes of 240x80 and 480x160, and density contours are depicted in Figure \ref{obl_shock_ref} for both schemes. As illustrated, the incident and reflected shocks are accurately resolved, with the shocks being less diffusive for higher-order schemes. Semi-log plots of density, momentum, and energy residuals versus the number of iterations are displayed in Figures \ref{obl_shock_ref_res} and \ref{obl_shock_ref_res_2}, demonstrating good convergence to steady-state solutions for both schemes.   

\subsubsection{Supersonic flow over a compression ramp}
\begin{figure}
		\begin{minipage}{1\textwidth}
			\centering
			\includegraphics[width=.48\linewidth]{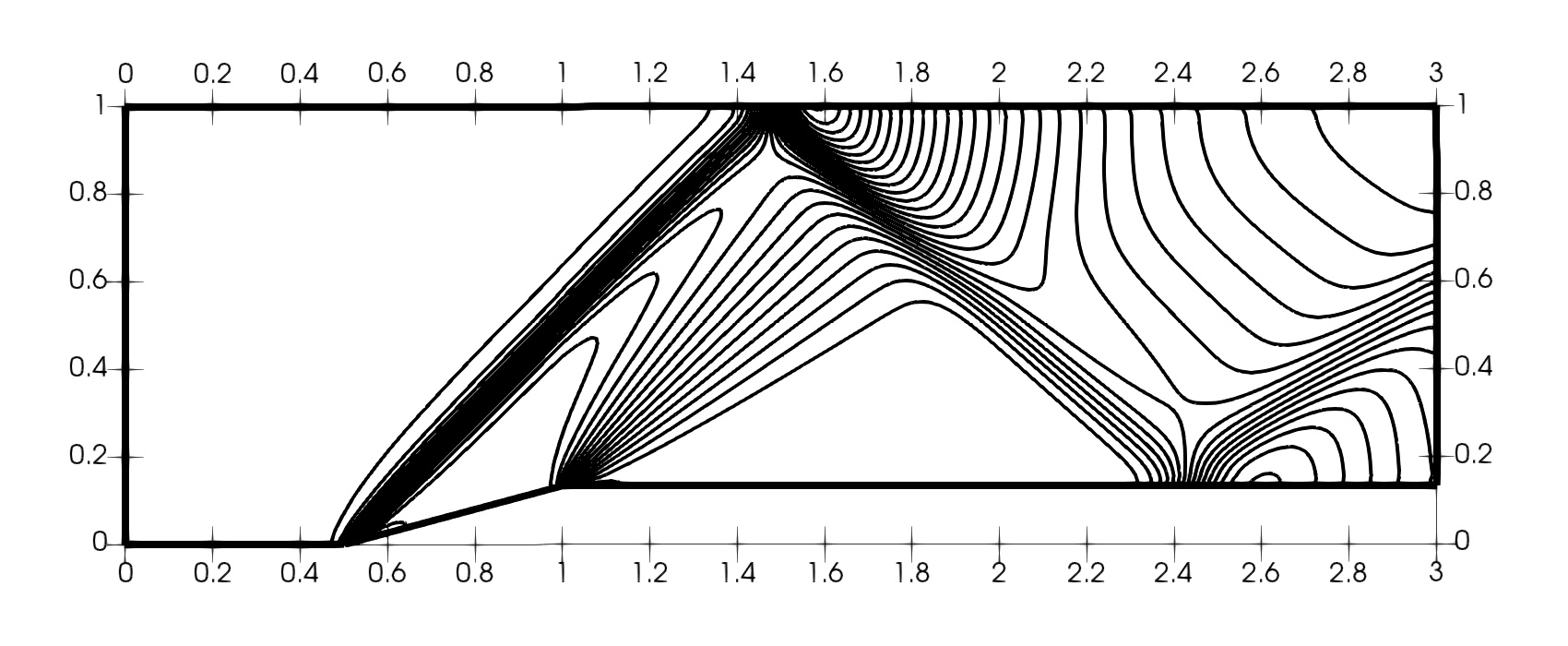}\quad
			\includegraphics[width=.48\linewidth]{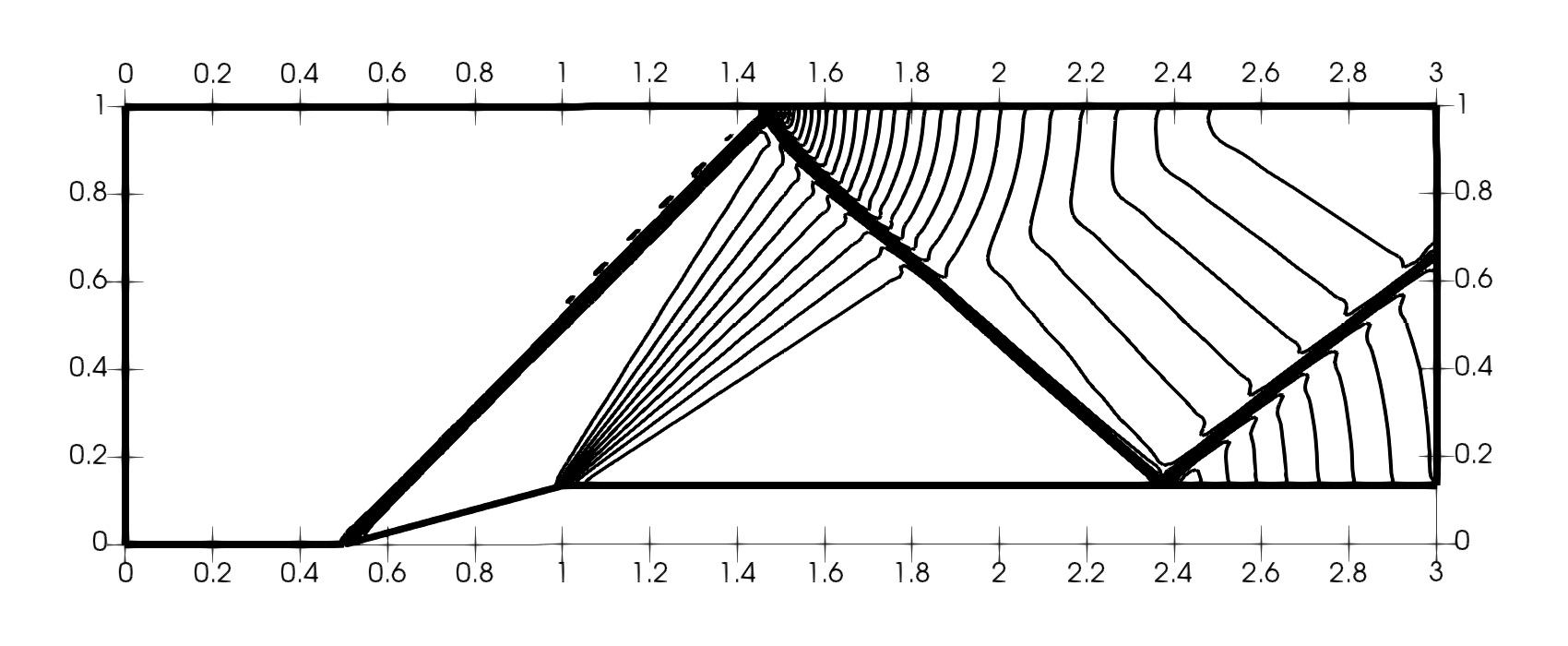}\quad \\
		\end{minipage}
		\caption{Pressure contours (50 in 0.6-2.4) for $15^o$ compression ramp with ES (left) and HES (right) schemes}
		\label{comp_ramp}
	\end{figure}

This test case, as studied by Levy et al. \cite{levy_1993}, involves a supersonic flow of Mach 2 over a compression ramp of $15^{0}$.  The domain for this problem is $[0,3]\times[0,1]$. Supersonic inflow boundary conditions are imposed on the left boundary with primitive variables $\mathbf{W}=[1.0,2,0,1.0/1.4]$. Flow tangency boundary conditions are applied on the top and bottom boundaries, while a supersonic outflow boundary condition is imposed on the right boundary. The entire domain is initialized with the inflow conditions. This steady problem exhibits an initial shock emanating from a ramp, further getting reflected from the top and bottom boundaries. An expansion fan forms at the end of the ramp and further interacts with the reflected shock. Figure \ref{comp_ramp} shows pressure contours, capturing the incident shock, reflected shock, an expansion fan and wave interactions with both ES and HES schemes.  
 
\subsubsection{Hypersonic flow over a half-cylinder} \label{half-cylinder}
\begin{figure}
		\begin{subfigure}{.3\textwidth}
			\centering
			\includegraphics[width=.48\linewidth]{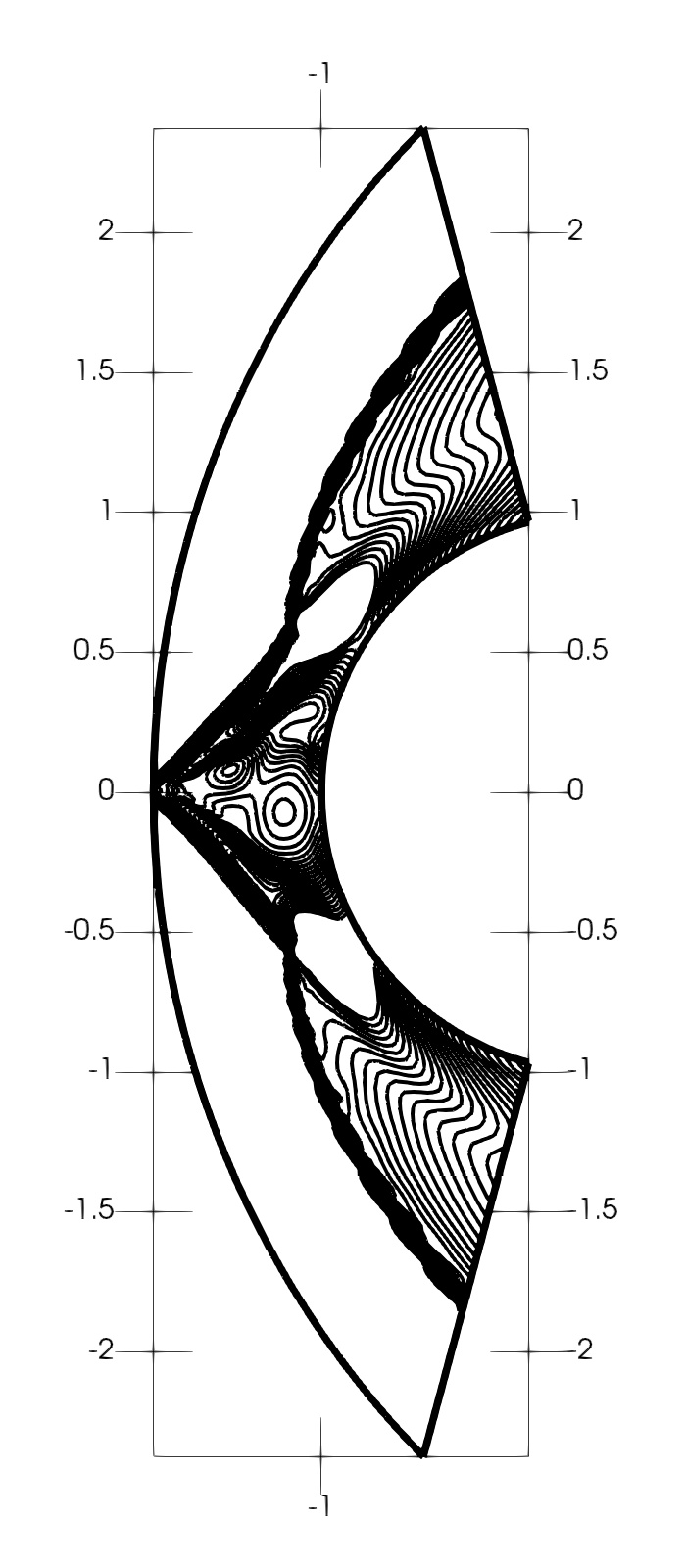}\quad
			\caption{}
			\label{half_cly_roe}
		\end{subfigure}
		\begin{subfigure}{.3\textwidth}
			\centering
			\includegraphics[width=.48\linewidth]{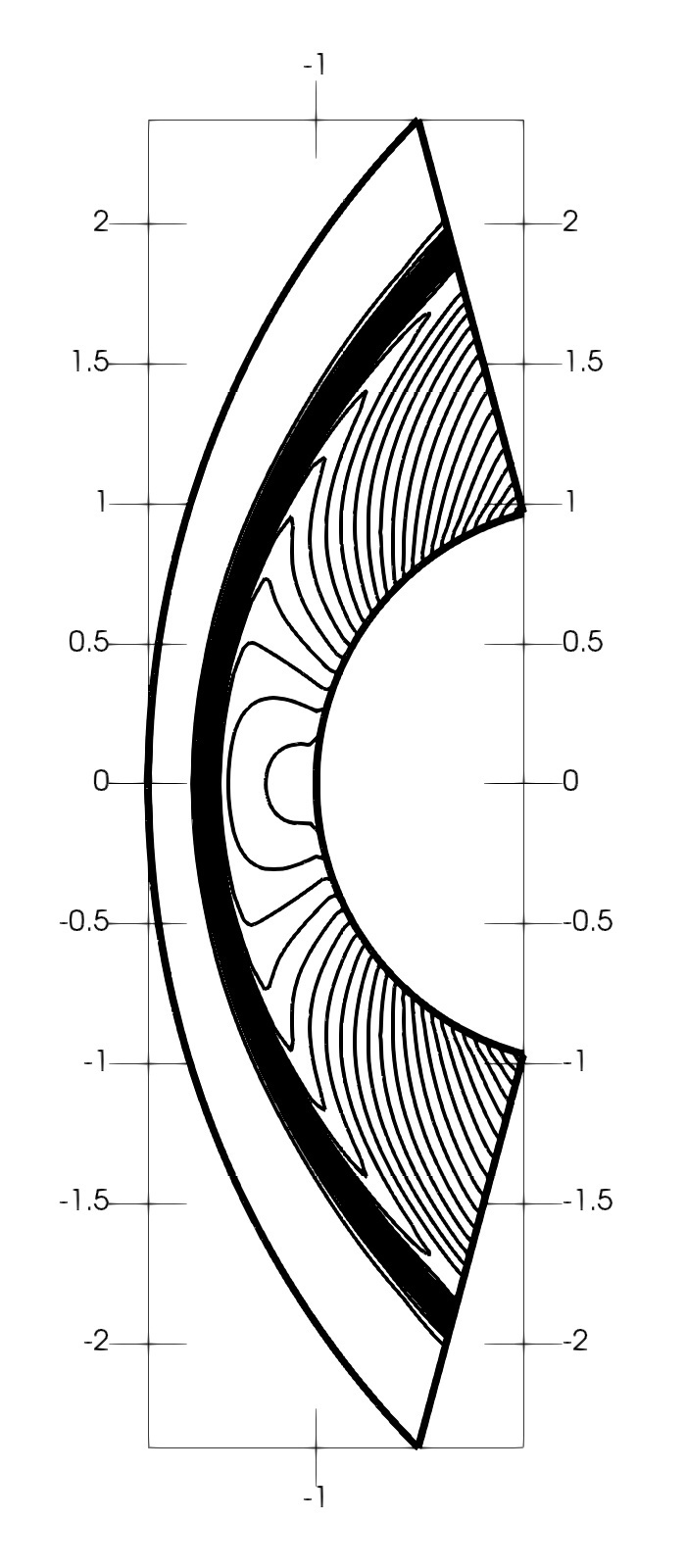}\quad
			\caption{}
			\label{half_cly_es}
		\end{subfigure}
		\begin{subfigure}{.3\textwidth}
			\centering
			\includegraphics[width=.48\linewidth]{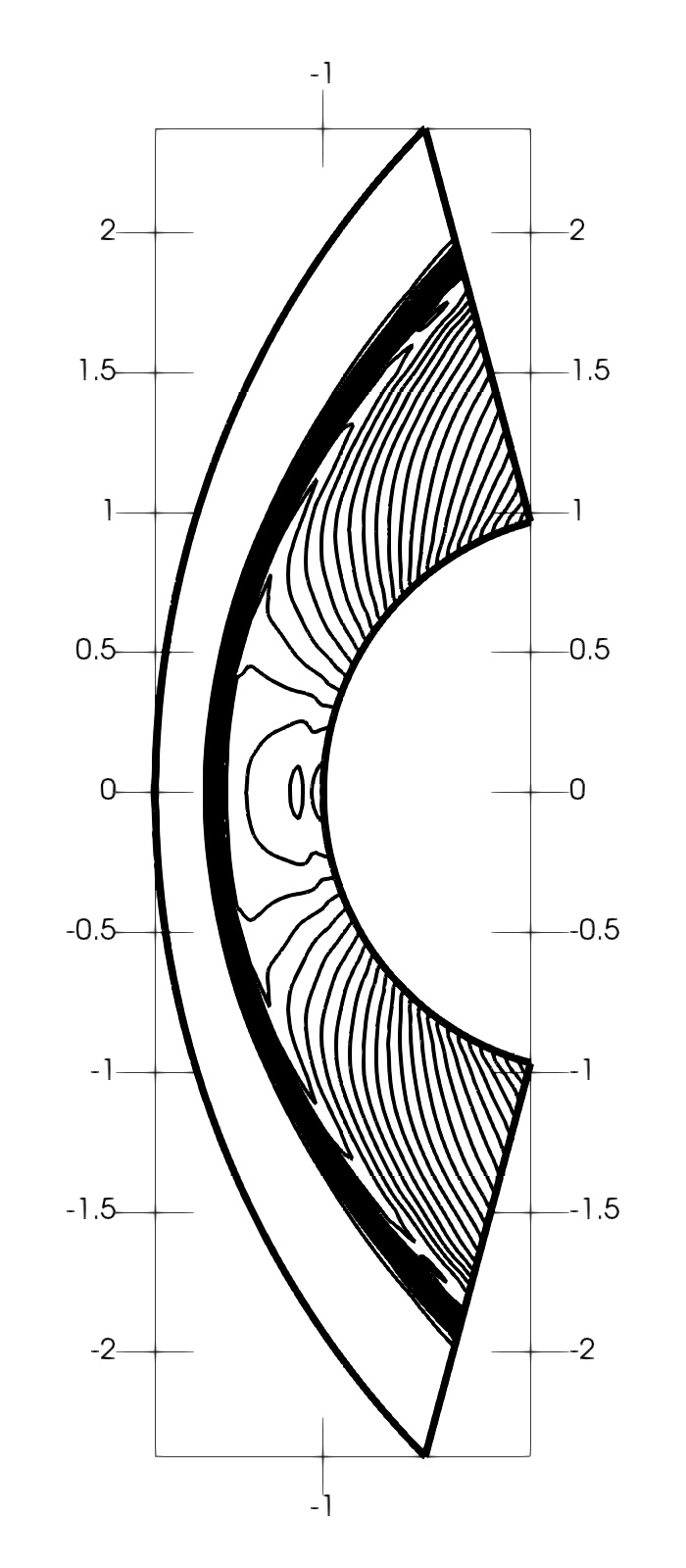}\quad
			\caption{}
			\label{half_cly_hes}
		\end{subfigure}
		\caption{Density contours (30 in 1.3-9.0) for hypersonic flow over cylinder on 320x40 grid with (a) Roe scheme, (b) ES scheme and (c) HES scheme }
		\label{half_cly}
	\end{figure}

	This test case involves hypersonic flow at Mach 20 over a half-cylindrical body, resulting in a pronounced bow shock. An undesirable phenomenon known as a {\em carbuncle shock} \cite{peery_1988} is observed in low-diffusive numerical schemes like Roe scheme, with an unphysical perturbation disrupting the bow shock on the stagnation line.  The problem becomes more severe when the mesh is aligned with the bow shock. This undesirable feature, which often occurs with Riemann solvers, has attracted considerable research and is counted together with several shock instabilities that plague the Riemann solvers \cite{quirk_1994}.  For generating the mesh, the mesh curvature is determined using the formulation by Huang et al. \cite{huang_2011}, given by the equations:
\begin{equation}
\begin{aligned}
x=&\frac{2 x_c-\sqrt{4x_c^2-4(1+\tan(\theta)^2)(x_c^2-r_i^2))}}{2(1+\tan(\theta)^2)} \\
y=&- \tan(\theta)x
\end{aligned}
\end{equation}
where $\theta=(j-1)*5\pi/(6M)-5\pi/12$, $x_c=1.8(N-i+1)/N$, and $r_i=1+2.4(N-i+1)/N$. Here, $N$ and $M$ are chosen as 40 and 320, respectively. The initial and inflow conditions (on the left boundary) are $\mathbf{W}=[1.0,20,0,1.0/1.4]$. A flow tangency boundary condition is imposed on the cylinder surface, while a supersonic outflow boundary condition is applied on the remaining periphery. Figure \ref{half_cly} depicts density contours obtained using different schemes. As illustrated in Figure \ref{half_cly_roe}, the Roe scheme exhibits a carbuncle shock, disrupting the shock structure, while the ES and HES schemes, shown in Figures \ref{half_cly_es} and \ref{half_cly_hes}, avoid this artifact.  
 
 \subsubsection{Forward-facing step in supersonic flow}
\begin{figure}
		\centering
		\begin{subfigure}{.48\textwidth}
			\includegraphics[width=\linewidth]{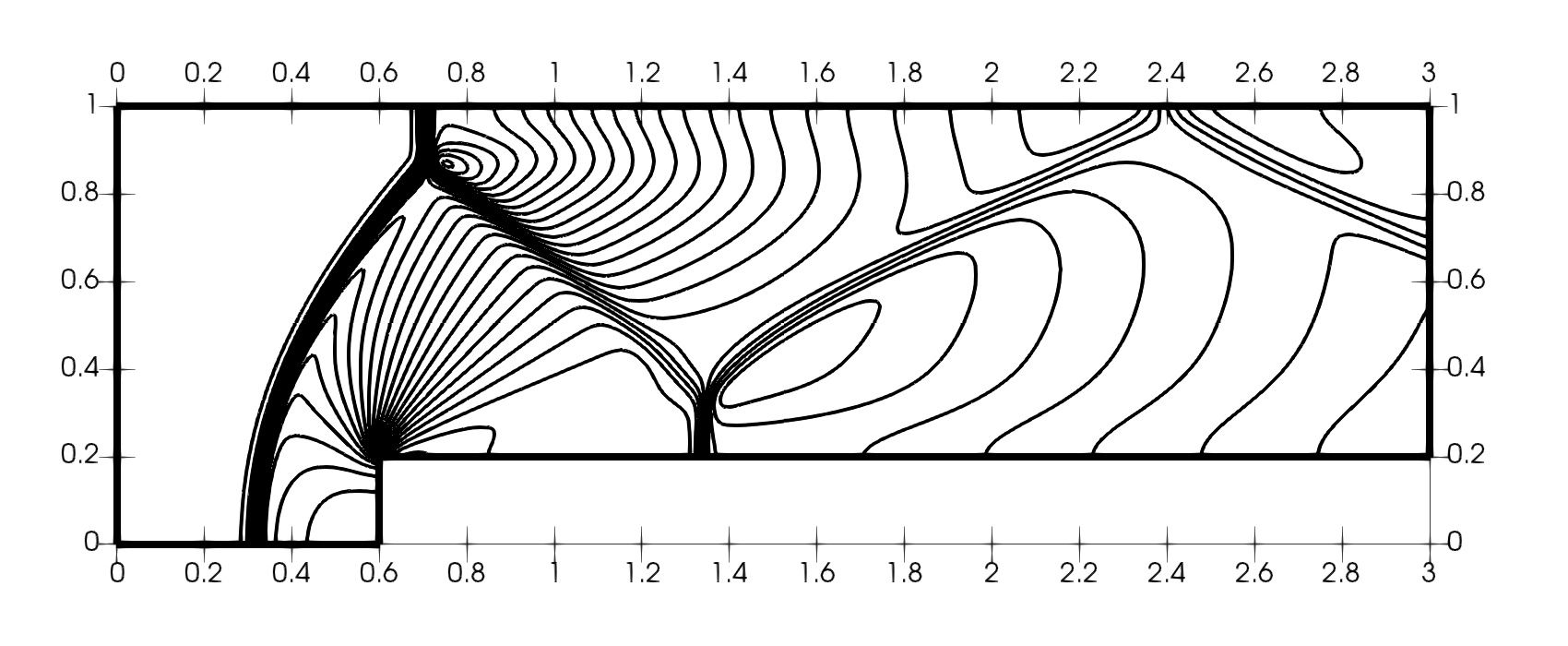}
		\end{subfigure}
		\begin{subfigure}{.48\textwidth}
			\includegraphics[width=\linewidth]{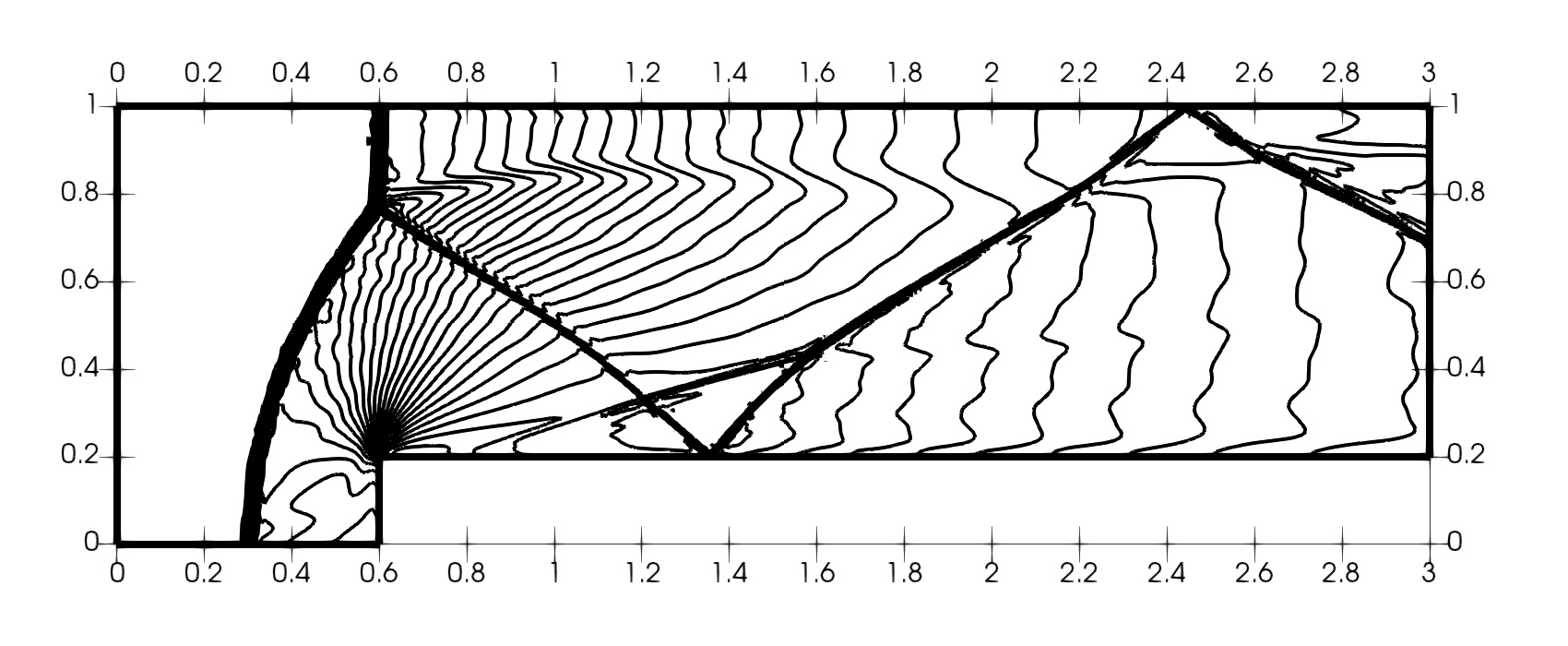}
		\end{subfigure}
		\caption{Density contours (30 in 0.3-6.8) for forward step on 480x160 grid using ES (left) and HES (right) scheme}
		\label{fwd_step}
	\end{figure}

	This test case, introduced by Emery \cite{emery_1968}, comprises a step facing a supersonic flow of Mach 3.   Domain is taken as $[0,3]\times[0,1]$ and inflow and initial conditions are: $\mathbf{W}=[1.0,3.0,0.0,1.0/1.4]$. The height of the step is 0.2 units and is located at a distance of 0.6 from the left boundary. Top and bottom boundaries are solid walls, with flow tangency boundary conditions being applied. At the left boundary supersonic inflow boundary condition is applied, and at the right boundary supersonic outflow boundary condition is enforced. This is an unsteady test case; the final solution is taken at $t=4\, s$, at which point the shocks move very slowly. Flow features include a bow shock and its reflection from top and bottom boundaries. Also, a powerful expansion fan is formed at the tip of the step, which interacts with the reflected shocks. A slipstream that travels downstream is also formed at the $\lambda$-shock (at the triple-point), approximately at $y=0.8$. The corner is a singularity in the flow field, and it induces a numerical boundary layer that ruins the solution downstream. "Corner fix" given in \cite{woodward_1984} is used to fix this issue. Results are shown in figure \ref{fwd_step}, and both the schemes can capture the discontinuities and wave interactions with good accuracy, with the HES being better than the ES scheme. No expansion shock is seens with either of the schemes in the region near the corner. For the ES scheme, the Mach stem is observed on the lower wall in the final solution. However, no Mach stem is observed for the HES scheme, and the Mach stem at the top wall is also captured at the correct location ($x=0.6$). Slipstream at $y=0.8$ just after the Mach stem on the top surface is also clearly captured, which is a feature associated with low diffusion schemes. 

    \subsubsection{Shock diffraction problem} 
\begin{figure}
		\centering
		\begin{subfigure}{.24\textwidth}
			\includegraphics[width=\linewidth]{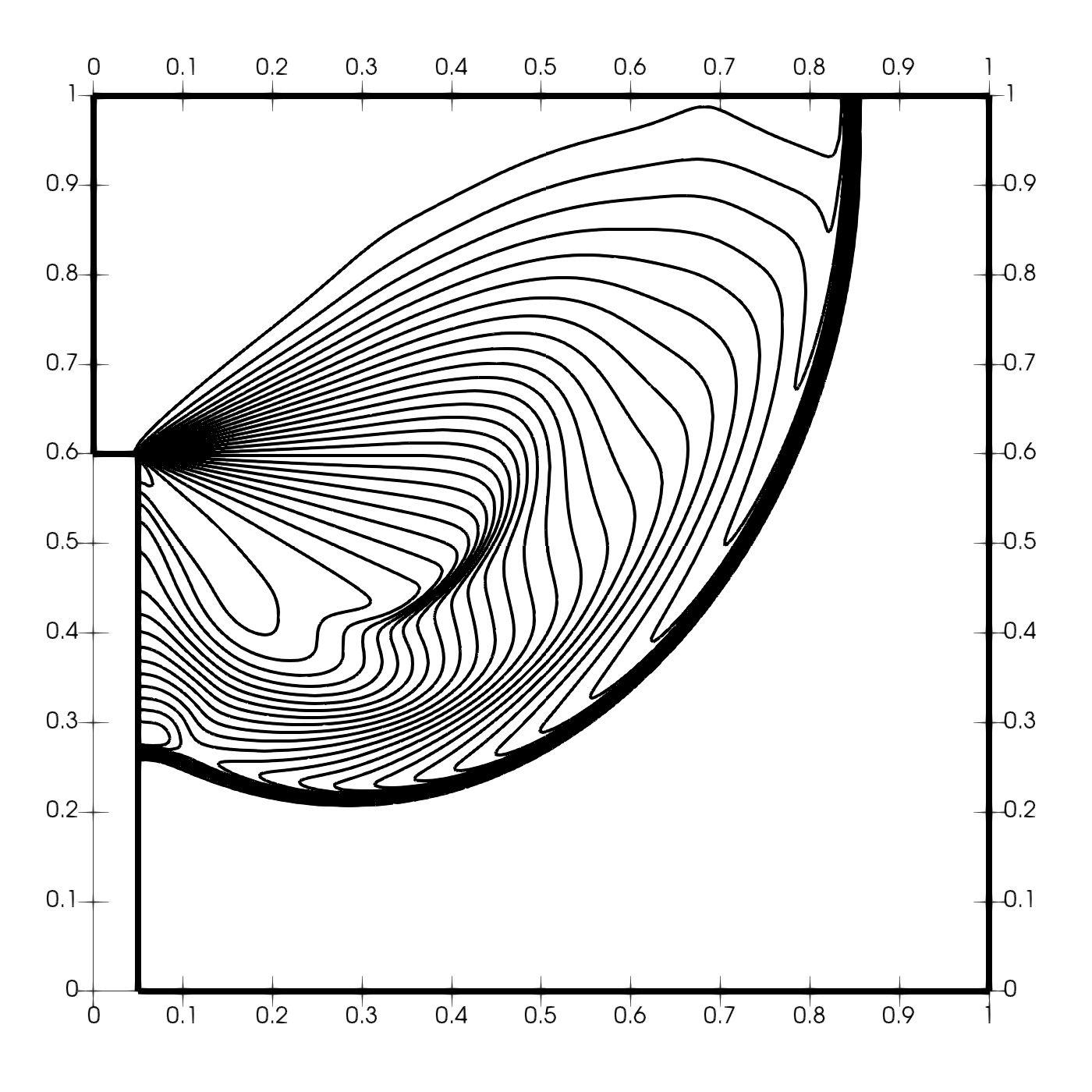}
			\caption{}
		\end{subfigure}
		\begin{subfigure}{.24\textwidth}
			\includegraphics[width=\linewidth]{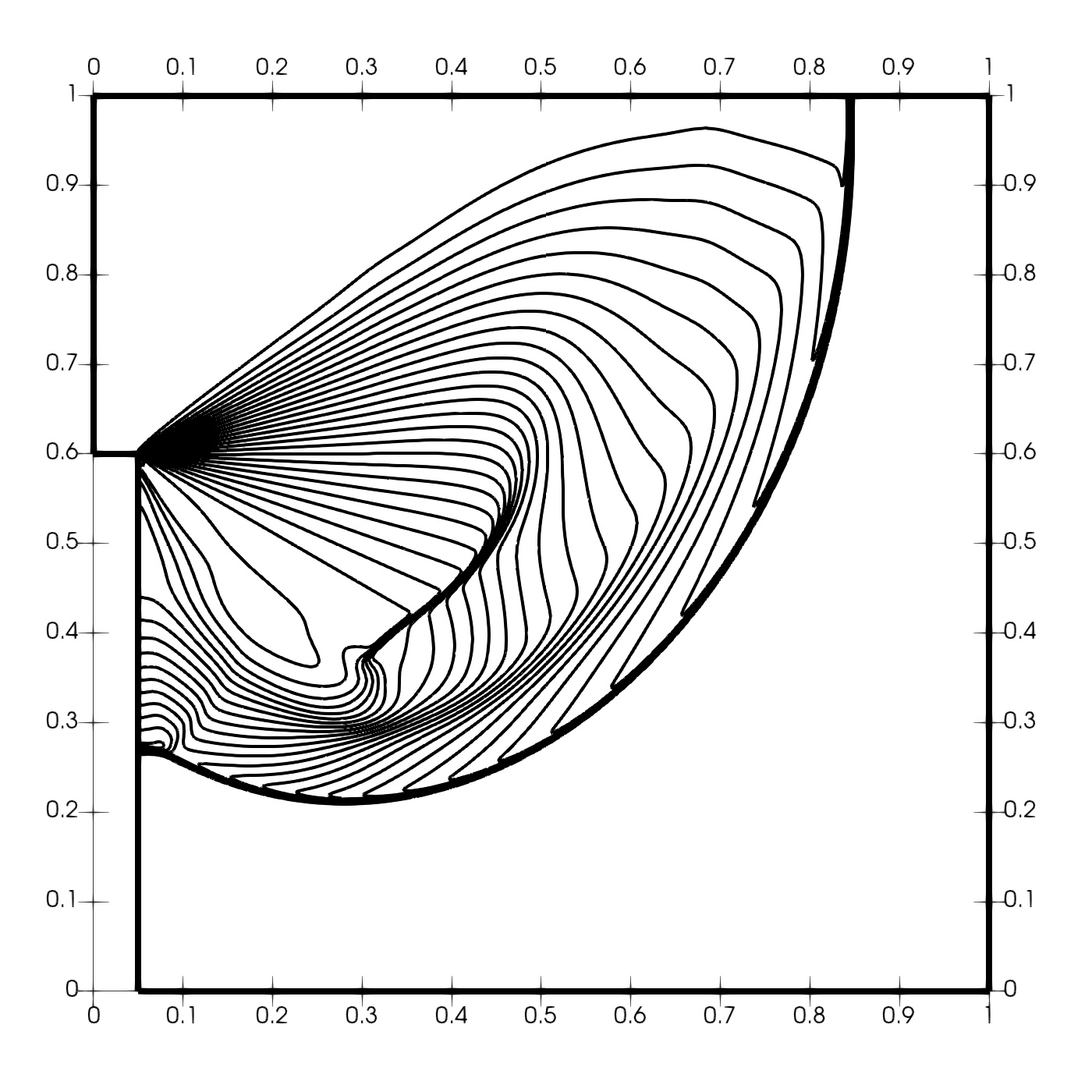}
			\caption{}
		\end{subfigure}	
		\begin{subfigure}{.24\textwidth}	
			\includegraphics[width=\linewidth]{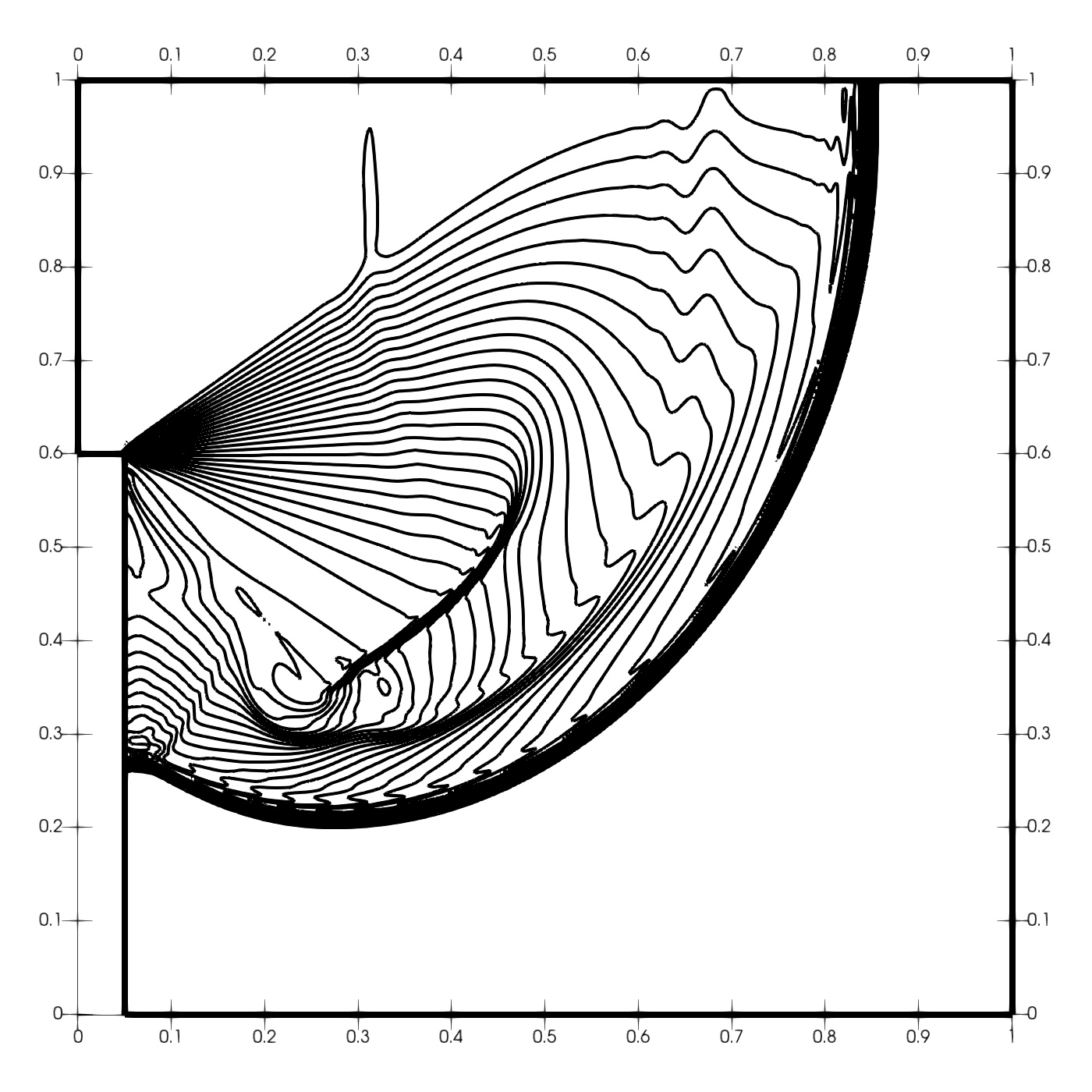}
			\caption{}
		\end{subfigure}	
		\begin{subfigure}{.24\textwidth}
			\includegraphics[width=\linewidth]{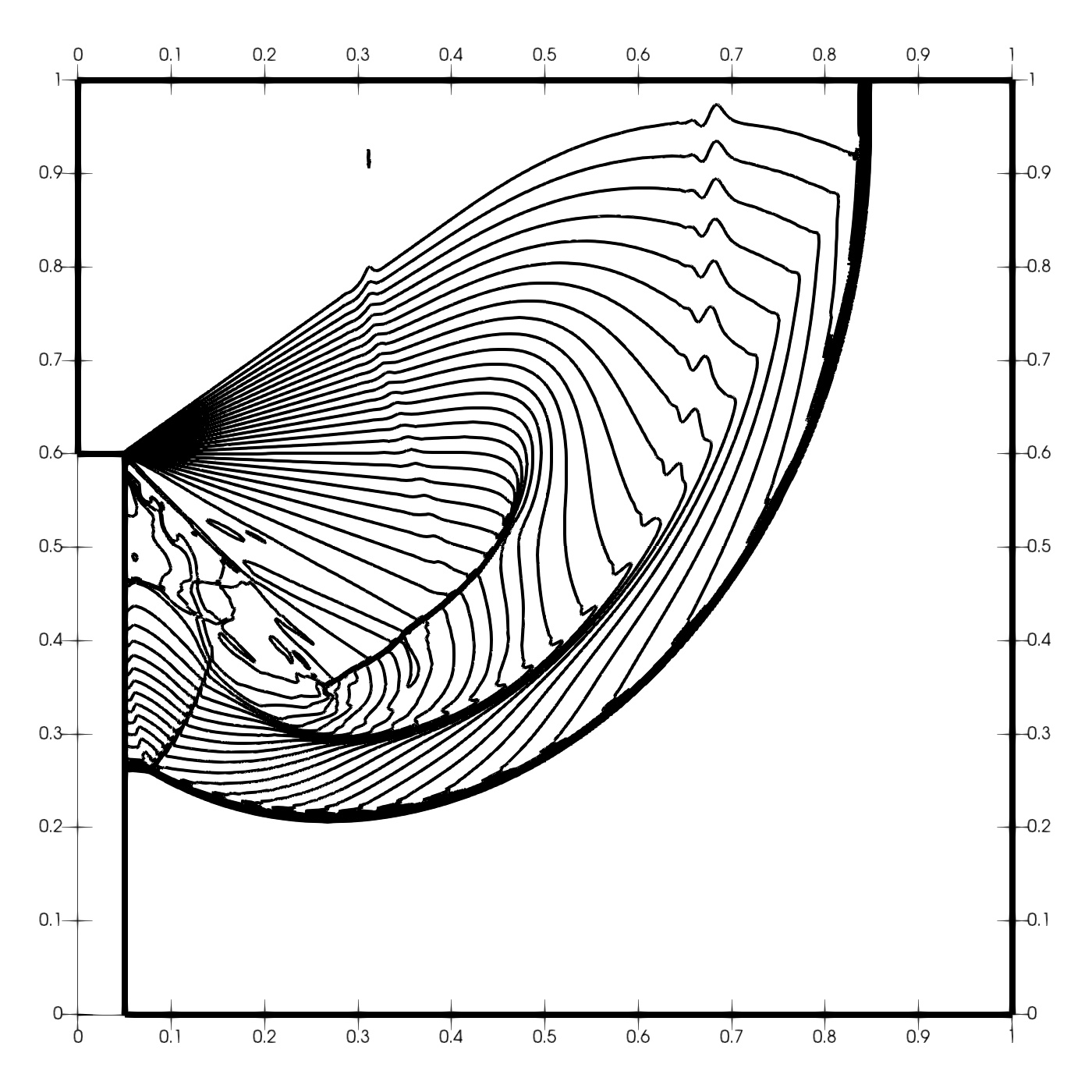}
			\caption{}
		\end{subfigure}
		\caption{Density contours (30 in 0.05-7.1) for backwards-facing step. (a) and (b) using the ES scheme on 400x400 and 1200x1200 grids, respectively. (c) and (d) using the HES scheme on 400x400 and 1200x1200 grids, respectively}
		\label{back_step}
	\end{figure}

	This test case involves a Mach 5.09 flow diffracting over a $90^{o}$ corner of a backward-facing step, forming a strong expansion wave. Density and pressure after the expansion are minimal, and thus, schemes that fail to preserve positivity fail in this test case, such as Roe's approximate Riemann solver \cite{quirk_1994}.  Domain for this test case is taken as $[0,1]\times[0,1]$ with the corner of length 0.05 units located at $y=0.6$. The domain is initialized with $\mathbf{W}=[1.4,0,0,1.0]$ for $x>0.05$ and the post-shock conditions for $x<0.05$. Supersonic inflow boundary condition is imposed on the left side from $y=0.6$ to $1.0$ and flow tangecy boundary condition is imposed on the surface of the step. Symmetry boundary condition is imposed on all the other boundaries. The solution is computed at a final time of $T=0.1561s$. Numerous schemes that suffer from shock instabilities do not preserve the shock structure in this test case, and an unphysical expansion shock often appears at the corner. Also, the rarefaction from the corner creates near-zero density and pressures, and schemes fail to preserve their positivity. Results from the ES and HES schemes can be seen in figure \ref{back_step} for different grids. The planar shock is preserved, and the entropy stable schemes do not produce expansion shocks or any other anomalies.    
    \subsubsection{Odd-even decoupling}
\begin{figure}
		\begin{subfigure}{.3\textwidth}
			\centering
			\includegraphics[width=\linewidth]{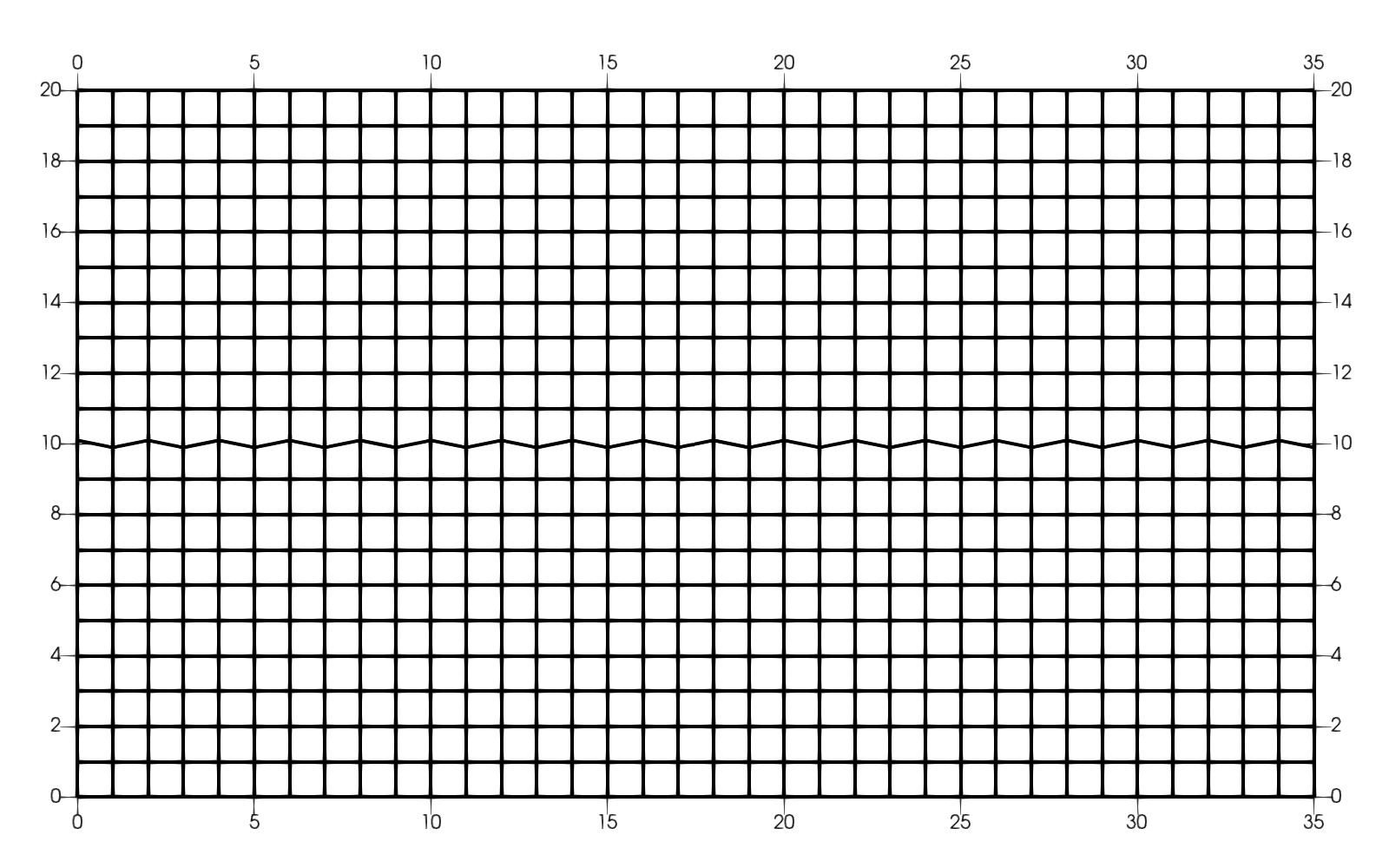}\quad
			\caption{}
			\label{odd_even_grid}
		\end{subfigure}
			\begin{subfigure}{.3\textwidth}
			\centering
			\includegraphics[width=\linewidth]{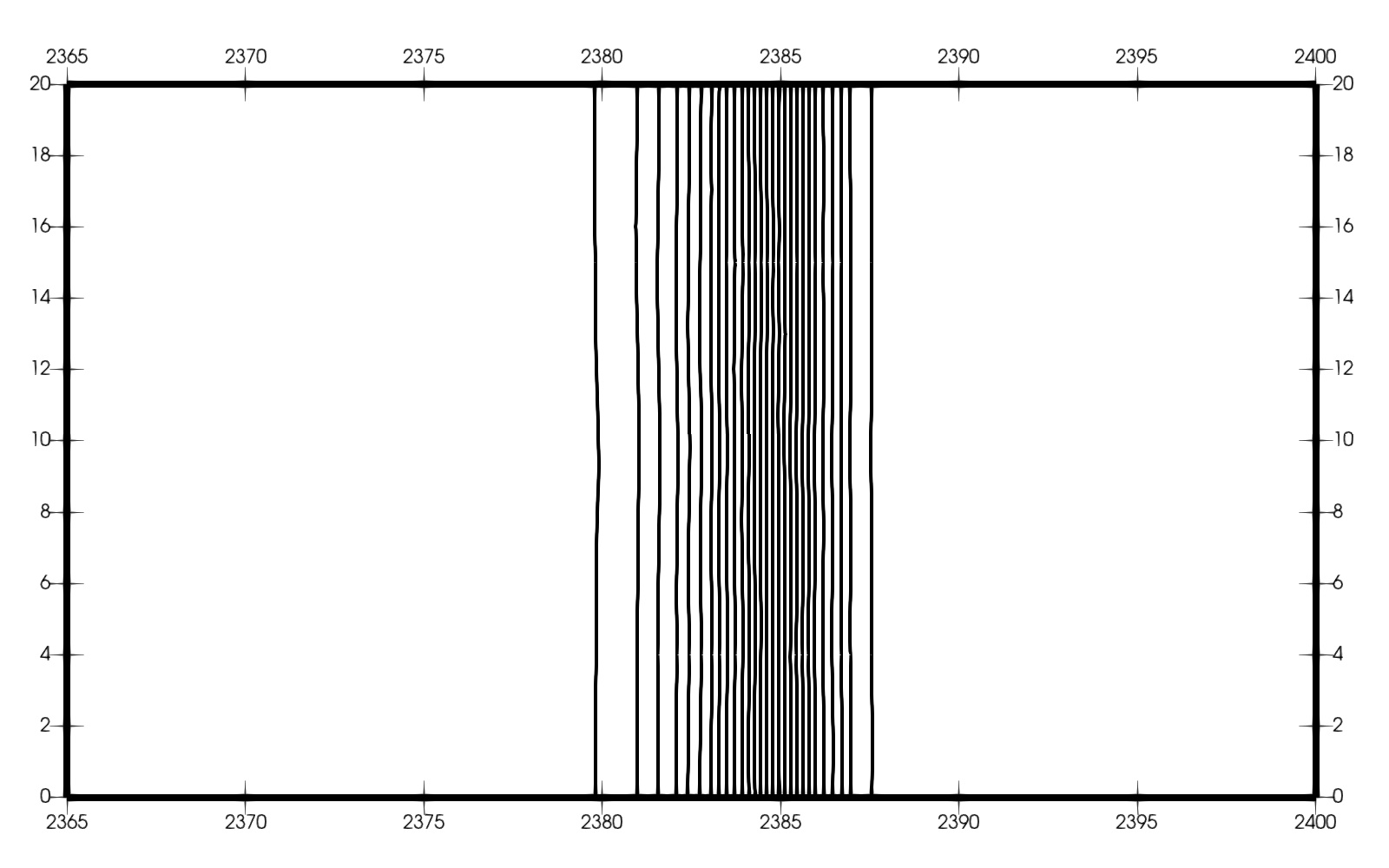}\quad
			\caption{}
			\label{odd_even_1}
		\end{subfigure}
			\begin{subfigure}{.3\textwidth}
			\centering
			\includegraphics[width=\linewidth]{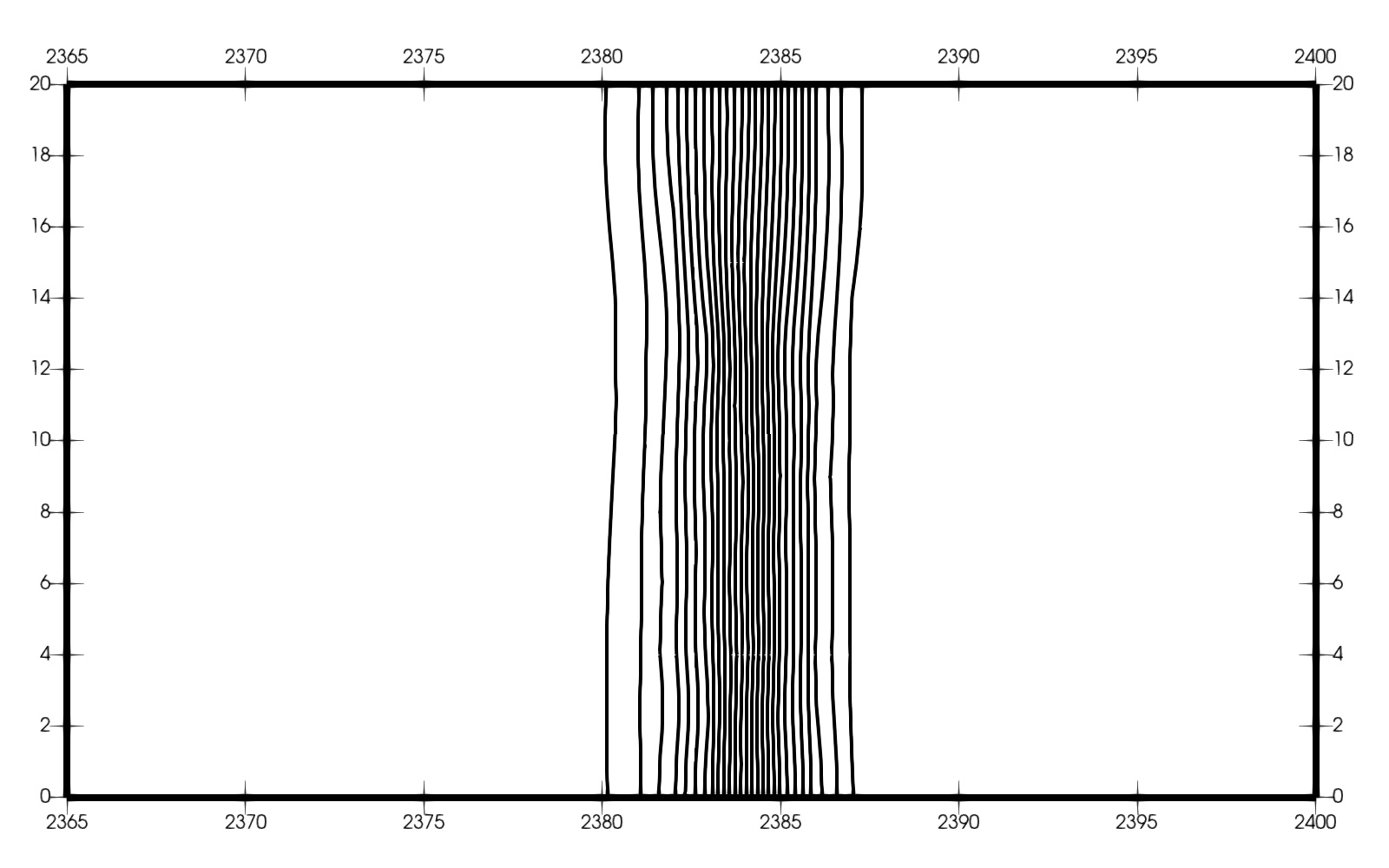}\quad
			\caption{}
			\label{odd_even_2}
		\end{subfigure}
		\caption{(a): Part of the grid (x=0 to 35) used for odd-even decoupling with centerline perturbations. (b): Contours of density (30 in 0.9-44) for odd-even decoupling using the ES scheme. (c): Contours of density (30 in 0.9-44) for odd-even decoupling using HES scheme}
		\label{odd_even}
	\end{figure}

	Another simulation anomaly is the grid-aligned planar shock structure deterioration caused by odd-even decoupling.  Given in \cite{quirk_1994}, in this test case, a planar shock travels in a long rectangular tube. The centre-line of the grid is perturbed slightly, given as 
	\begin{equation}
		y_{i,j}=
		\left\{ \begin{matrix}
			\begin{aligned} &y_{i,j}+0.1 \qquad \quad \text{for} \quad i \quad even,  \\ &y_{i,j}-0.1  \qquad \quad \text{for} \quad i \quad odd. \end{aligned}
		\end{matrix} \right.
	\end{equation}
Domain is taken as $[0,2400]\times[0,20]$ and is partly shown in figure \ref{odd_even_grid}. The domain is initialized with $\mathbf{W}=[1.0,0,0,1.0]$ for $x>5$, and post-shock conditions corresponding to Mach 20 are specified for $x<5$. It refers to a stronger shock than that given in \cite{quirk_1994}.  Flow tangency boundary conditions are applied at the top and bottom boundaries, and a supersonic outflow boundary condition was applied at the right boundary. The odd-evendecoupling phenomenon causes the shock structure to break down, and many schemes (like Riemann solvers) fail before reaching the final time of $T=330s$. This shock instability worsens with increasing the shock strength and magnitude of perturbations \cite{fleischmann_2020}. Results with ES and HES schemes are shown in figures \ref{odd_even_1} and \ref{odd_even_2}; with both the schemes, shock structure remains intact. This test case and the hypersonic flow over the half-cylinder (section 7.2.3) indicate that the schemes are free of any numerical shock instabilities. 
\subsubsection{Double Mach reflection}
This is another unsteady test case \cite{woodward_1984} where a shock of Mach 10 meets a wall at an angle of $30^o$, resulting in a complicated shock structure including a reflected shock and a slipstream. Domain for this test case is taken to be $[0,4]\times[0,1]$. The initial conditions for the problem are
	\begin{equation}
		\mathbf{W}=
		\left\{ \begin{matrix}
			\begin{aligned} & [1.4,0,0,1.0] \qquad \qquad \qquad \qquad \text{if} \quad y \leq \sqrt{3}(x-1/6)  \\ &[8.0,33\sqrt{3}/8,-4.125,116.5]  \quad \text{otherwise.}  \end{aligned}
		\end{matrix} \right.
	\end{equation}  

This is an unsteady benchmark test case introduced by Woodward \& Colella \cite{woodward_1984} to evaluate shock-capturing schemes. A Mach 10 shock, inclined at 60°, travels through air ($\gamma$ = 1.4) and reflects off a solid wall that begins at $x = 1/6$ along the bottom boundary. The shock is introduced not by a simple discontinuity, but through boundary-driven conditions: post-shock values are enforced on the bottom-left segment and the top boundary is dynamically updated to match the shock's angle and speed. This setup produces a complex, self-similar flow involving two triple points — each forming a Mach stem, a reflected shock, and a contact discontinuity (slipstream). One of the key challenges is resolving the high-density jet formed near the wall, which resembles shaped charge behavior and requires high-order methods like PPM to accurately capture. Lower-order methods typically fail to resolve this jet and the secondary reflected shock, highlighting the importance of accurate shock capturing and contact discontinuity resolution.

\begin{figure}
            \begin{subfigure}{.48\textwidth}
			\centering
			\includegraphics[width=\linewidth]{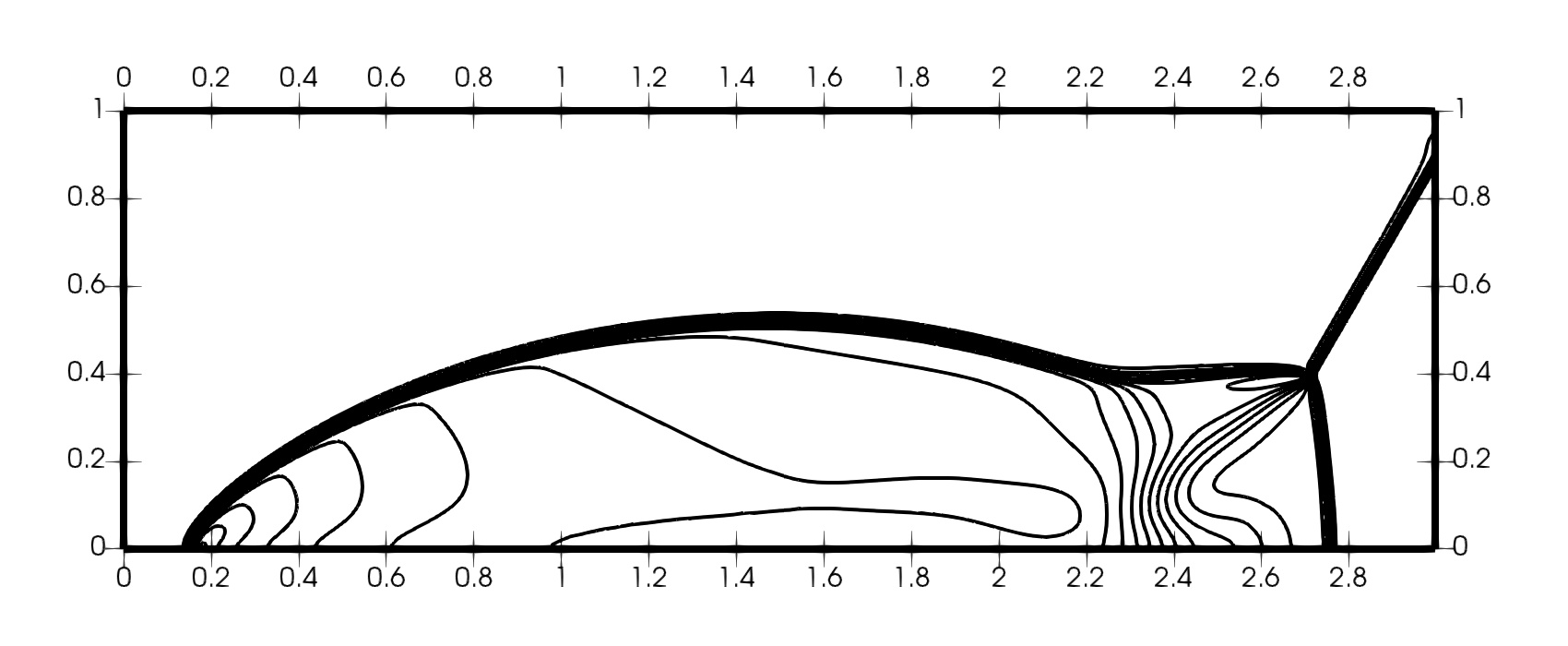}
            \end{subfigure}%
            \begin{subfigure}{.48\textwidth}
                \centering
			\includegraphics[width=\linewidth]{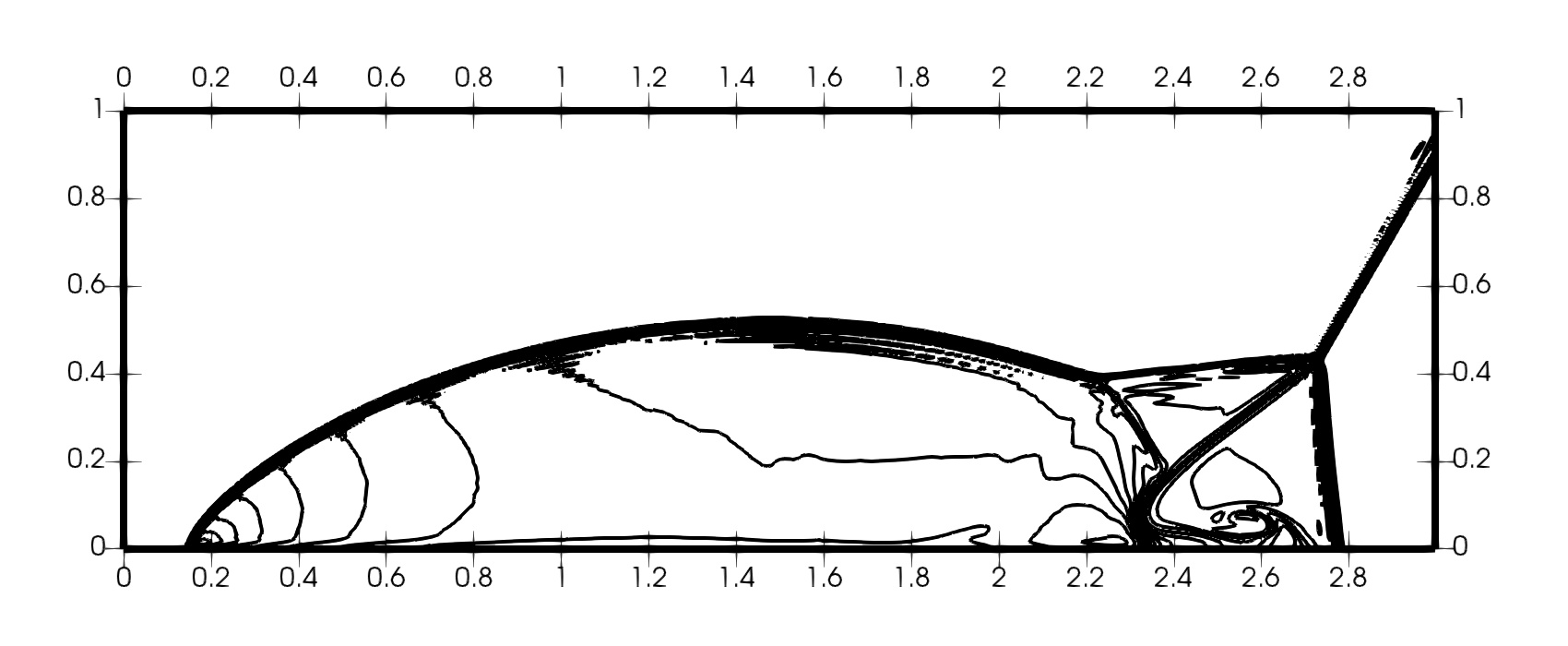}
            \end{subfigure}
            \caption{Density contours (30 in 1.3-22) for double mach reflection test case on a grid of 960x240 using ES (left) and HES (right) schemes. } 
            \label{mach_ref_1}
        \end{figure}

\begin{figure}
		\includegraphics[width=\linewidth]{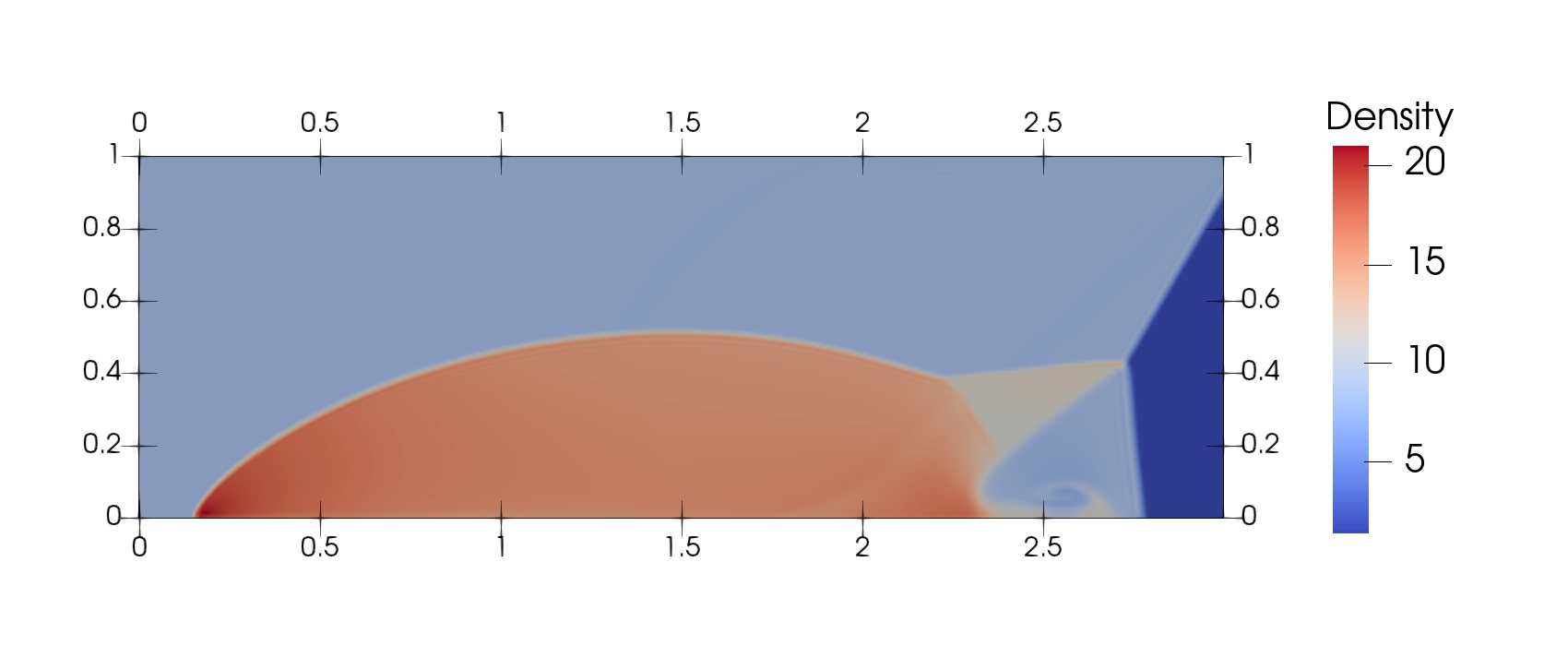}
            \caption{Color plot of density for the double mach reflection test case on a grid of 960x240 using the HES scheme.}
            \label{mach_ref_2}
	\end{figure}

Density contours of the solutions obtained using ES and HES schemes are given in figure \ref{mach_ref_1}. Prominent flow features like the Mach reflections and stems are captured well by both schemes.  However, only the higher order scheme is able to capture the complicated flow features such as the secondary reflected shock and jet with sufficient detail as shown in the density map in figure \ref{mach_ref_2}.    
 
\subsubsection{Shock - vortex filament interaction}  
\begin{figure}
		\begin{subfigure}{.32\textwidth}
			\includegraphics[width=\linewidth]{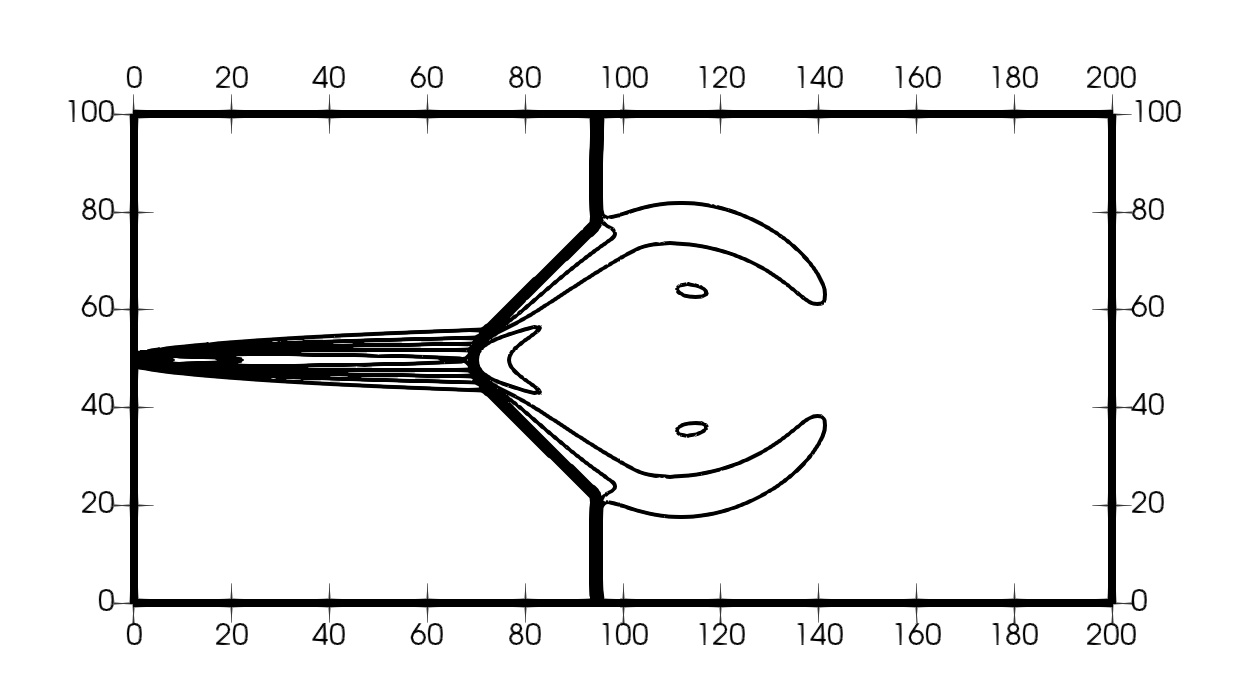}
			\caption{}
		\end{subfigure}
		\begin{subfigure}{.32\textwidth}
			\includegraphics[width=\linewidth]{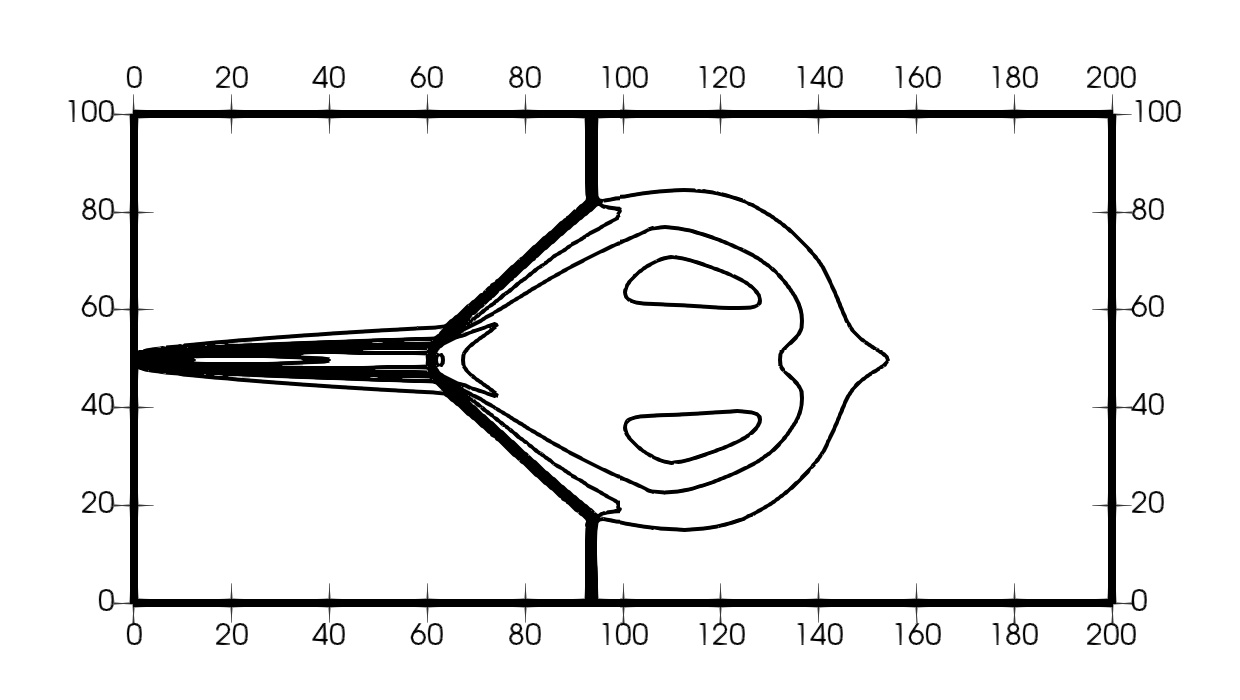}
			\caption{}
		\end{subfigure}
		\begin{subfigure}{.32\textwidth}
			\includegraphics[width=\linewidth]{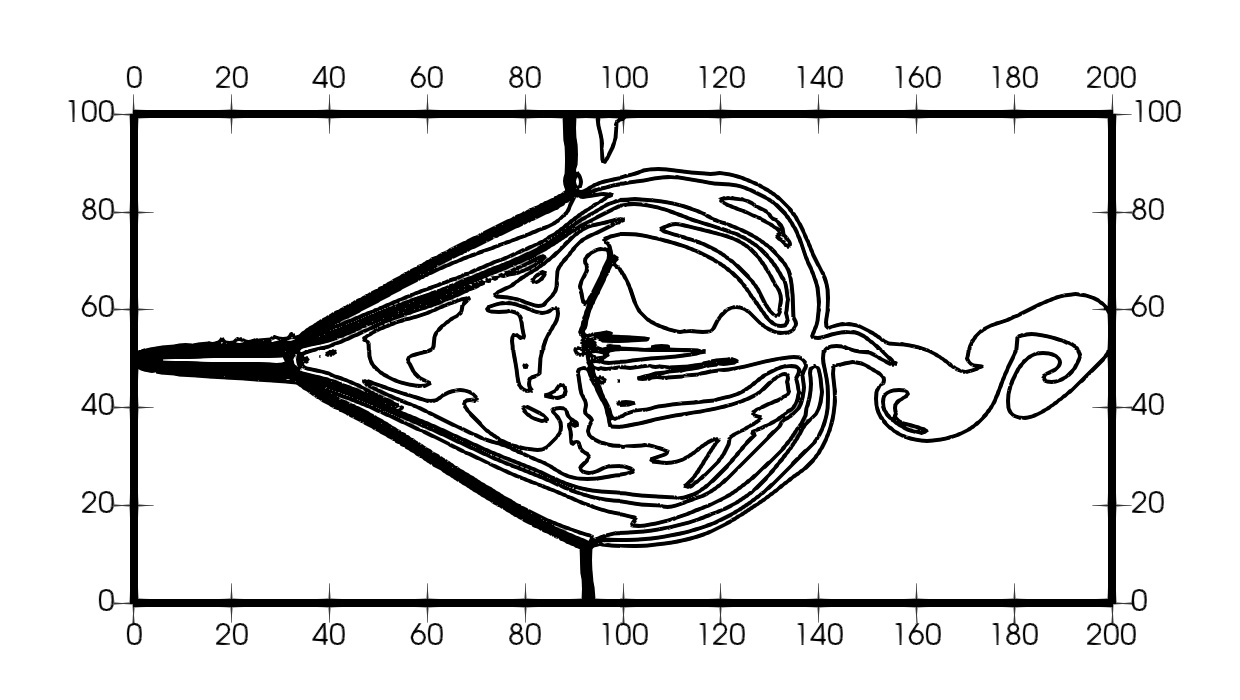}
			\caption{}
		\end{subfigure}
		\caption{entropy $(s=ln(p/\rho^\gamma))$ contours (15 in -0.6-4.2) for shock vortex filament interaction with (a) LLF (b) ES (c) HES schemes}
		\label{shock_vorfil_int}
	\end{figure}

In section \ref{half-cylinder}, the new schemes are demonstrated not producing any unphysical carbuncle shock, which is ultimately a numerical artifact. Recent studies have demonstrated that carbuncles can be triggered under certain physical conditions. Elling \cite{elling_2009} shows that a vortex filament interacting with a steady shock of high strengths will trigger a carbuncle structure that grows with time. This is true for many schemes, such as the Gudunov, Lax-Friedrichs, and  Osher schemes. Kemm \cite{kemm_2018} shows that schemes such as HLLE have very high numerical shear viscosity, preventing the formation of the physical carbuncles induced by vortex filament. Since this test case is a model for the shock-boundary layer interaction problem (vortex filament representing a numerical boundary layer), correct capturing of the induced carbuncle is essential and must not be suppressed. The domain for this problem is taken to be [0,200]x[0,100] with a steady shock of M=20 at x=100. At the left boundary, supersonic inflow conditions are imposed corresponding to $\mathbf{W}={1.0,20.0,0,1.0/1.4}$ everywhere except at a single cell in the centre where $\mathbf{W}={1.0,0,0,1.0/1.4}$ is imposed. The top and bottom boundaries are taken as solid walls, and the right boundary is considered to represent supersonic outflow. Simulation is run for a time of $T=20s$. As shown in figure \ref{shock_vorfil_int}, the physical carbuncle is captured for both schemes, and the flow features are much more pronounced for the HES scheme than for the ES scheme. Compared with the LLF scheme in figure \ref{shock_vorfil_int} (a), the carbuncle is captured more precisely with ES and HES schemes.  
\subsubsection{Flows NACA 0012 airfoil}
\begin{figure}
		\centering
		\begin{minipage}{.8\textwidth}
			\includegraphics[width=0.48\linewidth]{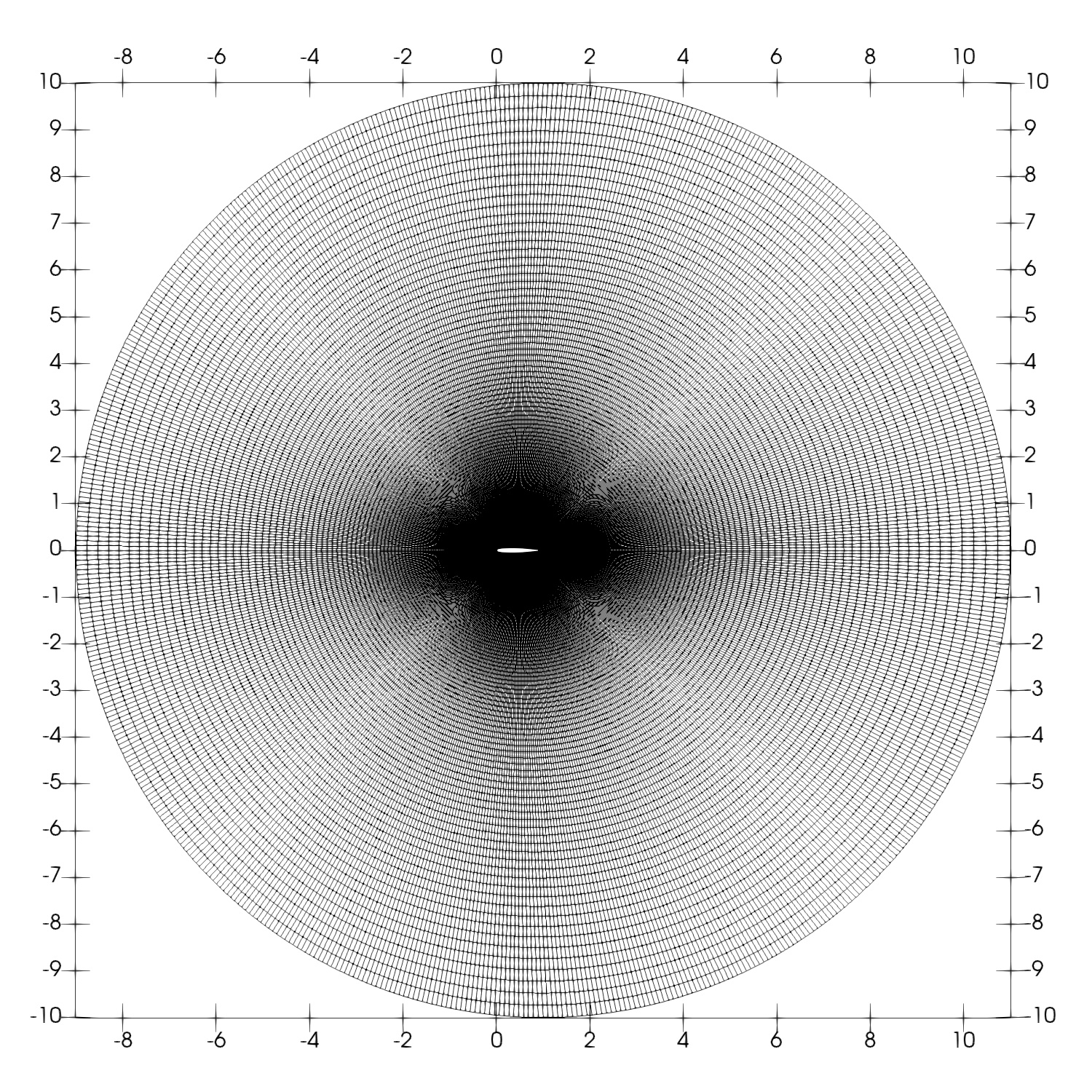}\quad
			\includegraphics[width=0.48\linewidth]{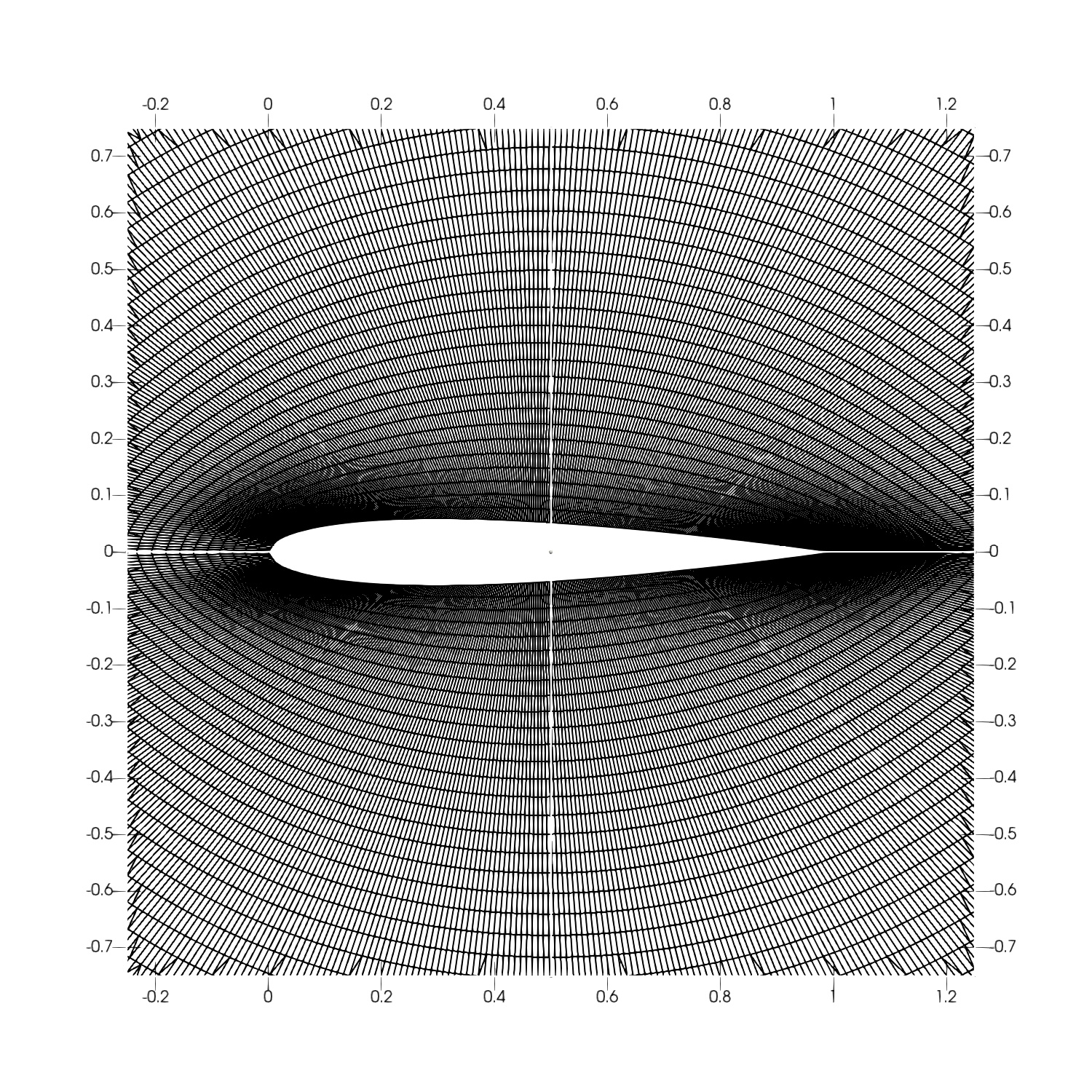}\quad
		\end{minipage}
		\caption{Structured mesh for NACA0012 airfoil}
		\label{airfoil_grid}
	\end{figure}

\begin{figure}
		\centering
		\captionsetup{justification=centering}
		\begin{subfigure}{.48\textwidth}
			\centering
			\includegraphics[width=\linewidth]{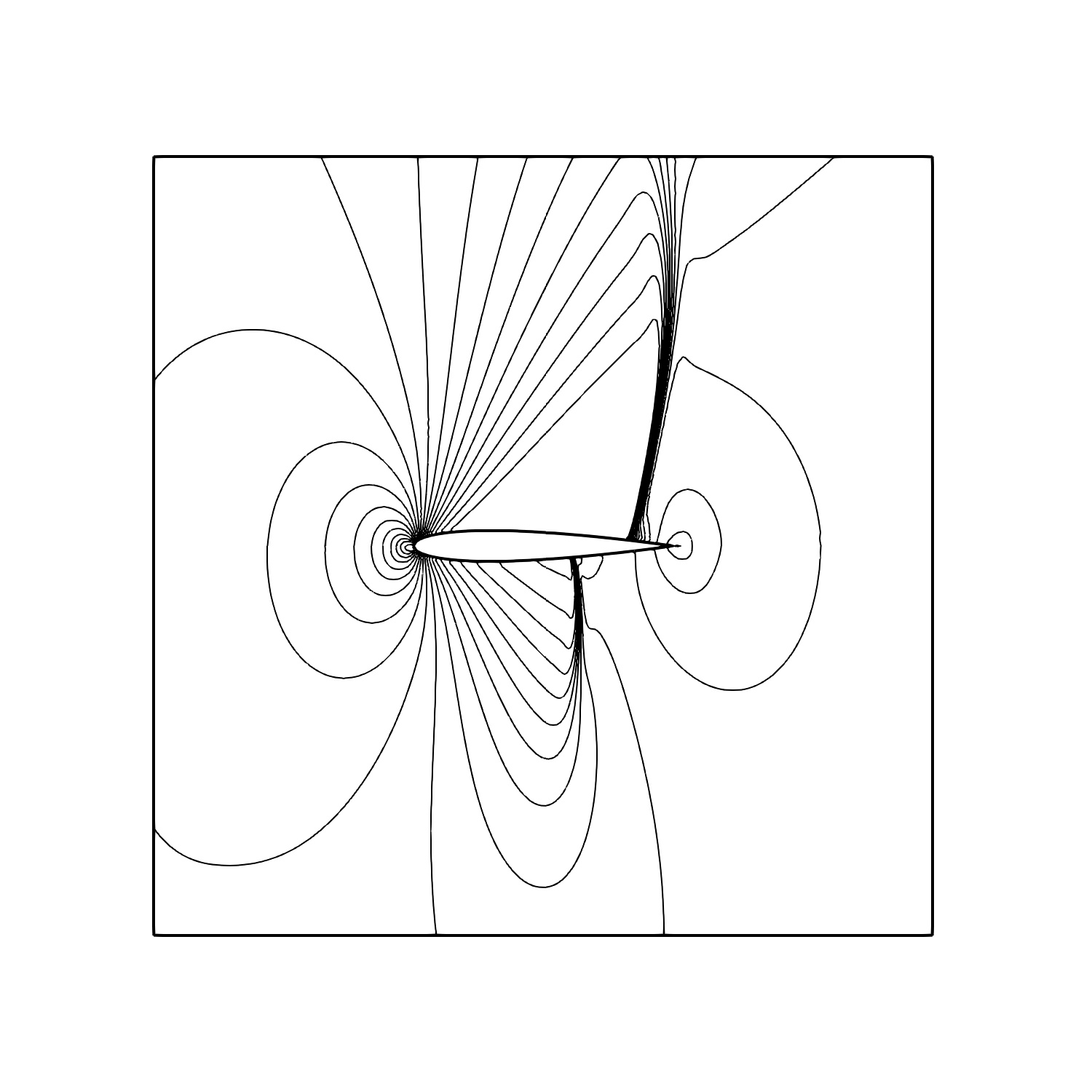}
			\caption{Transonic case with M=0.85 and $2^o$ angle of attack. Pressure contours (30 in 0.4-1.6) using HES Scheme.}
		\end{subfigure}
		\begin{subfigure}{.48\textwidth}
			\centering
			\includegraphics[width=\linewidth]{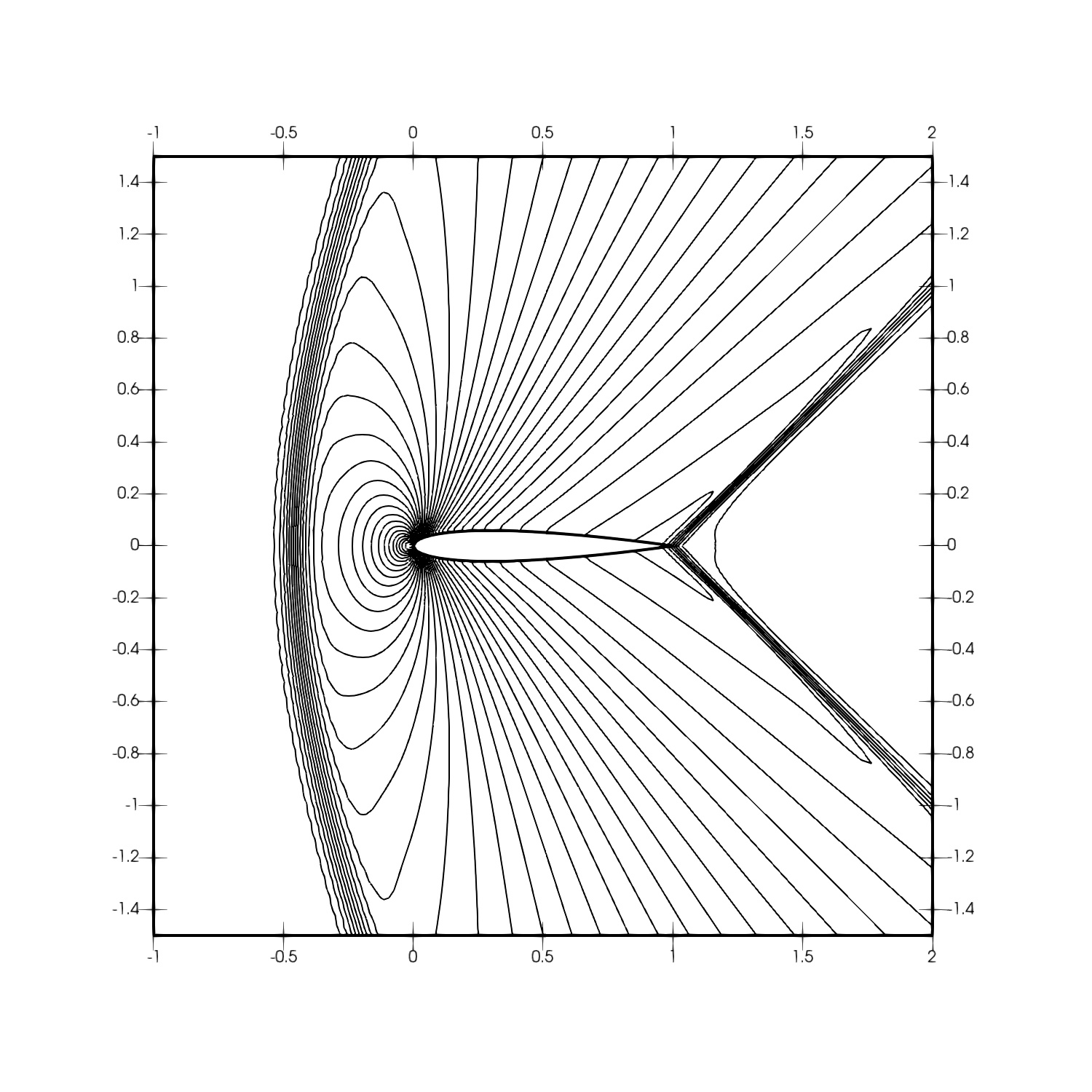}
			\caption{Transonic case with M=1.2 and $0^o$ angle of attack. Pressure contours (30 in 0.5-2.4) using HES Scheme.}
		\end{subfigure}
		\caption{NACA 0012 Airfoil test case}
		\label{airfoil_countour}
	\end{figure}

\begin{figure}
		\begin{subfigure}{0.24\textwidth}
			\includegraphics[width=\linewidth]{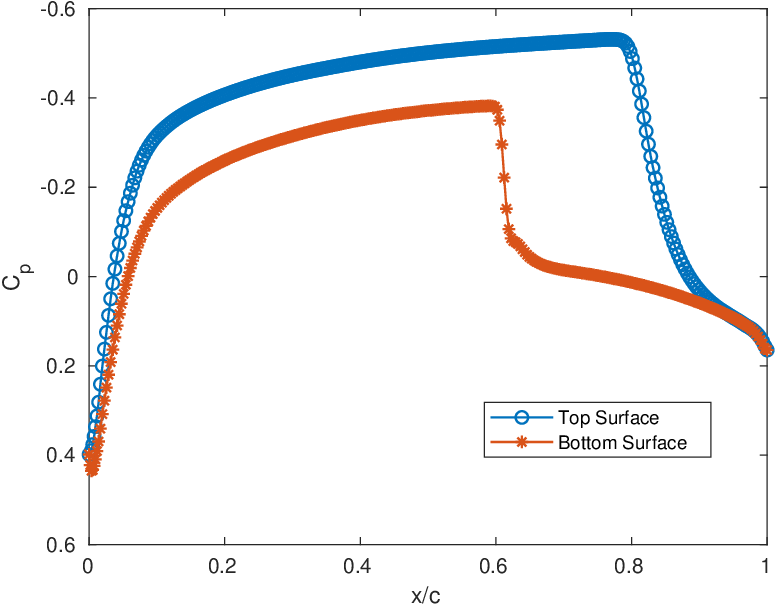}
			\caption{}
		\end{subfigure}
		\begin{subfigure}{0.24\textwidth}
			\includegraphics[width=\linewidth]{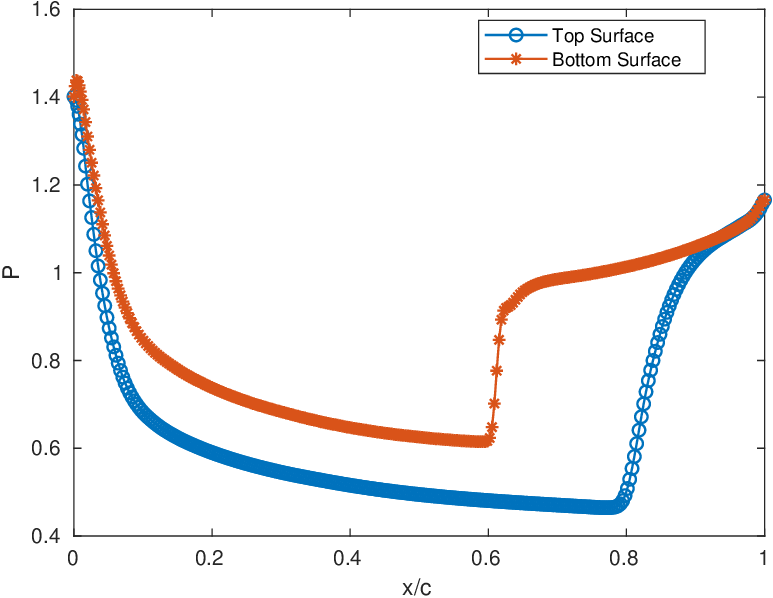}
			\caption{}
		\end{subfigure}
		\begin{subfigure}{0.24\textwidth}
			\includegraphics[width=\linewidth]{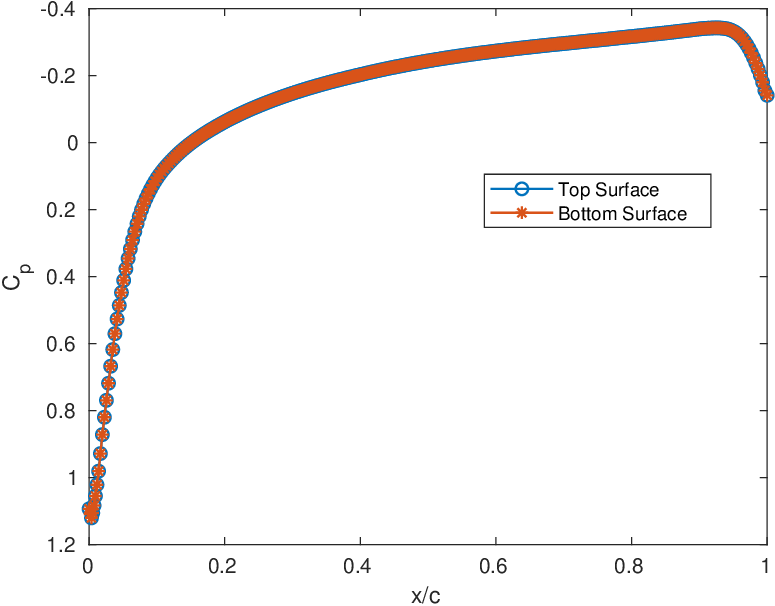}
			\caption{}
		\end{subfigure}
		\begin{subfigure}{0.24\textwidth}
			\includegraphics[width=\linewidth]{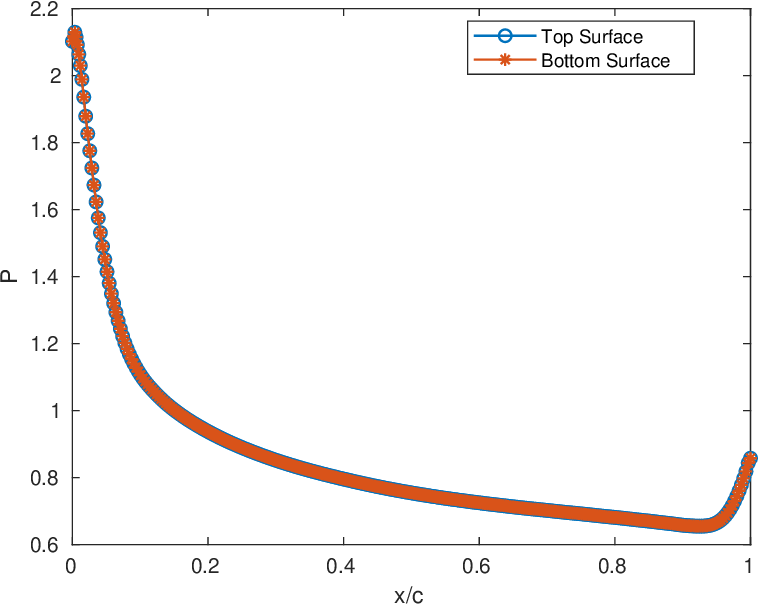}
			\caption{}
		\end{subfigure}
		\caption{Pressure and $C_p$ plots along top and bottom surfaces of airfoil using HES scheme : (a) and (b) Transonic case. (c) and (d) Supersonic case.}
		\label{airfoil_plots}
	\end{figure}

The new schemes are tested for flows over a symmetric NACA 0012 airfoil, in trans-sonic and supersonic regimes. The domain is taken as a circle of radius of 10 units centered at the trailing edge of an airfoil with chord of unit length. A structured mesh of 200x300 elements with 200 points on the airfoil surface is used for the simulation (see fig \ref{airfoil_grid}). Grid is stretched in radial direction by a 2.5\%. Far-field boundary conditions are imposed on the outer boundary, and flow tangency wall boundary conditions are imposed on the airfoil surface. Flow field is initialized with $\mathbf{W}={1.4,M \cos(\theta),M\sin(\theta),1.0/1.4}$, where M is the mach number of flow and $\theta$ is the angle of attack.  Flow is computed for the following conditions: M=0.85, $\theta=2^o$ and M=1.2, $\theta=0$. These computations are only performed for the hybrid entropy stable flux \eqref{es_final}. Pressure and pressure coefficient $C_p=(p-p_{\infty})/(\frac{1}{2} \rho_{\infty} u_{\infty}^2)$ for both test cases are shown in figure \ref{airfoil_plots}. For the transonic case, the shocks on upper and lower surfaces can be seen to be captured accurately in these plots.  For the supersonic csse, crisply captured bow shock and fish-tail shocks can be observed. Pressure contours for both schemes are shown in figure \ref{airfoil_countour} and the $C_P$ plots are shown in figure \ref{airfoil_plots}.   
		\label{mach_ref}


    \newpage 

    \section{Summary and Conclusions}
A family of structure-preserving numerical fluxes for the Euler equations that conserve entropy and preserve kinetic energy in a semi-discrete sense are introduced in this work. Beginning with a derivation based on the inherent structures within the energy equation, new entropy-conservative and kinetic energy-preserving fluxes that do not require logarithmic averages are proposed, thereby improving computational efficiency. These fluxes maintain second-order accuracy and are validated on representative smooth problems such as the Taylor-Green vortex and isentropic vortex convection, demonstrating their effectiveness in conserving physical invariants in the absence of shocks.

To extend applicability to problems involving discontinuities, a scalar numerical diffusion mechanism based on the Rankine–Hugoniot conditions is introduced. This diffusion strategy ensures entropy stability and kinetic energy damping while preserving steady contact discontinuities exactly. Furthermore, to reduce numerical diffusion in smooth regions, a hybrid entropy-stable (HES) scheme is developed using an entropy-distance-based shock sensor to selectively switch between low and high diffusion modes. The resulting scheme balances accuracy and stability across both smooth and discontinuous regimes.

Extensive one- and two-dimensional numerical experiments confirm the robustness, accuracy, and stability of the proposed methods. Notably, the schemes successfully handle classical and challenging test cases such as odd-even decoupling, double Mach reflection, shock-vortex interaction, and supersonic flow over airfoils, without exhibiting common numerical artifacts like expansion shocks or carbuncle-like instabilities.   
\\
\\
{\bf{CRediT author statement}} \\
{\bf{Kunal Bahuguna}}: Conceptualization, Methodology, Software, Validation, Formal analysis, Writing — original draft, Visualization. \\
{\bf{Ramesh Kolluru}}: Conceptualization, Validation, Investigation, Writing — review \& editing, Supervision. \\
{\bf{S. V. Raghurama Rao}}: Conceptualization, Validation, Investigation, Resources, Writing — review \& editing, Supervision. \\
\\
{\bf{Declaration of competing interest}}\\
The authors declare that they have no known financial interests or personal relationships with any other people or organisations that could influence the work presented here.\\

\section{Appendices}

\subsection{Appendix A: Entropy conservative flux using an optimisation procedure}
\renewcommand{\theequation}{A-\arabic{equation}}
\setcounter{equation}{0}
    Entropy-conservative fluxes given in section 3 can also be obtained from an optimisation problem defined as follows. 
    \begin{equation}
        \min || \mathbf{F^c}-\overline{\mathbf{F}} || \quad \text{such that} \quad \Delta \mathbf{V} \cdot \mathbf{F}^c = \Delta \psi
    \end{equation}
    The above quadratic optimisation with an equality constraint can be analytically solved using the method of Lagrange multipliers. The Lagrange function can be given as
    \begin{equation}
        L(\mathbf{F}^c,\lambda)=|| \mathbf{F^c}-\overline{\mathbf{F}} ||+\lambda (\mathbf{V} \cdot \mathbf{F}^c - \Delta \psi)
    \end{equation}
    where $\lambda$ is the Lagrange multiplier corresponding to the entropy conservative constraint. The solution to this optimisation problem can be found using equation 
    \begin{equation}
        \nabla \cdot L(\mathbf{F}^c,\lambda) = 0
    \end{equation}
    which gives
    \begin{subequations}
        \begin{align}
            \mathbf{F}^c=&\overline{\mathbf{F}}-\frac{\lambda}{2} \Delta \mathbf{V} \\
            \frac{\lambda}{2}=&\frac{\Delta \mathbf{V} \cdot \overline{\mathbf{F}} - \Delta \psi}{\Delta \mathbf{V} \cdot \Delta \mathbf{V}}
        \end{align}
    \end{subequations}
    Note that the $\lambda/2$ here is exactly the dissipation term in EC1 flux \eqref{EC1_flux}.
    A similar procedure can be followed to obtain an entropy conservative and kinetic energy preserving flux. The optimisation problem can be formulated with two constraints as follows.
    \begin{equation}
        \min || \mathbf{F^c}-\overline{\mathbf{F}} || \quad \text{such that} \quad \Delta \mathbf{V} \cdot \mathbf{F}^c = \Delta \psi \quad \text{and} \quad F_2^c=F_1^c \overline{u}+\overline{p}
    \end{equation}
    This, however, results in a complicated flux. To obtain a simple flux we can modify the above problem as
    \begin{equation}
        \min || F_3^c-\overline{F_3} || \quad \text{such that} \quad \Delta \mathbf{V} \cdot \mathbf{F}^c = \Delta \psi
    \end{equation}
    where we assume $F_1^c=\overline{F_1}$ and $F_2^c=F_1^c \overline{u}+\overline{p}$. Lagrange function for this can be given as
    \begin{equation}
        L({F_3}^c,\lambda)=|| F_3^c-\overline{F_3} ||+\lambda (\mathbf{V} \cdot \mathbf{F}^c - \Delta \psi)
    \end{equation}
    with the following solution obtained by taking $\nabla \cdot L(F_3^c,\lambda)=0$.
    \begin{subequations}
        \begin{align}
            F_3^c=&\overline{F_3}-\frac{\lambda}{2} \Delta V_3 \\
            \frac{\lambda}{2}=&\frac{\Delta \mathbf{V} \cdot \overline{\mathbf{F}} - \Delta \psi}{\Delta {V_3} \cdot \Delta {V_3}}
        \end{align}
    \end{subequations}
    This corresponds to the ECKEP flux \eqref{ECKEP_flux} obtained in section 3.
\subsection{Appendix B: Entropy Distance Positivity}
\renewcommand{\theequation}{B-\arabic{equation}}
\setcounter{equation}{0}
    For small changes we can write entropy distance $\Delta \mathbf{V} \cdot \Delta \mathbf{U}$ as
    \begin{equation}
        ED=d\mathbf{V}^T d\mathbf{U}=\left(\frac{d\mathbf{V}}{d\mathbf{U}} \cdot d\mathbf{U}\right)^T  d \mathbf{U}=d \mathbf{U}^T \left(\frac{d \mathbf{V}}{d \mathbf{U}}\right)^T d \mathbf{U}
    \end{equation}
    We know that $\mathbf{V}=\frac{d \eta(\mathbf{U})}{d\mathbf{U} }$ and $\frac{ d \mathbf{V}}{ d \mathbf{U}}=\frac{d^2 \eta(\mathbf{U})}{d\mathbf{U}^2}$  which is a positive definite matrix because of convexity of $\eta(\mathbf{U})$ and hence $ED \geq 0$ for any states $\mathbf{U}_1$ and $\mathbf{U}_2$.
	\bibliographystyle{elsarticle-num}
	\bibliography{EC_Paper_Citations}
\end{document}